\documentclass[11pt,leqno]{article}

\usepackage{amsmath,amsfonts,amscd,amssymb,theorem}

\long\def\comment#1\endcomment{}

\comment
\endcomment

\makeatletter
\begingroup
\gdef\th@dotted{\normalfont\itshape
  \def\@begintheorem##1##2{%
        \item[\hskip\labelsep \theorem@headerfont ##1\ ##2.]}%
\def\@opargbegintheorem##1##2##3{%
   \item[\hskip\labelsep \theorem@headerfont ##1\ ##2\ (##3).]}}
\endgroup
\makeatother

\theoremstyle{dotted}

\newtheorem{theorem}{Theorem}[section]
\newtheorem{lemma}[theorem]{Lemma}

\newtheorem{prop}[theorem]{Proposition}
\newtheorem{corr}[theorem]{Corollary}

\makeatletter
\begingroup
\gdef\th@upshape{\normalfont
  \def\@begintheorem##1##2{%
        \item[\hskip\labelsep \theorem@headerfont ##1\ ##2.]}%
\def\@opargbegintheorem##1##2##3{%
   \item[\hskip\labelsep \theorem@headerfont ##1\ ##2\ (##3).]}}
\endgroup
\makeatother

\theoremstyle{upshape}

\newtheorem{defn}[theorem]{Definition}
\newtheorem{remark}[theorem]{Remark}
\newtheorem{exa}[theorem]{Example}

\makeatletter
\renewcommand{\subsection}{\@startsection{subsection}{2}{0pt}{-3ex
plus -1ex minus -0.2ex}{-2mm plus -0pt minus
-2pt}{\normalfont\bfseries}} 
\renewcommand{\subsubsection}{\@startsection{subsubsection}{3}{0pt}{-3ex
plus -1ex minus -0.2ex}{-2mm plus -0pt minus
-2pt}{\normalfont\bfseries}} 
\makeatother

\makeatletter
\@addtoreset{equation}{section}
\makeatother

\newcommand{\cntrct}                
{\hspace{2pt}\raisebox{1pt}{\text{$\lrcorner$}}\hspace{2pt}}

\newcommand{\proof}[1][Proof.]{\smallskip\noindent{\em #1}}
\def\endproof{\hfill\ensuremath{\square}\par\medskip}

\def\eqref#1{\thetag{\ref{#1}}}

\let\latexref=\ref
\def\ref#1{{\normalfont{\latexref{#1}}}}

\newcommand{\wt}{\widetilde}
\newcommand{\wh}{\widehat}

\newcommand{\idot}{{\:\raisebox{1pt}{\text{\circle*{1.5}}}}}
\newcommand{\hdot}{{\:\raisebox{3pt}{\text{\circle*{1.5}}}}}
\newcommand{\Z}{{\mathbb Z}}
\newcommand{\Qq}{{\mathbb Q}}
\newcommand{\N}{{\mathbb N}}
\newcommand{\RR}{{\mathbb R}}
\newcommand{\CC}{{\mathbb C}}

\newcommand{\eps}{\varepsilon}
\renewcommand{\phi}{\varphi}

\newcommand{\Rr}{{\sf R}}
\newcommand{\Tt}{{\sf T}}

\def\dlim_#1{\displaystyle\lim_{#1}}

\newcommand{\hash}{\sharp}

\newcommand{\Hom}{\operatorname{Hom}}
\newcommand{\Ext}{\operatorname{Ext}}
\newcommand{\RHom}{\operatorname{RHom}}

\newcommand{\Fun}{\operatorname{Fun}}

\newcommand{\id}{\operatorname{\sf id}}
\newcommand{\Id}{\operatorname{\sf Id}}
\newcommand{\gr}{\operatorname{\sf gr}}

\newcommand{\A}{{\cal A}}
\newcommand{\D}{{\cal D}}
\newcommand{\M}{{\cal M}}
\newcommand{\C}{{\cal C}}
\newcommand{\B}{{\cal B}}

\newcommand{\Q}{{\cal Q}}
\newcommand{\T}{{\cal T}}

\newcommand{\R}{{\cal R}}

\newcommand{\I}{{\mathbb I}}

\newcommand{\Cycl}{\operatorname{Cycl}}
\newcommand{\Sets}{\operatorname{Sets}}

\newcommand{\Maps}{\operatorname{Maps}}

\newcommand{\Aut}{{\operatorname{Aut}}}

\newcommand{\Sthom}{\operatorname{\sf StHom}}

\newcommand{\amod}{{\text{\rm -mod}}}
\newcommand{\ppt}{{\sf pt}}

\newcommand{\Ab}{\operatorname{Ab}}

\newcommand{\copr}{{\textstyle\coprod}}

\newcommand{\LZ}{\Lambda\Z}
\newcommand{\wLZ}{\wh{\LZ}}
\newcommand{\LR}{\Lambda R}
\newcommand{\wLR}{\wt{\Lambda R}}
\newcommand{\LI}{\Lambda I}

\newcommand{\DeR}{\Delta R}
\newcommand{\wDer}{\wt{\Delta R}}

\newcommand{\DML}{\D\!\M\!\Lambda}
\newcommand{\DM}{{\cal D}{\cal M}}

\newcommand{\DLR}{\D\!\Lambda\text{{\upshape R}}}

\newcommand{\DF}{\D{\mathcal F}}

\newcommand{\FDM}{{\mathcal F}\D\!\M}

\newcommand{\wq}{\widehat{q}}

\newcommand{\bO}{\overline{O}}

\newcommand{\Ass}{\operatorname{\sf Ass}}

\newcommand{\Exp}{\operatorname{\sf Exp}}
\newcommand{\Div}{\operatorname{\sf Div}}
\newcommand{\Stab}{\operatorname{\sf Stab}}

\newcommand{\Supp}{\operatorname{\sf Supp}}

\newcommand{\holim}{\operatorname{\sf holim}}

\newcommand{\Sp}{\text{\rm-Sp}}
\newcommand{\spp}{\text{\rm-sp}}
\newcommand{\Top}{\text{\rm-Top}}

\newcommand{\Real}{\operatorname{\rm Real}}

\newcommand{\wPhi}{\widehat{\Phi}}
\newcommand{\wPsi}{\widehat{\Psi}}

\newcommand{\TC}{\operatorname{TC}}
\newcommand{\can}{\operatorname{\sf can}}

\newcommand{\HH}{\mathcal{H}}

\title{Cyclotomic complexes}

\author{D. Kaledin}

\begin{document}

\maketitle

\tableofcontents

\section*{Introduction.}

This paper is a sequel to \cite{Ka-ma}; the goal of both papers is
to try to understand how some purely homological notions used in
\cite{K} are related to Topological Cyclic Homology of
\cite{BHM}. This turned out to be a rather lengthy project, since
one has to construct appropriate homological counterparts of several notions
from stable homotopy theory. 

The paper \cite{Ka-ma} dealt with derived Mackey functors,
homological analogues of genuine $G$-equivariant spectra of
\cite{LMS}. While the abelian category of Mackey functors is very
well known in stable homotopy theory, and does play an important
role, its most naive derived generalization turned to be not quite
well-behaved. Thus a slightly different derived version of Mackey
functors was constructed and studied in \cite{Ka-ma}.

The present paper deals with cyclotomic spectra of \cite{BM},
especially as presented in \cite{HM}. To find a homological
counterpart for those, one first has to go beyond \cite{Ka-ma}: by
definition, cyclotomic spectra are equivariant with respect to the
circle group $S^1$, and \cite{Ka-ma} only dealt with finite groups
$G$. We cannot really construct a good homological analogue of all
$S^1$-spectra spectra, but we do construct a category $\DML(\Z)$ of
``cyclic Mackey functors'' which captures the part of the
equivariant stable category relevant to Topological Cyclic
Homology. We then introduces the triangulated category $\DLR(\Z)$ of
``cyclotomic complexes''.

Ideally, the relation between cyclotomic complexes and cyclotomic
spectra should be expressed by a commutative diagram
\begin{equation}\label{dgg}
\begin{CD}
\DLR(\Z) @>>> \Cycl\\
@VVV @VVV\\
\D(\Z) @>>> \Sthom
\end{CD}
\end{equation}
of ``brave new schemes'', understood for example as tensor
triangulated categories with some enhancement. Here $\Sthom$ is the
stable homotopy category, $\D(\Z)$ is the derived category of
abelian groups, and $\Cycl$ is the category of cyclotomic
spectra. The diagram should be ``almost Cartesian''. More precisely,
it should become Cartesian if we restrict our attension to the
subcategory $\Cycl^o \subset \Cycl$ of cyclotomic spectra $T$ with
trivial geometric fixed points $\Phi^{S^1}T$ with respect to the
whole group $S^1$.

At present, such a nice picture seems way beyond reach; besides the
obvious difficulties with making all the ``brave new'' notions
precise, it seems that up to now, no-one constructed $\Cycl$ even as
a triangulated category. Thus in practice, we restrict our attention
to the following two things:
\begin{enumerate}
\item we construct an equivariant homology cyclotomic complex
  $C_\idot(T)$ for every cyclotomic spectrum $T$; this ought to
  correspond to the top arrow in \eqref{dgg},
\item we construct a topological cyclic homology functor $\TC$ on
  the category $\DLR(\Z)$ in such a way that for any cyclotomic
  spectrum $T$, $\TC(C_\idot(T))$ is naturally identified with the
  homology of the spectrum $\TC(T)$.
\end{enumerate}
We then adopt a different perspective and give a completely
different and very simple description of the category $\DLR(\Z)$.
As it happens, cyclotomic complexes are essentially equivalent to
``filtered Dieudonn\'e modules'' of \cite{FL}. Filtered Dieudonn\'e
modules are rather simple linear-algebraic gadgets with a deep
meaning --- they give a $p$-adic counterpart of Deligne's notion of
a mixed Hodge structure, and the whole story acquires a distinctly
motivic flavour. A more detailed discussion of this is available in
\cite{icm}.

Filtered Dieudonn\'e modules arise naturally as the cristalline
cohomology of algebraic varieties over $\Z_p$, while cyclotomic
spectra appear as Topological Hochschild Homology spectra of ring
spectra $A$. In view of the equivalence we established, there are
many areas of intersection where one can compare the two
constructions. We did not attempt to do so in this paper; we stick
to pure linear algebra, and leave the geometric applications for
future research. The only comparison result that we prove says that
for profinitely complete cyclotomic complexes, the topological
cyclic homology $\TC$ in fact coincides with the syntomic cohomology
of \cite{FM} familiar in the theory of Dieudonn\'e modules.

\medskip

The paper is organized as follows. To begin the story, we need some
model for $S^1$-equivariant spaces and their homology; for better or
for worse, we have chosen to use the combinatorial approach using
A. Connes's category $\Lambda$. Section 1 contains necessary fact
about the category $\Lambda$ and its various cousins. In Section 2,
we construct cyclic Mackey functors. Section 3 deals with cyclotomic
complexes. Section 4 is the tecnical heart of the paper: here we
construct equivariant homology functors from $S^1$-spectra to
$\DML(\Z)$ and from cyclotomic spectra to $\DLR(\Z)$. Then in
Section 5, we forget all about topology and prove the comparison
theorem between cyclotomic complexes and filtered Dieudonn\'e
modules. Finally, in Section 6 we briefly discuss topological cyclic
homology, and prove the comparison theorems for $\TC$. Appendix
contains some technicalities, mostly from \cite{Ka-ma}.

\subsection*{Acknowledgements.} Many discussions on the subject with
G. Merzon were very inspirational and helpful. I owe a lot to
L. Hesselholt for his patient explanation about topology (this goes
both for this paper in particular, and for the whole project in
general). I am grateful for V. Vologodsky for explanations about
Dieudonn\'e modules. It is a pleasure to thank A. Beilinson,
V. Drinfeld, V. Ginzburg, D. Kazhdan, and J.P. May for their
interest in this work. A part of the paper was finished at the
Hebrew University of Jerusalem, and another part was done while
visiting the University of Chicago; the hospitality of both places
is gratefully acknowledged.

\section{Cyclic categories.}

\subsection{Connes' cyclic category.}\label{connes.subs}

Recall that A. Connes' {\em cyclic category} $\Lambda$ is a small
category whose objects $[n]$ are indexed by positive integers $n$,
$n \geq 1$. Maps between $[n]$ and $[m]$ can be defined in various
equivalent ways; for the convenience of the reader, we recall two of
these descriptions.

\smallskip

\noindent
{\em Topological description}. The object $[n]$ is thought of as a
``wheel'' -- a cellular decomposition of the circle $S^1$ with $n$
$0$-cells, called {\em vertices}, and $n$ $1$-cells, called {\em
edges}. A continuous map $f:S^1 \to S^1$ induces a map $\wt{f}:\RR
\to \RR$ of the universal covers; say that $f$ is {\em monotonous}
if $\wt{f}(a) \geq \wt{f}(b)$ for any $a,b \in \RR$, $a \geq
b$. Then morphisms from $[n]$ to $[m]$ in the category $\Lambda$ are
homotopy classes of monotonous continuous maps $f:[n] \to [m]$ which
have degree $1$ and send vertices to vertices.

\smallskip

\noindent
{\em Combinatorial description}. Consider the category
$\Lambda_{big}$ of totally ordered sets equipped with an
order-preserving endomorphism $\tau$. Let $[n] \in \Lambda_{big}$ be
the set $\Z$ with the natural linear order and endomorphism $\tau:\Z
\to \Z$, $\tau(a) = a + n$. Let $\Lambda_\infty \subset
\Lambda_{big}$ be the full subcategory spanned by $[n]$, $n \geq
1$. For any $[n],[m] \in \Lambda_\infty$, the set
$\Lambda_\infty([n],[m])$ is acted upon by the endomorphism $\tau$
(on the left, or on the right, by definition it does not matter). We
define the set of maps $\Lambda([n],[m])$ in the category $\Lambda$
by
\begin{equation}\label{lambda.inf}
\Lambda([n],[m])=\Lambda_\infty([n],[m])/\tau.
\end{equation}
Other descriptions are possible, see \cite[Chapter 6]{Lo} and
Appendix to \cite{FT}.

\smallskip

For any $[n] \in \Lambda$, the set $V([n])$ of vertices of the
corresponding decomposition of the circle can be naturally
identified with the set $\Lambda([1],[n])$ of maps from $[1]$ to
$[n]$, and the set $E([n])$ of edges can be identified with the set
$\Lambda([n],[1])$ -- in particular, $E(-)$ is a contravariant
functor (geometrically, the preimage of an edge is contained in
exactly one edge). The automorphism group $\Aut([n])$ is the cyclic
group $\Z/n\Z$ generated by the clockwise rotation; we will denote
the generator by $\sigma$. In the combinatorial description,
$\sigma$ corresponds to the map $\Z \to \Z$, $a \mapsto a+1$.

Given an integer $p \geq 2$, one can define a category $\Lambda_p$
by taking the same set of objects $[n]$, $n \geq 1$, and setting
$$
\Lambda_p([n],[m])=\Lambda_\infty([n],[m])/\tau^p.
$$
The category $\Lambda_p$ is intermediate between $\Lambda_\infty$
and $\Lambda$; in particular, the obvious projections
$\Lambda_\infty([n],[m])/\tau^p \to \Lambda_\infty([n],[m])/\tau$
together define a functor
$$
\pi_p:\Lambda_p \to \Lambda.
$$
The functor $\pi_p$ is a bifibration with fiber $\ppt_p =
[1/(\Z/p\Z)]$, the groupoid with one object and automorphism group
$\Z/p\Z$. On the other hand, $\Lambda_p([n],[m])$ can be identified
with the set of maps $f:[np] \to [mp]$ in $\Lambda$ such that $f
\circ \sigma^n = \sigma^m \circ f$; this gives a canonical functor
$$
i_p:\Lambda_p \to \Lambda
$$
such that on objects, we have $i_p([n]) = [np]$. Denote by
\begin{equation}\label{li.def}
\LI = \coprod_{p \geq 1}\Lambda_p
\end{equation}
the disjoint union of all the categories $\Lambda_p$, $p \geq
1$. Then the functors $i_p$ and $\pi_p$ can be considered together
as two functors
\begin{equation}\label{i.p.def}
i,\pi:\LI \to \Lambda.
\end{equation}
The category $\Lambda$ is self-dual: an equivalence $\Lambda \cong
\Lambda^{opp}$ sends every object to itself, and a morphism $[n] \to
[m]$ represented by a map $f:\Z \to \Z$ goes to the map represented
by $f_\hash:\Z \to \Z$,
\begin{equation}\label{la.du}
f_\hash(a) = \max\{b \in \Z| f(b) \leq a \}.
\end{equation}
In the topological description, the duality interchanged edges and
vertices and corresponds to taking the dual cellular decomposition.

For any $[n],[m] \in \Lambda$, the set of maps $\Lambda([n],[m])$ is
finite. The groups $\Aut([n])$ and $\Aut([m])$ act on
$\Lambda([n],[m])$ by compositions, and both these actions are
stabilizer-free. We will need the following slightly more general
fact.

\begin{lemma}\label{twt.lemma}
Assume given three integers $m$, $n$, $l$ such that $m,l \geq 1$, $n
\geq 2$, and a map $f:[nl] \to [m]$ in $\Lambda$ such that
$$
f \circ \sigma^{l} = \sigma^{l_1} \circ f
$$
for some integer $l_1$, $0 \leq l_1 < m$. Then $m = nl_1$.
\end{lemma}

\proof{} Use the combinatorial description of $\Lambda$. Then $f$ is
represented by an order-preserving map $\wt{f}:\Z \to \Z$ such that
\begin{equation}\label{twt}
\wt{f}(a + nl) = \wt{f}(a) + m, \qquad \wt{f}(a + l) = \wt{f}(a) +
l_1 + bm
\end{equation}
for any $a \in \Z$, where $b$ is a fixed integer independent of
$a$. Since $a \leq a+l \leq a+nl$, this implies $0 \leq l_1 + bm
\leq m$, so that either $l_1 = 0$ and $b=1$, or $b=0$. The first
case is impossible since $\sigma^l$ acts on $\Lambda([nl],[m])$
without fixed points. Thus $b=0$, and \eqref{twt} immediately
implies the claim.
\endproof

The category $\Lambda/[1]$ of objects $[n] \in \Lambda$ equipped
with a map $[n] \to [1]$ is naturally equivalent to the category
$\Delta$ of non-empty finite totally ordered sets: geometrically,
$\Lambda/[1]$ is the category of wheels with a fixed edge, and
removing this edge creates a canonical clockwise total order on the
set of vertices of the wheel. We thus have a natural discrete
fibration $\Delta \cong \Lambda/[1] \to \Lambda$ inducing a
cofibration $j^o:\Delta^{opp} \to \Lambda^{opp}$. Dually, the
category $[1]\backslash\Lambda$ of objects $[n] \in \Lambda$
equipped with a map $[1] \to [n]$ is equivalent to $\Delta^{opp}$,
so that we get a natural discrete fibration $j:\Delta \to
\Lambda^{opp}$ (geometrically, $\Delta^{opp}$ is the category of
wheels with a fixed vertex). The same constructions work for the
categories $\Lambda_n$, $n \geq 2$. In particular, we obtain a
canonical functor
$$
j_n:\Delta \to \Lambda^{opp}_n,
$$
and we have $\pi \circ j_n = j$ for any $n \geq 2$. Let $\Delta_n
\to \Delta$ be the bifibration obtained by the Cartesian square
$$
\begin{CD}
\Delta_n @>>> \Lambda^{opp}_n\\
@VVV @VV{\pi_n}V\\
\Delta @>{j_n}>> \Lambda^{opp}.
\end{CD}
$$
Then the functor $j_n$ gives a splitting $\Delta \to \Delta_n$ of
this bifibration, so that we have $\Delta_n \cong \Delta \times
\ppt_n$. In more down-to-earth terms, this means that the group
$\Z/n\Z$ acts on the functor $j_n$. We can also compose $j_n$ with
the embedding $i_n$; this results in a commutative diagram
$$
\begin{CD}
\Delta @>{j_n}>> \Lambda^{opp}_n\\
@V{r_n}VV @VV{i_n}V\\
\Delta @>{j}>> \Lambda^{opp},
\end{CD}
$$
where $r_n:\Delta \to \Delta$ is the {\em edgewise subdivision
  functor} given by
\begin{equation}\label{edge}
r_n([m]) = [n] \times [m],
\end{equation}
where $[m]$, $[n]$ are totally ordered sets with $m$ resp. $n$
elements, and $[n] \times [m]$ is given the left-to-right
lexicographical order.

\subsection{Cyclotomic category.}\label{lr.subs}

We now introduce the following definition based on the topological
description of the category $\Lambda$.

\begin{defn}
  The {\em cyclotomic category} $\LR$ is the small category with the
  same objects $[n]$, $n \geq 1$, as the category $\Lambda$. If we
  think of $[n]$ as a configuration of $n$ marked points on a circle
  $S^1$, then morphisms from $[n]$ to $[m]$ in the category $\LR$
  are homotopy classes of monotonous continuous maps $f:[n] \to [m]$
  which send marked points to marked points and have positive
  degree, $\deg f \geq 1$.
\end{defn}

The only difference with the category $\Lambda$ is that the maps are
allowed to have degree bigger than $1$. A typical new map is obtained
as follows. For every configuration of $n$ points on a circle and
any positive integer $l \geq 1$, consider the $l$-fold \'etale cover
$\pi_l:S^1 \to S^1$, and the configuration of $nl$ preimages of $n$
marked points. Then $\pi$ gives a well-defined map $\pi_{n,l}:[nl]
\to [n]$ in the category $\LR$. Moreover, every map $f:[m] \to [m]$
of degree $l$ in the category $\LR$ factors as $f = \pi_{n,l} \circ
f'$ for some $f':[m] \to [nl]$ of degree $1$, and such a
factorization is unique up to the action of the group $\Z/l\Z$ of
deck transformations of the covering $\pi_l:S^1 \to S^1$. Thus the
set $\LR_l([m],[n])$ of degree-$l$ maps from $[m]$ to $[n]$ is
naturally identified with the quotient
\begin{equation}\label{lr.quo}
\LR_l([m],[n]) = \Lambda([m],[nl])/(\Z/l\Z)
\end{equation}
by the action of the group $\Z/l\Z$ generated by $\sigma^n:[nl] \to
[nl]$. In particular, $\LR_l([m],[n])$ is finite for every $[m]$,
$[n]$ and $l$.

\begin{defn}
A map $f:[m] \to [n]$ in the category $\LR$ is {\em horizontal} if
it is of degree $1$. A map $f:[m] \to [n]$ of some degree $l \geq 1$
is {\em vertical} if the map $f':[m] \to [nl]$ in the decomposition
$f = \pi_{n,l} \circ f$ is invertible.
\end{defn}

It follows from the discussion above that vertical and horizontal
maps form a factorization system on $\LR$ in the sense of
Definition~\ref{facto.def}. The subcategory $\LR_h$ formed by
horizontal maps is by definition equivalent to $\Lambda$. Moreover,
for any group $G$, let $O_G$ be the category of {\em finite
$G$-orbits} -- that is, finite sets equipped with a transitive
$G$-action. Then the subcategory $\LR_v \subset \LR$ formed by
vertical maps is obviously equivalent to the orbit category $O_\Z$
-- the equivalence sends $[n] \in \LR$ to the orbit $\Z/n\Z$ (all
finite $\Z$-orbits are of this form).

\begin{lemma}\label{lr.cart}
For any pair of a horizontal map $h:[m_1] \to [m]$ and a vertical
map $v:[m_2] \to [m]$ in $\LR$, there exists a Cartesian square
$$
\begin{CD}
[m_{12}] @>{h_1}>> [m_2]\\
@V{v_1}VV @VV{v}V\\
[m_1] @>{h}>> [m]
\end{CD}
$$
with horizontal $h_1$ and vertical $v_1$.
\end{lemma}

\proof{} Clear. \endproof

Composing the functors $i$ and $\pi$ of \eqref{i.p.def} with the
natural embedding $\Lambda \cong \LR_h \hookrightarrow \LR$, we
obtain functors
$$
\wt{i},\wt{\pi}:\LI \to \LR.
$$
Moreover, the quotient maps $\pi_{n,l}$, $n,l \geq 1$ taken together
define a vertical map
\begin{equation}\label{i.p.v.def}
\wt{v}:\wt{i} \to \wt{\pi}.
\end{equation}
Let $\wt{\LI}$ be the category of vertical maps $v:[m] \to [m']$ in
$\LR$, with maps from $v_1:[m_1] \to [m_1']$ to $v_2:[m_2] \to
[m_2']$ given by commutative squares
\begin{equation}\label{lzred.maps}
\begin{CD}
[m_1] @>{f}>> [m_2]\\
@V{v_1}VV @VV{v_2}V\\
[m'_1] @>{f'}>> [m'_2]
\end{CD}
\end{equation}
with horizontal $f$, $f'$. Then sending $a \in \LI$ to
$\wt{v}:\wt{i}(a) \to \wt{\pi}(a)$ defines a functor
\begin{equation}\label{li}
\LI \to \wt{\LI}.
\end{equation}

\begin{lemma}\label{li.lemma}
The functor \eqref{li} is an equivalence of categories.
\end{lemma}

\proof{} Clear.\endproof

Sending a wheel $[n] \in \LR$ to the set $V([n])$ of its vertices
defines a functor $\LR \to \Sets$. We let
\begin{equation}\label{der.deor}
\begin{CD}
\DeR @>\wt{j}>> \LR^{opp}
\end{CD}
\end{equation}
be the discrete fibration corresponding to $V$ by the Grothendieck
construction.

\begin{lemma}\label{der.le}
The functor $\delta = \deg \circ \wt{j}:\DeR \to [1/\N^*]$ is a
cofibration, with fiber $\Delta$, and transition functor $r_m$
corresponding to $m \in \N^*$ given by the edgewise subdivision
functor.
\end{lemma}

\proof{} By definition, $\DeR$ is opposite to the full subcategory
in the slice category $[1]\backslash\LR$ spanned by horizontal maps
$h:[1] \to [n]$, $[n] \in \LR$. Moreover, it inherits from $\LR$ the
vertical/horizontal factorization system. One now immediately
deduces that vertical maps are Cartesian with respect to $\delta$,
while the fiber of $\delta$ is spanned by horizontal maps.
\endproof

\subsection{Extended categories.}\label{lz.subs}

Let $\N^*$ be the monoid of positive integers $l \geq 1$ with
respect to multiplication, and let $[1/\N^*]$ be the category with
one object $1$ and
$$
\Hom_{[1/\N^*]}(1,1) = \N^*.
$$
Sending a map to its degree gives then a functor
\begin{equation}\label{deg}
\deg:\LR \to [1/\N^*].
\end{equation}
This functor has a section, the fully faithful embedding
$\alpha:[1/\N^*] \to \LR$ which sends $1$ to $[1] \in \LR$ (any map
$[1] \to [1]$ is uniquely determined by its degree). Moreover, let
$$
I = 1\backslash{}[1/\N*]
$$
be the category of objects $a \in [1/\N^*]$ equipped with a map $1
\to a$ (the slice category). Equivalently, $I$ is $\N^*$ considered
as a partially ordered set with order given by divisibility, and
turned into a category in the standard way. Then we have a natural
cofibration
\begin{equation}\label{I.N}
I \to [1/\N^*]
\end{equation}
whose fiber is the set $\N^*$ considered as a discrete category. By
the Grothen\-dieck construction, this corresponds to a functor
$[1/\N^*] \to \Sets$ sending $1$ to $\N^*$, in other words, to an
action of the monoid $\N^*$ on itself; the action is by right
multiplication.

Let now $\N^*$ act on itself both on the right and on the left, and
let $\I$ be corresponding category cofibered over $[1/\N^*] \times
[1/\N^*]$ with fiber $\N^*$. Equivalently, $\I$ is obtained by the
Cartesian square
$$
\begin{CD}
\I @>>> I\\
@VVV @VVV\\
[1/\N^*] \times [1/\N^*] @>>> [1/\N^*],
\end{CD}
$$
where the bottom map is induced by the product map $\N^* \times \N^*
\to \N^*$. Composing the cofibration $\I \to [1/\N^*] \times
[1/\N^*]$ with the projection onto the right multiple $[1/\N^*]$, we
obtain a cofibration
\begin{equation}\label{II.N}
\I \to [1/\N^*],
\end{equation}
with the fiber $I$. We also have a natural Cartesian functor $I \to
\I$; on fibers, it is given by the inclusion of the discrete
category $\N^*$ into $I$ (which is nothing but $\N^*$ considered as
a partially ordered set). Explicitly, objects in $\I$ are positive
integers $n \geq 1$, and morphisms are generated by morphisms
\begin{equation}\label{f.r}
F_l,R_l:n \to nl
\end{equation}
for any $n,l \geq 1$, subject to relations
$$
F_n \circ F_m = F_{nm}, \quad R_n \circ R_m = R_{nm}, \quad F_n
\circ R_m = R_m \circ F_n
$$
for any $n,m \geq 1$. This is opposite to the category introduced by
T. Goodwillie in \cite{good}.

\begin{defn}
The {\em extended cyclic category} $\LZ$ is given by
$$
\LZ = \LR \times_{[1/\N^*]} I,
$$
and the {\em extended cyclotomic category} $\wLR$ is given by
$$
\wLR = \LR \times_{[1/\N^*]} \I,
$$
where in both cases, $\LR \to [1/\N^*]$ is the degree functor $\deg$
of \eqref{deg}, and $I \to [1/\N^*]$ resp. $\I \to [1/\N^*]$ are the
cofibrations \eqref{I.N} resp. \eqref{II.N}.
\end{defn}

The section $\alpha:[1/\N^*] \to \LR$ of the degree functor
$\deg:\LR \to [1/\N^*]$ induces a functor
\begin{equation}\label{alpha}
\wt{\alpha}:\I \to \wLR.
\end{equation}
By definition, we have cofibrations 
\begin{equation}\label{lz.cofib}
\lambda:\LZ \to \LR, \qquad \wt{\lambda}:\wLR \to \LR,
\end{equation}
with fibers identified with $\N^*$ resp. $I$. Explicitly, objects in
either $\LZ$ or $\wLR$ are given by pairs $\langle [n] \in \LR, l
\in \N^*\rangle$; we will denote sich a pair by $[n|m]$. A map from
$[n|m]$ to $[n',m']$ in $\LZ$ resp. $\wLR$ is a map $f:[n] \to [n']$
in $\LR$ such that $m'=m\deg f$, resp. $m' = lm\deg f$ for some
integer $l \geq 1$. The vertical/horizontal factorization system
then induces an analogous vertical/horizontal factorization systems
on $\LZ$ and $\wLR$. Denote by $\LZ_v,\LZ_h \subset \LZ$,
$\wLR_h,\wLR_h \subset \wLR$ the subcategories spanned by vertical,
resp. horizontal maps. Then $\LZ_h$ and $\wLR_h$ decompose as
\begin{equation}\label{lz.dec.h}
\LZ_h = \coprod_{m \geq 1}\LZ_h^m \cong \N^* \times
\Lambda,\qquad\qquad \wLR_h = I \times \Lambda,
\end{equation}
where $\LZ_h^m$ is the full subcategory spanned by the objects
$[n|m]$, $n \geq 1$, and for every $m$, the category $\LZ_h^m$ is
naturally equivalent to $\Lambda$. On the other hand, the category
$\LZ_v$ decomposes as
\begin{equation}\label{lz.dec.v}
\LZ_v = \coprod_{n \geq 1}\LZ_v^n \cong \coprod_{n \geq 1}O_{\Z/n\Z},
\end{equation}
where $\LZ_h^n$ is the full subcategory spanned by objects $[n'|m']$
with $n=n'm'$, and for every $n$, the category $\LZ_h^n$ is
naturally equivalent to the category $O_{\Z/n\Z}$. These
decompositions are induced by the identifications $\LR_h \cong
\Lambda$, $\LR_v \cong O_{\Z}$.

\subsection{More on extended cyclic category.}\label{wlz.subs}

Here are some more simple properties of the category $\LZ$. By
Yoneda, the category $\LZ$ is fully and faithfully embedded into the
category $\Fun(\LZ^{opp},\Sets)$. Restricting to $\Lambda = \LZ_h^1
\subset \LZ$, we obtain a functor
$$
Y:\LZ \to \Fun(\Lambda^{opp},\Sets).
$$
This gives an alternative purely combinatorial description of the
category $\LZ$, since we have the following result.

\begin{lemma}\label{Y.le}
The functor $Y$ is fully faithful.
\end{lemma}

\proof{} By definition, $Y([n|1])$ is the functor
$h_{[n]}:\Lambda^{opp} \to \Sets$ represented by $[n] \in \Lambda$,
and by \eqref{lr.quo}, $Y([n|m])$ for $m \geq 2$ is the quotient
$h_{[nm]}/(\Z/m\Z)$ of $h_{[nm]}$ by the action of the cyclic group
$\Z/m\Z \subset \Aut([nm])$. We then have
\begin{align*}
\Hom(Y([n|m]),Y([n'|m'])) &=
\Hom(h_{[nm]},h_{[n'm']}/(\Z/m'\Z))^{\Z/m\Z}\\ 
& = \left(\Lambda([nm],[n'm'])/(Z/m'\Z)\right)^{\Z/m\Z},
\end{align*}
and by Lemma~\ref{twt.lemma}, this is non-empty only if $m' = lm$
for some $l$, and coincides with
$$
\Lambda([n,n'l])/(\Z/l\Z).
$$
Again by \eqref{lr.quo}, this coincides with $\LZ([n|m],[n'|m']) =
\LR_l([n],[n'l])$.
\endproof

\begin{lemma}\label{pb.lemma}
For any pair of a horizontal map $h:[n_1|m] \to [n|m]$ and a
vertical map $v:h:[n'|m'] \to [n|m]$, there exists a cartesian
square
\begin{equation}\label{pb}
\begin{CD}
[n'_1|m'] @>{v'}>> [n_1|m]\\
@V{h'}VV @VV{h}V\\
[n'|m'] @>{v}>> [n|m]
\end{CD}
\end{equation}
in $\LZ$ with horizontal $h'$ and vertical $v'$.
\end{lemma}

\proof{} Use Lemma~\ref{lr.cart} and notice that $\deg v_1 = \deg
v$.
\endproof

By virtue of this Lemma, we can define a new category $\wLZ$ as
follows: the objects are the same as in $\LZ$, morphisms from $c$ to
$c'$ are given by isomorphism classes of diagrams
\begin{equation}\label{wh.lz}
\begin{CD}
c @<v<< c_1 @>h>> c'
\end{CD}
\end{equation}
with vertical $v$ and horizontal $h$. Composition is given by
pullbacks. Note that a diagram \eqref{wh.lz} has no non-trivial
automorphisms. Therefore $\wLZ$ has a natural factorization
system, with horizontal resp. vertical maps represented by a diagram
with invertible $v$ resp. $h$. As before, we denote the
corresponding subcategories by $\wLZ_h,\wLZ_v \subset
\wLZ$. We have decomposition
$$
\wLZ_v = \coprod_{n \geq 1}O_{\Z/n\Z}^{opp}, \qquad \LZ_h = \N^*
\times \Lambda
$$
induced by \eqref{lz.dec.v} and \eqref{lz.dec.h}, so that we have
$\wLZ_v \cong \LZ_v^{opp}$ and $\wLZ_h \cong \LZ_h$. Note that the
equivalence $\Lambda \cong \Lambda^{opp}$ of \eqref{la.du} gives an
equivalence
\begin{equation}\label{lz.h.opp}
\LZ_h \cong \N^* \times \Lambda \to \LZ_h^{opp} \cong
\N^* \times \Lambda^{opp}.
\end{equation}

\begin{lemma}
There exists en equivalence of categories
$$
\wLZ \cong \LZ^{opp}
$$
which restricts to the identity functor $\wLZ_v \cong
\LZ_v^{opp} \to \LZ_v^{opp} \subset \LZ^{opp}$ and to the
equivalence $\wLZ_h \cong \LZ_h \to \LZ_h^{opp} \subset
\LZ^{opp}$ of \eqref{lz.h.opp}.
\end{lemma}

\proof{} The shape of the equivalence $\wLZ \cong \LZ^{opp}$ is
prescribed by the conditions: it is identical on objects and on
vertical morphisms, and it sends a horizontal morphism $h$ in some
$\LZ_h^m \subset \LZ$, $\LZ_h^m \cong \Lambda$ to $h_\hash$. To see
that this is consistent, we have to check that for any Cartesian
square \eqref{pb}, the diagram
$$
\begin{CD}
[n'_1|m'] @>{v'}>> [n_1|m]\\
@A{h'_\hash}AA @AA{h_\hash}A\\
[n'|m'] @>{v}>> [n|m]
\end{CD}
$$
is commutative. This immediately follows from the construction: both
maps $h$ and $h'$ can be represented by the same map $[n_1m] \to
[nm]$.
\endproof

\subsection{Homological properties.}

The geometric realization $|\Lambda|$ of the nerve of the category
$\Lambda$ is homotopy equivalent to $BU(1)$, the classifying space
of the unit circle group $U(1)$. In particular, the cohomology of
the category $\Lambda \cong \Lambda^{opp}$ with coefficients in the
constant functor $\Z$ is given by
$$
H^\hdot(\Lambda^{opp},\Z) = H^\hdot(\Lambda,\Z) \cong \Z[u],
$$
where $u$ is a generator of degree $2$. To see $u \in
H^2(\Lambda,\Z)$ explicitly, associate to a wheel $[n] \in \Lambda$
its cellular cohomology complex $C^\hdot([n],\Z)$. Since
topologically, every wheel is the circle $S^1 \cong U(1)$, we have
an exact sequence
$$
\begin{CD}
0 @>>> \Z @>>> C^0([n],\Z) @>>> C^1([n],\Z) @>>> \Z @>>> 0
\end{CD}
$$
for any $[n] \in \Lambda$. This sequence depends functorially on
$[n]$, so that we obtain an exact sequence
\begin{equation}\label{4.term}
\begin{CD}
0 @>>> \Z @>{b_0}>> j_*\Z @>{B}>> j^o_!\Z @>{b_1}>> \Z @>>> 0
\end{CD}
\end{equation}
of functors in $\Fun(\Lambda^{opp},\Z)$, where we have identified
$$
\begin{aligned}
C^1([n],\Z) &\cong \Z[E([n])] \cong (j^o_!\Z)([n]),\\
C^0([n],\Z) &\cong \Z[V([n])]^* \cong (j_*\Z)([n]),
\end{aligned}
$$
with $j:\Delta \to \Lambda^{opp}$, $j^o:\Delta^{opp} \to
\Lambda^{opp}$ as in Subsection~\ref{connes.subs}.  The exact
sequence \eqref{4.term} represents by Yoneda an element $u \in
\Ext^2(\Z,\Z)$; this is the generator of the polynomial algebra
$\Ext^\hdot(\Z,\Z) = H^\hdot(\Lambda^{opp},\Z)$.

The cellular cohomology complexes are also functorial with respect
to maps of higher degrees, so that the exact sequence \eqref{4.term}
extends to the category $\Fun(\LR,\Z)$. The extended sequence takes
the form
\begin{equation}\label{4.term.R}
\begin{CD}
0 @>>> \Z @>{b_0}>> \wt{j}_*\Z @>{B}>> E @>{b_1}>> \deg^*\Z(1)
@>>> 0,
\end{CD}
\end{equation}
where $\wt{j}$ is as in \eqref{der.deor}, $E \in \Fun(\LR^{opp},\Z)$
denotes the functor $[n] \mapsto C^1([n],\Z)$, and $\Z(1) \in
\Fun([1/\N^*],\Z)$ is the functor corresponding to $\Z$ with every
$n \in \N^*$ acting by multiplication by $n$.

Using the functors $i$, $\pi$ of \eqref{i.p.def}, we can pull back
the sequence \eqref{4.term} to the category $\LI^{opp}$ in two ways,
and the map \eqref{i.p.v.def} induces a map between these
pullbacks. After fixing a positive integer $n \geq 2$ and
restricting to $\Lambda_n^{opp} \subset \LI^{opp}$, what we obtain
is a commutative diagram
\begin{equation}\label{n.ti}
\begin{CD}
0 @>>> \Z @>>> i_n^*j_*\Z @>{i_n^*B}>> i_n^*j^{opp}_!\Z @>>> \Z @>>>
0\\
@. @V{\id}VV @VV{v_n}V @VV{v_n}V @VV{n\id}V\\
0 @>>> \Z @>>> \pi^*_nj_*\Z @>{\pi^*_nB}>> \pi^*_nj^{opp}_!\Z
@>>> \Z @>>> 0,
\end{CD}
\end{equation}
where the vertical maps $\eta_n$ are isomorphisms. The geometric
realization $|\Lambda_n|$ of the category $\Lambda_n$ has the same
homotopy type as $|\Lambda|$, and the functor $i_n:\Lambda_n \to
\Lambda$ induces a homotopy equivalence of the realizations. The
first line in \eqref{n.ti} represents by Yoneda the generator $u \in
H^2(\Lambda_n^{opp},\Z)$ of the cohomology algebra
$$
H^\hdot(\Lambda_n^{opp},\Z) \cong H^\hdot(\Lambda^{opp},\Z) \cong
\Z[u].
$$
On the other hand, the functor $\pi_n:\Lambda_n \to \Lambda$ does
{\em not} induce a homotopy equivalence of realizations: on the
level of realizations $|\Lambda_n| \cong |\Lambda| \cong BU(1)$, the
map induced by $\pi_n$ corresponds to the $n$-fold covering $U(1)
\to U(1)$. the second line in \eqref{n.ti} represents the elemend
$\pi_n^* = nu \in H^2(\Lambda_n^{opp},\Z)$.

\subsection{Homological vanishing.}

We will also need some cohomological vanishing results on small
categories $\DeR$ and $\LR^{opp}$. First recall that if we take any
$p \geq 1$ and let $r_p:\Delta \to \Delta$ be the edgewise
subdivison functor \eqref{edge}, then for any $M \in \D(\Delta,k)$,
the natural map
$$
H^\hdot(\Delta,r_p^*M) \to H^\hdot(\Delta,M)
$$
is an isomorphism (see e.g. \cite[Lemma 1.14]{K}, but in fact the
statement is very well-known). By adjunction, this means that the
natural map
\begin{equation}\label{edge.2}
r_{p!}k \to k
\end{equation}
is a quasiisomorphism for any $p \geq 1$ (here $k$ denotes the
constant functors).

\begin{lemma}\label{no.h}
Let $h:\Delta \cong \DeR_h \to \DeR$ be the natural embedding, and
assume given an object $M \in \D(\DeR^{opp},\Z)$ such that
$H^\hdot(\Delta,M)=0$. Then $H^\hdot(\DeR,M) = 0$.
\end{lemma}

\proof{} Let $\delta:\DeR \to [1/\N^*]$ be the cofibration of
Lemma~\ref{der.le}. It suffices to prove that $R^\hdot\delta_*M =
0$. Equivalently,
$$
H^\hdot(\wDer,\kappa^*M)=0,
$$
where $\wDer$ is the category of objects $[n] \in \DeR$ equipped
with a map $\delta([n]) \to 1$ --- in other words, a number $m \in
\N^*$ --- and $\kappa:\wDer \to \DeR$ is the forgetful functor. For
any $l \geq 1$, let $\iota_l:\Delta \to \wDer$ be the embedding
sending $[n]$ to itself equipped with the number $l$. We also have
the forgetful functor $\delta:\wDer \to I$ induced by $\delta:\DeR
\to [1/\N^*]$, and by Lemma~\ref{der.le}, this is a cofibration with
transition functors $r_n$. Therefore the functor
$i_l^*i_{m!}:\D(\Delta,k) \to \D(\Delta,k)$ is trivial unless $m=lp$
for some $p$, and isomorphic to $r_{p!}$ if this is the case. Then
by \eqref{edge.2}, the adjunction map $i_{l!}k \to k$ is an
isomorphism at $\iota_p(\Delta) \subset \wDer$ for all $p$ dividing
$l$. Therefore
$$
k \cong \lim_{\to}i_{l!}k,
$$
so that it suffices to prove that
$$
\RHom^\hdot(i_{l!}k,\kappa^*M) = 0
$$
for any $l$. By adjunction, this means
$H^\hdot(\Delta,\iota_l^*\kappa^*M) = 0$, and we have $\kappa \circ
\iota_l = \id$.
\endproof

\begin{corr}\label{no.h.cor}
For any $M \in \Fun(\LR^{opp},k)$, we have
$$
H^\hdot(\LR^{opp},M \otimes E) = 0,
$$
where $E \in \Fun(\LR^{opp},k)$ is as in \eqref{4.term.R}.
\end{corr}

\proof{} As in \cite[Lemma 1.10]{K}, \eqref{4.term.R} implies by
devissage that it suffices to prove that
$$
H^\hdot(\DeR,\wt{j}^*(M \otimes E)) = 0,
$$
and by Lemma~\ref{no.h}, it suffices to prove that
$H^\hdot(\Delta,j^*h^*(M \otimes E))=0$. This is dual to \cite[Lemma
  1.10]{K}.
\endproof

By Corollary~\ref{no.h.cor}, \eqref{4.term.R} induces a long exact
sequence of cohomology
\begin{equation}\label{connes.lr}
\begin{CD}
H^\hdot(\LR^{opp},M) @>>> H^\hdot(\DeR,\wt{j}^*M) @>>>\\
@>>> H^{\hdot-1}(\LR^{opp}, M \otimes \deg^*\Z(1)) @>>>
\end{CD}
\end{equation}
for any $M \in \Fun(\LR^{opp},\Z)$.

\begin{prop}\label{profini}
Assume given a functor $M \in \Fun(\LR^{opp},\Z)$, and assume
further that $M$ is profinitely complete. Then the natural map
$$
H^\hdot(\LR^{opp},M) \to H^\hdot(\DeR,\wt{j}^*M)
$$
is an isomorphism.
\end{prop}

\proof{} By the projection formula, we have
$$
H^\hdot(\LR^{opp},M \otimes \deg^*\Z(1)) \cong
H^\hdot([1/\N^*],\Z(1) \otimes R^\hdot\deg_*M),
$$
and by \eqref{connes.lr}, we have to prove that this is
tautologically $0$. Since $R^\hdot\deg_*$ commutes with profinite
completion, it suffices to prove the following.

\begin{lemma}
For any profinitely complete $M \in \Fun([1/\N^*],\Z)$, we have
$$
H^\hdot([1/\N^*],M \otimes \Z(1)) = 0.
$$
\end{lemma}

\proof{} Every profinitely complete abelian group $M$ decomposes as
$$
M = \prod_p M_p,
$$
where the product is over all primes $p$, and $M_p$ is pro-$p$
complete. Therefore we may assume that $M$ is pro-$p$complete for
some prime $p$. Decompose the monoid $\N^*$ into the product $\N^* =
\N \times \N_{\{p\}}^*$, where $\N_p \subset \N$ consists of all
positive integers prime to $p$, and $\N \subset \N^*$ is the monoid
of all power $p^n$, $n \geq 0$. Then by the K\"unneth formula, it
suffices to prove that
$$
H^\hdot([1/\N],\Z(1) \otimes M) = 0.
$$
But for any $M' \in \Fun([1/\N],\Z)$, $H^\hdot([1/\N],M')$ can
obviously be computed by the two-term complex
$$
\begin{CD}
M' @>{\id - t}>> M',
\end{CD}
$$
where $t:M' \to M'$ is the action of the generator $1 \in \N$. In our
case $M' = M \otimes \Z(1)$, $M'$ is pro-$p$ complete, and $t$ is
divisible by $p$. Therefore $\id - t$ is invertible.
\endproof

\section{Cyclic Mackey functors.}\label{mack.sec}

\subsection{The quotient category description.}

Consider the wreath product $\LZ\wr\Gamma$ of the enhanced cyclic
category $\LZ$ with the category $\Gamma$ of finite sets;
explicitly, $\LZ\wr\Gamma$ can be identified with the full
subcategory in $\Fun(\Lambda^{opp},\Sets)$ spanned by finite
disjoint unions of objects $[n|m] \in \LZ \subset
\Fun(\Lambda^{opp},\Sets)$ (in particular, we have a natural full
embedding $\LZ \subset \LZ\wr\Gamma$). A map
$$
\coprod_{s \in S}[n_s|m_s] \to \coprod_{s' \in S'}[n_{s'}|m_{s'}]
$$
between two such unions indexed by finite sets $S$, $S'$ consists of
a map $f:S \to S'$ and a map $f_s:[n_s|m_s] \to [n_{f(s)}|m_{f(s)}]$
for any $s \in S$. Say that a map $\langle f,\{f_s\}\rangle$ is {\em
vertical} if $f_s$ is vertical for any $s$. Say that the map is {\em
  horizontal} if $f$ is invertible, and $f_s$ is horizontal for any
$s$. Then vertical and hozirontal maps obviously form a
factorization system on $\LZ\wr\Gamma$, and we have the following.

\begin{lemma}\label{pb.quo}
For any vertical map $v:[a] \to [b]$ in $\LZ\wr\Gamma$ and any map
$f:[b'] \to [b]$, there exists a Cartesian square
$$
\begin{CD}
[a'] @>{v'}>> [b']\\
@V{f'}VV @VV{f}V\\
[a] @>{v}>> [b]
\end{CD}
$$
with vertical $v'$.
\end{lemma}

\proof{} Since vertical and horizontal maps form a factorization
system on $\LZ\wr\Gamma$, we may assume that $f$ is either
horizontal or vertical. In the first case, the claim follows from
Lemma~\ref{pb.lemma}. In the second case, it suffices to prove that
the category $(\LZ\wr\Gamma)_v = \LZ_v\wr\Gamma$ has pullbacks. This
folows from \eqref{lz.dec.v}, since for every $m$, the category
$O_{\Z/m\Z}\wr\Gamma$ has pullbacks (it is equivalent to the
category of finite sets equipped with an action of the group
$\Z/m\Z$).
\endproof

Our definition of (derived) cyclic Mackey functors mimicks the
defition of derived Mackey functors given in \cite[Section 3]{Ka-ma}
(and more specifically, in Subsection~3.4). For any two objects
$c,c' \in \LZ$, we let $Q^\wr\LZ(c,c')$ be the category of diagrams
\begin{equation}\label{domik}
\begin{CD}
c @<{v}<< c_1 @>{f}>> c'
\end{CD}
\end{equation}
in $\LZ\wr\Gamma$ with vertical $v$. Maps from a diagram $c \gets
c_1 \to c'$ to a diagram $c \gets c_2 \to c'$ are given by such maps
$g = \langle g,\{g_s\} \rangle:c_1 \to c_2$ that $g$ commutes with
$f$ and with $v$, and each of the component maps $g_s$ is
invertible.

We obviously have $Q^\wr\LZ(c,c') = Q(\LZ(c,c'))\wr\Gamma$, where
$Q\LZ(c,c') \subset Q^\wr\LZ(c,c')$ is the subcategory of diagrams
\eqref{domik} with $c_1 \in \LZ \subset \LZ\wr\Gamma$, and
invertible maps between them. This identification gives the
projection functor $\rho_{c,c'}:Q^\wr(c,c') \to \Gamma$ sending a
diagram $c \gets c_1 \to c'$ to the finite set $S$ of components of
the object $c_1 = \copr_{s \in S}[n_s|m_s] \in \LZ\wr\Gamma$. We let
$$
T_{c,c'} = \rho^*T \in \Fun(Q^\wr(c,c'),\Z)
$$
be the functor from $\Q^\wr(c,c')$ to the category of abelian groups
obtained by pullback from the functor $T \in \Fun(\Gamma,\Z)$, $T(S)
= \Z[S]$.

By Lemma~\ref{pb.quo}, for any $c,c',c'' \in \LZ$ we have a natural
functor
$$
m_{c,c',c''}:Q^\wr(c,c') \times Q^\wr(c',c'') \to Q^\wr(c,c'')
$$
given by pullback. This operation is associative, so that we have a
$2$-category $\Q\LZ$ with the same objects as $\LZ$, and with
morphism categories $Q^\wr(-,-)$. Analogously to
\cite[(3.7)]{Ka-ma}, we have natural maps
$$
\mu_{c,c',c''}:T_{c,c'} \boxtimes T_{c',c''} \to m_{c,c',c''}^*T_{c,c''},
$$
and these maps are associative on triple products. Therefore by
\cite[Subsection 1.6]{Ka-ma}, we have an $A_\infty$-category
$\B_\idot$ with the same objects as $\LZ$, with morphisms given by
$$
\B_\idot(c,c') = C_\idot(Q^\wr(c,c'),T_{c,c'}),
$$
the bar complex of the category $Q^\wr(c,c')$ with coefficients in
the functor $T_{c,c'}$, and with the compositions induced by the
functors $m_{c,c'}$ and the maps $\mu_{c,c',c''}$.

For any abelian category $\Ab$, consider the derived category
$\D(\B_\idot^{opp},\Ab)$ of $A_\infty$-functors from the opposite
category $\B^{opp}_\idot$ to the category of complexes of objects in
$\Ab$. By definition, the category $\LZ$ is embedded into the
$2$-category $\Q\LZ$ --- the embedding functor $q:\LZ \to \Q\LZ$ is
identical on objects, and sends a morphism into a diagram
\eqref{domik} with $v = \id$. For any $c,c' \in \LZ$, the
restriction $q^*T_{q(c),q(c')}$ is the constant functor $\Z$, so
that by restriction, we obtain a natural functor
\begin{equation}\label{q.eq}
q^*:\D(\B^{opp}_\idot,\Ab) \to \D(\LZ^{opp},\Ab).
\end{equation}
Let $h:\LZ_h \to \LZ$ be the natural embedding; composing $q^*$ with
$h^*$, we obtain a restriction functor
$$
\D(\B^{opp}_\idot,\Ab) \to \D(\LZ_h^{opp},\Ab).
$$

\begin{defn}
A {\em cyclic Mackey functor} with values in an abelian category
$\Ab$ is an $A_\infty$-functor $M \in \D(\B^{opp}_\idot,\Ab)$
whose restriction $h^*q^*M \in \D(\LZ_h^{opp},\Ab)$ is locally
constant in the sense of Definition~\ref{norm.defn}.
\end{defn}

Cyclic Mackey functors form a full triangulated subcategory in the
category $\D(\B^{opp}_\idot,\Ab)$ which we will denote by $\DML(\Ab)
\subset \D(\B^{opp}_\idot,\Ab)$. In fact, in this paper we will only
need the case $\Ab=k\amod$, the category of modules over a
commutative ring $k$; to simplify notation, we will denote
$\DML(k\amod) = \DML(k)$.

By definition, for every positive integer $m \geq 1$, the embedding
$O_{\Z/m\Z} \cong \LZ_v^m \subset \LZ$ of \eqref{lz.dec.v} extends
to a $2$-functor
\begin{equation}\label{i.m}
i_m:\overline{\Q}^\wr O_{\Z/m\Z} \to \Q\LZ,
\end{equation}
where $\Q^\wr O_{\Z/m\Z}$ is the $2$-category of \cite[Subsection
  3.4]{Ka-ma}, and $\overline{\Q}^\wr O_{\Z/m\Z} \subset \Q^\wr
O_{\Z/m\Z}$ is the full subcategory spanned by $\Z/m\Z$-orbits. This
$2$-functor is compatible with the coefficients $T_{c,c'} \in
\Fun(Q^\wr(c,c'),\Z)$, so it extends to an $A_\infty$-functor
$$
\wt{\B}^{O_{\Z/m\Z}}_\idot \to \B_\idot,
$$
where $\wt{\B}^{O_{\Z/m\Z}}_\idot$ is as in \cite[Subsection
  3.5]{Ka-ma}. The category $\wt{\B}^{O_{\Z/m\Z}}_\idot$ is by
definition self-dual; thus by restriction, for any commutative ring
$k$ we obtain a natural functor
$$
i_m^*:\D(B^{opp}_\idot,k) \to \DM(\Z/m\Z,k),
$$
where $\DM(\Z/m\Z,k)$ is the category of derived $\Z/m\Z$-Mackey
functors constructed in \cite{Ka-ma}. We also have the left-adjoint
functor
$$
i_{m!}:\DM(\Z/m\Z,k) \to \D(B^{opp}_\idot,k).
$$

\subsection{Geometric fixed points.}

Consider now the category $\wLZ$ of Subsection~\ref{wlz.subs}. This
category $\wLZ$ is also embedded into $\Q\LZ$ --- the embedding
functor $\wq:\wLZ \to \Q\LZ$ is identical on objects, and sends a
diagram \eqref{wh.lz} to the corresponding diagram \eqref{domik} (we
recall that diagrams \eqref{wh.lz} have no automorphisms, so there
are no choices involved in this construction). As in the case of the
functor $q$ of \eqref{q.eq}, the functor $\wq$ is compatible with
the coefficients $T_{c,c'}$, so that we have the restriction functor
$$
\wq^*:\D(B^{opp}_\idot,k) \to \D(\wLZ^{opp},k).
$$
For any positive integer $m \geq 1$, the $2$-functor $i_m$ restricts
to the embedding $O_{\Z/m\Z} \to \wLZ$, so that we have a functorial
isomorphism
\begin{equation}\label{bc.1}
i_m^* \circ \wq^* \cong \wq_m^* \circ i_m^*,
\end{equation}
where $\wq_m^*$ is the restriction with respect to the $2$-functor
$\wq_m:O_{\Z/m\Z} \to \Q^\wr(O_{\Z/m\Z})$ considered in
\cite[Subsection 5.3]{Ka-ma} (there denoted by $q$),

Recall now that we have the decompositions \eqref{lz.dec.v},
\eqref{lz.dec.h}, and the identification $\LZ_h \cong \wLZ_h \cong
\N^* \times \Lambda$. Define a functor
$$
\nu:\Fun(\LZ^{opp},k) = \Fun(\wLZ_h^{opp},k) \to \Fun(\wLZ^{opp},k)
$$
by setting $\nu(M)([n|m]) = M([n|m])$ for any $M \in
\Fun(\wLZ_h^{opp},k)$, $[n|m] \in \wLZ$, with 
\begin{equation}\label{nu.eq}
v(M)(v \circ h) = 
\begin{cases}
M(v \circ h), &\quad v\text{ invertible},\\
0, &\quad v\text{ not invertible},
\end{cases}
\end{equation}
where $f = v \circ h$ is the horizontal/vertical factorization of a
map $f$ in $\wLZ_h$. It is easy to see that this is well-defined,
and gives an exact functor $\nu$. By abuse of notation, we will
denote by the same letter $\nu$ its extension to the derived
categories. We have a left-adjoint functor $\phi:\Fun(\wLZ^{opp},k)
\to \Fun(\LZ_h^{opp},k)$, and its derived functor
$$
L^\hdot\phi:\D(\wLZ^{opp},k) \to \D(\LZ_h^{opp},k)
$$
is left-adjoint to $\nu:\D(\LZ_h^{opp},k) \to \D(\wLZ^{opp},k)$.

\begin{defn}\label{Phi.def}
The {\em (geometric) fixed point functor} 
$$
\Phi:\D(B^{opp}_\idot,k) \to \D(\LZ_h^{opp},k)
$$
is given by
$$
\Phi = L^\hdot\phi \circ \wq^*.
$$
\end{defn}

We note that by definition, the fixed points functor $\Phi$ actually
splits into components
$$
\Phi_m:\D(B^{opp}_\idot,k) \to \D((\LZ_h^m)^{opp},k) \cong
\D(\Lambda^{opp},k)
$$
numbered by positive integers $m \geq 1$. We will need some
compatibility results between $\Phi$ and geometric fixed point
functors for derived Mackey functors constructed in
\cite{Ka-ma}. Namely, fix a positive integer $m \geq 1$, and
consider the embedding $i_m$ of \eqref{i.m} and the corresponding
restriction functor $i_m^*$. Let $\bO_{\Z/m\Z}$ be the category of
$\Z/m\Z$-orbits and invertible maps between them, and define the
functor
$$
\nu_m:\Fun(\bO_{\Z/m\Z},k) \to \Fun(O_{\Z/m\Z},k)
$$
by the same formula \eqref{nu.eq} as the functor $\nu$ ($\nu_m(M)(v)
= M(v)$ for invertible $v$, and $0$ otherwise). Let
$\phi_m:\Fun(O_{\Z/m\Z},k) \to \Fun(\bO_{\Z/m\Z},k)$ be its
left-adjoint functor, with the derived functor $L^\hdot\phi_m$. We
then have an obvious isomorphism
\begin{equation}\label{bc.2}
i_m^* \circ \nu \cong \nu_m \circ i_m^*.
\end{equation}
By adjunction, the isomorphisms \eqref{bc.1} and \eqref{bc.2} induce
base change maps
\begin{equation}\label{bc.m}
i_m^* \circ L^\hdot\phi \to L^\hdot\phi_m \circ i_m^*, \qquad i_{m!}
\circ \wq_m^* \to \wq^* \circ i_{m!},
\end{equation}
where by abuse of notation, $i_m$ denotes both the $2$-functor
\eqref{i.m}, and its restrictions
$$
i_m:O_{\Z/m\Z} \to \wLZ^{opp}, \qquad i_m:\bO_{\Z/m\Z} \to
\wLZ^{opp}_h
$$
to $\wq_m(O_{\Z/m\Z})$ and $\wq_m(\bO_{\Z/m\Z})$.

\begin{lemma}\label{bc.le}
For any positive integer $m \geq 1$, the base change maps of
\eqref{bc.m} are isomorphisms.
\end{lemma}

\proof{} Let $\wh{h}:\wLZ^{opp}_h \to \wLZ^{opp}$, $h_m:\bO_{\Z/m\Z}
\to O_{\Z/m\Z}$ be the tautological embeddings. Then we obviously
have $\wh{h}^* \circ \nu \cong \id$, $h_m^* \circ \nu_m^* \cong
\id$, so that by adjunction,
$$
L^\hdot\phi \circ \wh{h}_! \cong \id, \qquad L^\hdot\phi_m \circ
h_{m!}  \cong \id.
$$
On the other hand, we have the horizontal/vertical factorization
system on $\wLZ^{opp} \cong \LZ$, and Lemma~\ref{fact.bc}
yields
$$
i_m^* \circ \wh{h}_! \cong h_{m!} \circ i_m^*.
$$
Therefore the map
$$
i_m^* \circ L^\hdot\phi \circ \wh{h}_! \to L^\hdot\phi_m \circ i_m^*
\circ \wh{h}_!
$$
is an isomorphism, so that the first of the maps \eqref{bc.m}
becomes an isomorphism after evaluating at any object $E \in
\D(\wLZ^{opp},k)$ of the form $E = \wh{h}_!E'$, $E' \in
\D(\wLZ^{opp}_h,k)$. Since the category $\D(\wLZ^{opp},k)$ is
generated by objects of this form, the map $i_m^* \circ L^\hdot\phi
\to L^\hdot\phi_m \circ i_m^*$ is itself an isomorphism.

For the second of the maps \eqref{bc.m}, let $\wh{h}' = \wq \circ
\wh{h}$, $h'_m = \wq_m \circ h_m$. Then again, the
horizontal/vertical factorization system on $\wLZ\wr\Gamma$ shows
that every diagram \eqref{domik} decomposes as
$$
\begin{CD}
c @<{v}<< c_1 @>{v_1}>> c_2 @>{h}>>c'
\end{CD}
$$
with vertical $v$, $v_1$ and horizontal $h$, and this implies that
$$
i_m^* \circ \wh{h}'_* \cong h'_{m*} \circ i_m^*.
$$
The proof is the same as in Lemma~\ref{fact.bc}, except that in
\eqref{fac.eq}, one replaces the $\Hom$-sets $\wLZ^{opp}(-,-)$,
$\wLZ^{opp}_v(-,-)$ with the $\Hom$-categories $Q^\wr\wLZ(-,-)$,
resp. $Q^\wr(O_{\Z/m\Z})(-,-)$. Therefore the base change map
$$
i_m^*(\wq^*E) \to \wq_m^*(i_m^*E)
$$
is an isomorphism for any $E \in \D(\wLZ^{opp},k)$ of the form $E =
\wh{h}'_*E'$, $E' \in \D(\wLZ^{opp}_h,k)$. Objects of this form
generate the derived category, so that $i_m^* \circ \wq^* \cong
\wq_m^* \cong i_m^*$; by adjunction, this gives the claim.
\endproof

\begin{corr}\label{phi.cor}
For any $M \in \D(B^{opp}_\idot,k)$ and any positive integers $m,n \geq 1$, we
have a natural isomorphism
$$
\Phi_m(M)([n]) = \Phi(M)([n|m]) \cong
\Phi^{[(\Z/nm\Z)/(\Z/m\Z)]}i_{mn}^*M,
$$
where $[(\Z/nm\Z)/(\Z/m\Z)] \in O_{\Z/nm\Z}$ is understood as a
$\Z/nm\Z$-orbit, and $\Phi^c$, $c \in O_{\Z/nm\Z}$ is the geometric
fixed points functor of \cite[Subsection 5.1]{Ka-ma}. Moreover, for
any derived Mackey functor $M \in \DM(\Z/mn\Z,k)$, we have a natural
isomorphism
\begin{equation}\label{phi.bc}
\Phi_m(i_{mn!}M) \cong i_{mn!}(\Phi^{[(\Z/nm\Z)/(\Z/m\Z)]}M).
\end{equation}
\end{corr}

\proof{} Immediately follows from Lemma~\ref{bc.le} and
\cite[Proposition 6.5]{Ka-ma}.
\endproof

\subsection{Coalgebra description.}\label{mack.coa.subs}

We will now adapt the constructions of \cite[Section 6]{Ka-ma} and
use the fixed point functor $\Phi$ to obtain another description of
the category $\DML(k)$ of cyclic Mackey functors.

\subsubsection{$A_\infty$-coalgebra.}
Consider the cyclotomic category $\LR$ of Subsection~\ref{lr.subs}
with its vertical/horizontal factorization system. For any object
$[m] \in \LR$, let $[m]\backslash\LR$ be the category of objects
$[m'] \in \LR$ equipped with a map $[m] \to [m']$. The factorization
system on $\LR$ induces a factorization system on $[m]\backslash\LR$
for any $[m]$.

\begin{lemma}\label{pb.quo.1}
For any vertical map $v:[a] \to [b]$ in $[m]\backslash\LR$ and any
map $f:[b'] \to [b]$, there exists a Cartesian square
$$
\begin{CD}
[a'] @>{v'}>> [b']\\
@V{f'}VV @VV{f}V\\
[a] @>{v}>> [b]
\end{CD}
$$
with vertical $v'$.
\end{lemma}

\proof{} The same argument as in Lemma~\ref{pb.quo} shows that the
wreath product category $\LR\wr\Gamma$ has fibered
products. Therefore so does the slice category
$[m]\backslash(\LR\wr\Gamma)$. It remains to notice that the
embedding $[m]\backslash\LR \to [m]\backslash(\LR\wr\Gamma)$ has a
left-adjoint (any $S \in \LR\wr\Gamma$ is a disjoint union of
objects in $\LR$, and any map $[m] \to S$ factors uniquely through
one of these objects).
\endproof

Now, for any morphism $f:[m'] \to [m]$ in $\LR$ and for any integer
$n \geq 0$, consider all diagrams
\begin{equation}\label{v.n}
\begin{CD}
[m'] @>{g}>> [m_0] @>{v_0}>> \dots @>{v_{n-1}}>> [m_n] @>{v_n}>> [m]
\end{CD}
\end{equation}
in $\LR$ such $v_i$, $0 \leq i \leq n$ is vertical, $v_n$ is not
invertible, and $f = v_n \circ \dots \circ v_0 \circ g$. Let
$V_n(f)$ be the groupoid of all such diagrams and isomorphisms
between them. Since $\LR$ is small, $V_n(f)$ is small for any $f$,
$n$.

Assume that $f = f^{(1)} \circ f^{(2)}$ for some $[m''] \in \LR$ and
some maps $f^{(1)}:[m''] \to [m]$, $f^{(2)}:[m'] \to [m'']$. Then
for any $l \geq 1$ and any diagram $\alpha \in V_n(f)$ of the form
\eqref{v.n}, we can apply Lemma~\ref{pb.quo.1} and form a
commutative diagram
$$
\begin{CD}
[m'] @>{g}>> [m_0] @>{v_0}>> \dots @>{v_{n-1}}>> [m_n] @>{v_n}>> [m]\\
@| @A{f^{(1)}_0}AA @. @A{f^{(1)}_n}AA @A{f^{(1)}}AA\\
[m'] @>{g'}>> [m'_0] @>{v'_0}>> \dots @>{v_{n-1}'}>> [m'_n] @>{v_n'}>> [m''],
\end{CD}
$$
where $f^{(2)} = v_n' \circ \dots \circ v_0' \circ g'$, and all the
commutative squares are Cartesian squares in the category
$[m']\backslash\LR$. For any $i$, $0 \leq i \leq n$, we have a
natural vertical map $\nu_i = v'_n \circ \dots \circ v_i:[m'_i] \to
[m'']$. Take the minimal $i$ such that $\nu_i$ is an isomorphism,
let $\alpha^{(2)}$ be the diagram
$$
\begin{CD}
[m'] @>{g'}>> [m'_0] @>{v'_0}>> \dots @>{v_{i-2}'}>> [m'_{i-1}]
@>{\nu_i \circ v_{i-1}'}>> [m'']
\end{CD}
$$
in $V_i(f^{(2)})$, and let $\alpha^{(1)}$ be the diagram
$$
\begin{CD}
[m''] @>{f^{(1)}_i \circ \nu_i^{-1}}>> [m_i] @>{v_i}>> \dots
@>{v_n}>> [m]
\end{CD}
$$
in $V_{n-i}(f^{(1)})$. Then sending $\alpha$ to $\alpha^{(1)} \times
\alpha^{(2)}$ gives a well-defined functor
\begin{equation}\label{copr.vn}
V_n(f) \to \coprod_{0 \leq i \leq n}V_i(f^{(1)}) \times
V_{n-i}(f^{(2)}).
\end{equation}
This construction is obviously associative: for any $l$-tuple
$f_1,\dots,f_l$ of composable maps in $\LR$, we can compose the
functors \eqref{copr.vn} and obtain a functor
\begin{equation}\label{copr.eqa}
V_n(f_1 \circ \dots \circ f_l) \times I_l \to
\coprod_{n_1+\dots+n_l=n}V_{n_1}(f_1) \times \dots \times
V_{n_l}(f_l),
\end{equation}
where $I_l$ is the $l$-th groupoid of the monoidal category operad
of Definition~\ref{I.n.ope}. These functors are compatible with the
natural operad structure on $I_\idot$ in the obvious sense.

For any $i$, $1 \leq i \leq n$, forgetting the object $[m_i]$ in a
diagram \eqref{v.n} gives a functor $\delta_i:V_n(f) \to
V_{n-1}(f)$, and these functors satisfy the relations between
simplicial face maps (not only up to an isomorphism, but on the
nose). Therefore we can define a bicomplex $\T_{\idot,\idot}(f)$ by
setting
\begin{equation}\label{t.bi}
\T_{\idot,\idot}(f) = C_\idot(V_\idot(f),\Z),
\end{equation}
the bar complex of the groupoid $V_\idot(f)$ with coefficients in
the constant functor $\Z$. One differential in the bicomplex
\eqref{t.bi} comes from the bar complex, and the other one is given
by $d = d_1 - d_2 + \dots \pm d_n$, where $d_i$ is the map induced
by the functor $\delta_i$. The coproduct operations \eqref{copr.eqa}
strictly commute with the functors $\delta_i$, so that we have
canonical operations
$$
C_\idot(I_l,\Z) \otimes \T_\idot(f_1 \circ \dots \circ f_l)
\to \T_\idot(f_1) \otimes \dots \otimes \T_\idot(f_n),
$$
again compatible with the asymmetric operad structure on
$C_\idot(I_\idot,\Z)$. Fixing a map $\Ass_\infty \to
C_\idot(I_\idot,\Z)$, as in the Appendix, we equip the collection
$\T_\idot(-)$ with a structure of a $\LR$-graded
$A_\infty$-coalgebra in the sense of \cite[Subsection 1.5.4]{Ka-ma}.

For future use, we record right away some elementary properties of
the $\LR$-graded $A_\infty$-coalgebra $\T_\idot$.

\begin{lemma}\label{t.norm}
\begin{enumerate}
\item The $A_\infty$-coalgebra $\T_\idot$ is augmented, and for any
horizontal map $f$, we have $\T_l(f) = 0$ for $l \geq 1$.
\item For any composable maps $f_1$, $f_2$ such that $f_1$ is
  horizontal, the coproduct map
$$
b_2:\T_\idot(f_2 \circ f_1) \to \T_\idot(f_1) \otimes \T_\idot(f_2)
\cong \T_\idot(f_2)
$$
is an isomorphism.
\item Assume given an $n$-tuple $f_1,\dots,f_n$ of composable maps
  in $\LR$, and assume that $f_i$ is horizontal for $i \geq 2$, and
  $n \geq 3$. Then the corresponding $A_\infty$-operation
$$
b_n:\T_\idot(f_n \circ f_1) \to \T_\idot(f_1) \otimes \dots \otimes
\T_\idot(f_n) \cong \T_\idot(f_n)
$$
is equal to $0$.
\end{enumerate}
\end{lemma}

\proof{} \thetag{i} is obvious: for any $f$, $V_0(f)$ is by
definition a single point, and if $f$ is horizontal, $V_l(f)$ is
empty for any $l \geq 1$. To see \thetag{ii}, note that since $f_1$
is horizontal, the only non-trivial term in the coproduct
\eqref{copr.vn} is a map
\begin{equation}\label{copr.hori}
V_l(f_2 \circ f_1) \to V_l(f_2) \times V_0(f_1) = V_l(f_2),
\end{equation}
and sending a diagram $\alpha \in V_n(f_2)$ of the form \eqref{v.n}
to the diagram
$$
\begin{CD}
[m'] @>{g \circ f_1}>> [m_0] @>{v_0}>> \dots @>{v_{n-1}}>> [m_n] @>{v_n}>> [m]
\end{CD}
$$
defines a map $V_n(f_2) \to V_n(f_2 \circ f_1)$ which is strictly
inverse to \eqref{copr.hori}. Moreover, this inverse map
construction is obviously strictly associative, so that for any
$n$-tuple $f_1,\dots,f_n$ of composable maps with horizontal
$f_2,\dots,f_n$ and $n \geq 3$, we obtain a single map
$$
V_l(f_n) \to V_n(f_n \circ \dots \circ f_1).
$$
Therefore the coproduct map \eqref{copr.vn} is also strictly
associative, so that the map \eqref{copr.eqa} factors through the
natural map $I_n \to \ppt$. This by definition means that the higher
$A_\infty$-operation $b_n$ vanishes, which proves
\thetag{iii}.
\endproof

\begin{lemma}\label{T.taut}
Let $h:\LR_h \cong \Lambda \to \LZ$, $i:\LR_v \cong O_{\Z} \to \LZ$
be the tautological embeddings. Then $h^*\T_\idot$ is the trivial
$\Lambda$-graded $A_\infty$-coalgebra, $h^*\T_\idot(f) \cong \Z$ for
any map $f$ in $\Lambda$, while $i_m^*\T_\idot$ is isomorphic to the
$O_{\Z}$-graded $A_\infty$-coalgebra of \cite[Subsection
6.3.3]{Ka-ma}.
\end{lemma}

\proof{} The first claim immediately follows from
Lemma~\ref{t.norm}. For the second claim, note that if $f$ is
vertical, all the diagrams \eqref{v.n} consist of vertical maps,
thus coincide with the diagrams used in \cite[Subsection
6.3.3]{Ka-ma}.
\endproof

\subsubsection{The comparison theorem.}

Consider now the natural cobifration $\lambda:\LZ \to \LR$ of
\eqref{lz.cofib}, and let $\lambda^*\T_\idot$ be the $\LZ$-graded
coalgebra obtained by pullback. For any ring $k$, consider the
derived category $\D(\LZ,\lambda^*\T_\idot,k)$ of
$A_\infty$-comodules over $\lambda^*\T_\idot$. By
Lemma~\ref{T.taut}, the pullback $h^*\lambda^*\T_\idot$ with respect
to the tautological embedding $h:\LZ_h \to \LZ$ is the trivial
$\LZ_h$-graded $A_\infty$-coalgebra, so that we have a natural
pullback functor
$$
h^*:\D(\LZ,\lambda^*\T_\idot,k) \to \D(\LZ_h^{opp},k).
$$
Let $\DML_\T(k) \subset \D(\LZ,\lambda^*\T_\idot,k)$ be the full
subcategory spanned by objects $M$ whose restriction $h^*M$ is
locally constant in the sense of Definition~\ref{norm.defn}. We want
to show that the category $\DML_\T(k)$ is naturally equivalent to
the category $\DML(k)$ of $k\amod$-valued cyclic Mackey functors.

\medskip

To construct a comparison functor between these two categories, let
$V_l([n|m])$, $[n|m] \in \LZ$, $l \geq 0$ be the groupoid of diagrams
$$
\begin{CD}
[n_1|m_1] @>{v_1}>> \dots @>{v_{n-1}}>> [n_l|m_l] @>{v_n}>> [n|m]
\end{CD}
$$
in $\LZ$ with vertical $v_1,\dots,v_n$ and non-invertible $v_n$. Let
$\sigma_l:V_l([n|m]) \to \wLZ$ be the functor which sends such a
diagram to $[n_1|m_1] \in \wLZ$, or to $[n|m]$ is $l = 0$. For any
$A_\infty$-functor $E_\idot$ from $\B_\idot$ to the category of
complexes of $k$-modules, let $\Phi_\idot^{[n|m]}(E_\idot)$ be the
total complex of the triple complex
$$
C_\idot(V_\idot([n|m]),\sigma_l^*\wq^*E_\idot),
$$
where two differentials are induced by the differentials in
$E_\idot$ and in the bar complex, and the third one is as in
\cite[Subsection 6.3.1]{Ka-ma}. Then the same construction as in
\cite[Subsection 6.3.2]{Ka-ma} shows that the collection
$\Phi_\idot^{[n|m]}(E_\idot)$ has a natural structure of an
$A_\infty$-comodule over $\lambda^*\T_\idot$, so that we obtain a
functor
\begin{equation}\label{Phi.dml}
\Phi_\idot:\D(B^{opp}_\idot,k) \to \D(\LZ,\lambda^*\T_\idot,k).
\end{equation}
Choose an integer $m \geq 1$, and consider the embedding
$i_m:O_{\Z/m\Z} \cong \LZ^{m}_v \to \LZ$ for some $m \geq 1$. The
composition $\lambda \circ i_m:O_{\Z/m\Z} \to \LR$ factors through
the natural embedding $O_{\Z/m\Z} \subset O_{\Z} = \LR_v \to
\LR$. Therefore by Lemma~\ref{T.taut}, the $O_{\Z/m\Z}$-graded
$A_\infty$-coalgebra $i_m^*\lambda^*\T_\idot$ is isomorphic to the
$A_\infty$-coalgebra of \cite[Subsection 6.3.3]{Ka-ma}, and we have
a natural pullback functor
$$
i_m^*:\D(\LZ,\lambda^*\T_\idot,k) \to \DM(\Z/m\Z,k),
$$
where $\DM(\Z/m\Z,k)$ is the category of $k$-valued derived Mackey
functors of \cite{Ka-ma} for the group $\Z/m\Z$. We have an obvious
isomorphism
\begin{equation}\label{h.i}
i_m^* \circ h^* \cong h_m^* \circ i_m^*,
\end{equation}
where $h_m:\bO_{\Z/m\Z} \to O_{\Z/m\Z}$ is the tautological
embedding, and $h_m^*:\DM(\Z/m\Z,k) \to \D(\bO_{\Z/m\Z},k)$ is the
corresponding pullback functor. We note that by \cite[Lemma
6.18]{Ka-ma}, the functor $h_m^*$ admits a right-adjoint functor
$h_{m*}$. Moreover, by construction, the comparison functor
$\Phi_\idot$ of \eqref{Phi.dml} restricts to the corresponding
functor for $O_{\Z/m\Z}$ --- that is, we have a canonical
isomorphism
\begin{equation}\label{Phi.i.m}
\Phi^{[m]}_\idot \circ i_m^*E \cong i_m^* \circ \Phi_\idot,
\end{equation}
where on the left-hand side, $\Phi^{[m]}_\idot$ is the functor
$\Phi_\idot$ of \cite[Theorem 6.17]{Ka-ma}, for $\C = O_{\Z/m\Z}$.

\begin{lemma}\label{h.bc}
The functor $h^*:\D(\LZ,\lambda^*\T_\idot,k) \to \D(\LZ_h^{opp},k)$
admits a right-adjoint
$$
h_*:\D(\LZ_h^{opp},k) \to \D(\LZ,\lambda^*\T_\idot,k),
$$
and for every $m \geq 1$, the base change map $i_m^* \circ h_* \to
h_{m*} \circ i_m^*$ induced by \eqref{h.i} is an isomorphism.
\end{lemma}

\proof{} Immediately follows from Lemma~\ref{cofree} and
Lemma~\ref{cofree.bc}.
\endproof

\begin{lemma}\label{reso}
Assume given a triangulated subcategory $\D'$ in the category
$\D(\LZ,\lambda^*\T_\idot,k)$ which is closed with respects to
arbitrary products and contains all the objects $h_*M$, $M \in
\D(\LZ^{opp}_h,k)$. Then $\D' = \D(\LZ,\lambda^*\T_\idot,k)$.
\end{lemma}

\proof{} For any integer $n \geq 1$, let $h_n:\Lambda \cong \LZ_h^n
\subset \LZ_h$ be the embedding of the $n$-th component of the
decomposition \eqref{lz.dec.h}. For any $M \in \D(\LZ_h^{opp},k)$,
let the support $\Supp(M) \subset \N^*$ be the set of all such $n
\geq 1$ that $h_l^*M$ is non-zero for some $l$ dividing $n$. For any
$M \in \D(\LZ,\lambda^*\T_\idot,k)$, let $\Supp M = \Supp
h^*M$. Note that for any two integers $n,n' \geq 1$, $h_{n'}^* \circ
h_* \circ h_{n*} = 0$ unless $n$ divides $n'$, and $h_n^* \circ h_*
\circ h_{n*} \cong \id$ (by Lemma~\ref{h.bc}, we may check both
statements after applying the functors $i_m^*$, $m \geq 1$, and then
they immediately follow from \cite[Lemma 6.18]{Ka-ma}). Therefore in
particular, for any $n \geq 1$ and any $M \in
\D(\LZ,\lambda^*\T_\idot,k)$, we have
$$
\Supp(h_*h_{n*}h_n^*h^*M) \subset \Supp(M).
$$
Moreover, given an object $M \in \D(\LZ,\lambda^*\T_\idot,k)$, let
$M^{[1]}$ be the cone of the adjunction map
$$
M \to h_*(h^*M)_m,
$$
where $m$ is the smallest integer in $\Supp(M)$. Then
$$
\Supp(M^{[1]}) = \Supp(M) \setminus \{m\}.
$$
By induction, let $M^{[n]} = (M^{[n-1]})^{[1]}$ for any $n >
1$. Then we have natural maps $M^{[n]} \to M$, and their cones
$\wt{M}^{[n]}$ lie in $\D'$ and form an inverse system. We have a
compatible system of maps $\eta_n:M \to \wt{M}^{[n]}$, $n \geq
1$. By induction, for every $n > n' \geq 1$, the transition map
$\wt{M}^{[n+1]} \to \wt{M}^{[n]}$ becomes an isomorphism after
applying $h_{n'}^*h^*$, and so does the map $\eta_n$. Let
$$
\wt{M} = \holim \wt{M}^{[n]},
$$
where $\holim$ is defined by the telescope construction, and let
$\eta:M \to \wt{M}$ be the natural map. Then for every $n \geq 1$,
the inverse system $h_n^*h^*\wt{M}^{[n']}$ stabilizes for $n' > n$,
and $h_n^*h^*(\eta)$ is an isomorphism. Thus $\eta$ itself is an
isomorphism. But by construction, $\wt{M}$ lies in the category
$\D'$.
\endproof

\begin{lemma}\label{Phi.phi}
The composition
$$
h^* \circ \Phi_\idot:\D(B^{opp}_\idot,k) \to \D(\LZ_h^{opp},k)
$$
is isomorphic to the fixed points functor $\Phi$ of
Definition~\ref{Phi.def}.
\end{lemma}

\proof{} By construction, $h^*\Phi_\idot(E)$ for any $E \in
\D(B^{opp}_\idot,k)$ only depends on the restriction $\wq^*E \in
\D(\wLZ^{opp},k)$ --- the construction of the functor $\Phi_\idot$
also gives a functor $\phi_\idot:\D(\wLZ^{opp},k) \to
\D(\LZ^{opp}_h,k)$, and we have $h^* \circ \Phi_\idot \cong
\phi_\idot \circ \wq^*$. For an object $E \in \Fun(\wLZ^{opp},k)$,
the degree-$0$ homology of the complex $\phi_\idot(E)$ is easily
seen to be isomorphic to $\phi(E)$, and this isomorphism is
functorial in $E$. Thus by the universal property of the derived
functor, it extends to a map
$$
e:\phi_\idot \to L^\hdot\phi.
$$
We have to prove that $e:\phi_\idot(E) \to L^\hdot\phi(E)$ is an
isomorphism for any $E \in \D(\wLZ^{opp},k)$ of the form $E =
\wq^*E'$, $E' \in \D(B^{opp}_\idot,k)$. To do this, it suffices to
prove that
$$
i_m^*(e):i_m^* \circ h^* \circ \Phi_\idot \to i_m \circ \Phi
$$
is an isomorphism for any $m \geq 1$. By \eqref{Phi.i.m} and
\eqref{h.i}, the left-hand side is isomorphic to $h_m^* \circ
\Phi^{[m]}_\idot \circ i_m^*$, and by \cite[Lemma 6.15]{Ka-ma},
$h_m^* \circ \Phi^{[m]}$ is isomorphic to the direct sum of the
functors $\Phi^{[(\Z/nm\Z)/(\Z/m\Z)]}$ of
Corollary~\ref{phi.cor}. To finish the proof, it suffices to invoke
Corollary~\ref{phi.cor}.
\endproof

\begin{prop}\label{mack.equi}
The functor $\Phi_\idot$ of \eqref{Phi.dml} is an equivalence of
categories, and it identifies $\DML(k) \subset \D(\B^{opp}_\idot,k)$
with $\DML_\T \subset \D(\LZ,\lambda^*\T_\idot,k)$
\end{prop}

\proof{} As in Lemma~\ref{cofree}, let $\D' \subset
\D(\LZ,\lambda^*\T_\idot,k)$ be the subcategory of objects $M$ in
the category $\D(\LZ,\lambda^*\T_\idot,k)$ such that the functor
$$
\Hom(\Phi_\idot(-),M)
$$
from $\D(B^{opp}_\idot,k)$ to $\D(k)$ is representable. The
geometric fixed points functor $\Phi$ of Subsection obviously has a
right-adjoint. Therefore by Lemma~\ref{Phi.phi}, $\D'$ satisfies all
the assumptions of Lemma~\ref{reso}. Thus $\D'$ is the whole
category $\D(\LZ,\lambda^*\T_\idot,k)$, and $\Phi_\idot$ admits a
right-adjoint functor
$$
\Phi_\idot^{-1}:\D(\LZ,\lambda^*\T_\idot,k) \to \D(B^{opp}_\idot,k).
$$
By Lemma~\ref{Phi.phi}, the composition $\Phi_\idot^{-1} \circ h_*$
is right-adjoint to the geometric fixed points functor $\Phi$ of
Definition~\ref{Phi.def}. Moreover, for every $m \geq 1$, the
functor $\Phi^{[m]}_\idot$ of \eqref{Phi.i.m} is an equivalence of
categories by \cite[Theorem 6.17]{Ka-ma}. Let $\Phi^{-1}_m$ be the
inverse equivalence. Then by Corollary~\ref{phi.cor}, the base
change map
$$
i_m^* \circ \Phi_\idot^{-1} \circ h_* \to \Phi^{-1}_m \circ h_{m*}
\circ i_m^*,
$$
being adjoint to the direct product of isomorphisms \eqref{phi.bc},
is itself an isomorphism. Then by Lemma~\ref{h.bc}, the base change
map 
$$
i_m^*(\Phi_\idot^{-1}(M)) \to \Phi^{-1}_m(i_m^*(M))
$$
is an isomorphism for any $M \in \D(\LZ,\lambda^*\T_\idot,k)$ of the
form $M= h_*M'$, $M' \in \D(\LZ^{opp},k)$, and by Lemma~\ref{reso},
this implies that it is an isomorphism for any $M$. Thus we have
$$
i^*_m \circ \Phi_\idot \circ \Phi^{-1} \cong \Phi^{[m]} \circ
\Phi^{-1}_m \circ i_m^*,
\qquad
i^*_m \circ \Phi^{-1} \circ \Phi_\idot \cong \Phi^{-1}_m \circ
\Phi^{[m]} \circ i_m^*
$$
for any $m \geq 1$. Since $\Phi^{[m]}_\idot$ and $\Phi^{-1}_m$ are
mutually inverse equivalences of categories, this means that the
adjunction maps $\Id \to \Phi^{-1} \circ \Phi_\idot$, $\Phi^{-1}
\circ \Phi_\idot \to \Id$ become isomorphisms after restricting to
$\LZ_v$. Since this restriction is obviously a conservative functor,
we conclude that $\Phi_\idot$ and $\Phi^{-1}$ are mutually inverse
equivalences of categories. To prove that $\Phi_\idot$ identifies
$\DML(k)$ and $\DML_\T(k)$, note that $h^* \circ \Phi_\idot \cong
h^*$.
\endproof

\subsection{Restriction and corestriction.}\label{addi.subs}

We will now make some additional observations on cyclic Mackey
functors for future use. First of all, recall that we have
$2$-functors $q:\LZ \to \Q\LZ$, $\wq:\wLZ \to \Q\LZ$, and the agree
on horizontal maps, $q \circ h \cong \wq \circ \wh{h}$, so that we have
an isomorphism $h^* \circ q^* \cong \wh{h}^* \circ \wq^*$.

\begin{lemma}
The base change map
$$
\wh{h}_! \circ h^* \to \wq^* \circ q_!
$$
induced by the isomorphism $h^* \circ q^* \cong \wh{h}^* \circ \wq^*$
is itself an isomorphism.
\end{lemma}

\proof{} As in the proof of Lemma~\ref{bc.le}, this follows by the
same argument as in the proof of Lemma~\ref{fact.bc}.
\endproof

In particular, this Lemma shows that we have natural identifications
$$
h^* \circ \Phi_\idot \circ q_! \cong \Phi \circ q_! \cong
L^\hdot\phi \circ \wq^* \circ q_! \cong L^\hdot\phi \circ \wh{h}_!
\circ h^*,
$$
and since $L^\hdot\phi \circ \wh{h}_!$ is adjoint to $\wh{h}^* \circ \nu =
\Id$, the right-hand side is just $h^*$. In fact more is true: the
composition
$$
\Phi_\idot \circ q_!:\D(\LZ^{opp},k) \to \D(\LZ,\lambda^*\T_\idot,k)
$$
is naturally isomorphic to the corestriction functor $\xi^*$ with
respect to the augmentation map of the augmented $\LZ$-graded
$A_\infty$-coalgebra $\lambda^*\T_\idot$ (to construct an
isomorphism, resolve a functor $E \in \Fun(\LZ^{opp},k)$ by functors
of the form $i_{m!}E_m$, $E_m \in \Fun(O_{\Z/m\Z}^{opp},k)$ and
apply \cite[Lemma 6.20]{Ka-ma}). Therefore the corestriction functor
$\xi^*$ has a right-adjoint functor $\xi_*$ given by
$$
\xi_* = q^* \circ \Phi_\idot^{-1}:\D(\LZ,\lambda^*\T_\idot,k) \to
\D(\LZ^{opp},k).
$$
We also have 
\begin{equation}\label{xi-xi}
h^* \circ \xi_* \circ \xi^* \cong h^* \circ q^* \circ q_! \cong
\wh{h}^* \circ \wq^* \circ q_! \cong \wh{h}^* \circ \wh{h}_! \circ h^*.
\end{equation}
To compute the right-hand side more effectively, it is useful to
consider the category $\LI \cong \wt{\LI}$ of
Lemma~\ref{li.lemma}. Consider also the product $\LI \times \N^*$,
and define the projections $i,\pi:\LI \times \N^* \to \wLZ_h^{opp}
\cong \LZ_h \cong \Lambda \times \N^*$ by
\begin{equation}\label{i.p.red}
i = i \times \id, \qquad \pi = \pi \times \rho_m \text{ on }
\Lambda_m \times \N^* \subset \LI \times \N^*,
\end{equation}
where $\rho_m:\N^* \to \N^*$ is the map given by multiplication by
the integer $m \geq 1$.

\begin{lemma}\label{i-pi.lemma}
We have a natural isomorphism
$$
\wh{h}^* \circ \wh{h}_! \cong \pi_! \circ i^*.
$$
\end{lemma}

\proof{} Under the identification $\wLZ^{opp} \cong \LZ$, $\wh{h}$
goes to the tautological embedding $h:\LZ_h \to \LZ$. Let
$\overline{\LI}$ be the category of vertical maps $v:a \to a'$ in
$\LZ$, with maps between them given by commutative squares
$$
\begin{CD}
a_1 @>{f}>> a_2\\
@V{v_1}VV @VV{v_2}V\\
a'_1 @>{f'}>> a_2'
\end{CD}
$$
with horizontal $f$ (and arbitrary $f'$). Let $s:\overline{\LI} \to
\LZ_h$, $t:\overline{\LI} \to \LZ$ be the functors sending a map to
its source resp. target. Then $t$ is a cofibration, and $s$ has a
left adjoint $\iota:\LZ_h \to \overline{\LI}$ sending $a \in \LZ$ to
its identity map. Then $\wh{h} = t \circ \iota$, so that
$$
h_! \cong t_! \circ \iota_! \cong t_! \circ s^*.
$$
It remains to notice that we have a natural Cartesian square
$$
\begin{CD}
\wt{\LI} \times \N^* @>{h}>> \overline{\LI}\\
@V{\pi}VV @VV{t}V\\
LZ_h @>{h}>> \LZ,
\end{CD}
$$
so that $h^* \circ t_! \cong \pi_! \circ h^*$ by base change, and we
have $s \circ h = i$.
\endproof

\section{Cyclotomic complexes.}

We can now introduce the main subject of the paper, a cyclotomic
complex. Fix a commutative ring $k$. Consider the cyclotomic
category $\LR$ and the $\LR$-graded $A_\infty$-coalgebra $\T_\idot$
of Subsection~\ref{mack.coa.subs}. Let $\wt{\lambda}:\wLR \to \LR$
be the cofibration of \eqref{lz.cofib}, and let $h:\wLR_h \cong
\Lambda \times I \to \wLR$ be the natural embedding. By
Lemma~\ref{t.norm}, the restriction $h^*\wt{\lambda}^*\T_\idot$ is
the trivial $\wLR_h$-graded $A_\infty$-coalgebra,
$h^*\wt{\lambda}^*\T_\idot(f) = \Z$ for any map $f$ in $\wLR_h$, so
that we have a restriction functor
$$
h^*:\D(\wLR,\wt{\lambda}^*\T_\idot,k) \to \D(\LR_h,k).
$$

\begin{defn}\label{cyclo.defn}
A {\em cyclotomic complex} over $k$ is an $A_\infty$-comodule
$M_\idot$ over $\wt{\lambda}^*\T_\idot$ with values in the category
$k\amod$ such that $h^*M_\idot \in \D(\wLR_h,k)$ is locally constant
in the sense of Definition~\ref{norm.defn}.
\end{defn}

The derived category of cyclotomic complexes over $k$ will be
denoted denoted by $\DLR(k)$.

\subsection{Normalized $\LR$-graded coalgebras.}

Definition~\ref{cyclo.defn} is short enough, but it is only as
explicit and amenable to computations as the $A_\infty$-coalgebra
$\T_\idot$ (that is, not particularly). In this Section, we will
provide more explicit descriptions of the categories $\DLR(k)$. In
the process of doing so, we will also obtain a more convenient
description of the categories $\DML(k)$ of
Section~\ref{mack.sec}. We start with the following reduction
similar to \cite[Subsection 7.5,7.6]{Ka-ma}.

\begin{defn}
Assume given a finite group $G$. A complex $E_\idot$ of
$\Z[G]$-modules is {\em strongly acyclic} with respect to $G$ if for
any subgroup $H \subset G$, $H \neq G$, and any $\Z[H]$-module $V$
we have
$$
\lim_{\overset{n}{\to}} H^\hdot(H,V \otimes F^nE_\idot) = 0,
$$
where $F^\hdot E_\idot$ is the stupid filtration on the complex
$E_\idot$. A map $f:E_\idot \to E'_\idot$ between $\Z[G]$-modules is
a {\em strong quasiisomorphism} with respect to $G$ if its cone is
strongly acyclic.
\end{defn}

\begin{defn}
A $\LR$-graded $A_\infty$-coalgebra $\R_\idot$ is called {\em
normalized} if has the properties \thetag{i}-\thetag{iii} of
Lemma~\ref{t.norm}.
\end{defn}

In particular, the coalgebra $\T_\idot$ is normalized (by
Lemma~\ref{t.norm}). For any normalized $\LR$-graded
$A_\infty$-colagebra $\R_\idot$ and any map $f:[m] \to [m']$ is the
category $\LR$, the complex $\R_\idot(f)$ is equipped, by
definition, with an action of the cyclic group $\Aut([m])$.

\begin{defn}\label{str.qua}
An $A_\infty$-map $\xi:\R_\idot \to \R'_\idot$ between normalized
$\LR$-graded $A_\infty$-coalgebras is a {\em strong
quasiisomorphism} if for any map $f:[m] \to [m']$, the corresponding
map $\xi:\R_\idot(f) \to \R'_\idot(f)$ is a strong quasiisomorphism
with respect to the subgroup $\Aut(f) \subset \Aut([m])$ consisting
of such $g \in \Aut([m])$ that $f \circ g = f$.
\end{defn}

\begin{prop}\label{redu.prop}
For any strong $A_\infty$-quasiisomorphism $\xi:\R_\idot \to
\R'_\idot$ between normalized $\LR$-graded $A_\infty$-coalgebras and
any commutative ring $k$, the corestriction functors
$$
\begin{aligned}
\xi^*:&\D(\wLR,\wt{\lambda}^*\R_\idot,k) \to
\D(\wLR,\wt{\lambda}^*\R'_\idot,k), \\
\xi^*:&\D(\LZ,\lambda^*\R_\idot,k) \to
\D(\LZ,\lambda^*\R'_\idot,k)
\end{aligned}
$$
between the derived categories of $A_\infty$-comodules is an
equivalence of categories.
\end{prop}

\proof{} Both Lemma~\ref{h.bc} and Lemma~\ref{reso} hold for any
normalized $\LR$-graded $A_\infty$-coalgebra instead of $\T_\idot$,
with the same proof. Moreover, they also hold for $\wLR$ instead of
$\LZ$ (again with the same proof). Thus as in the proof of
Proposition~\ref{mack.equi}, the functors $\xi^*$ admits a
right-adjoint functors
$$
\begin{aligned}
\xi_*:&\D(\wLR,\wt{\lambda}^*\R'_\idot,k) \to
\D(\wLR,\wt{\lambda}^*\R_\idot,k), \\
\xi_*:&\D(\LZ,\lambda^*\R'_\idot,k) \to
\D(\LZ,\lambda^*\R_\idot,k).
\end{aligned}
$$
Moreover, for any integer $m \geq 1$, we have a natural embedding
$i_m:O_{\Z/m\Z} \cong \LR_v^m \to \LR$ and the corresponding
restriction functors
\begin{align*}
i_m^*:&\D(\wLR,\wt{\lambda}^*\R_\idot,k) \to
\D(O_{\Z/m\Z},i_m^*\wt{\lambda}^*\R_\idot,k),\\
i_m^*:&\D(\wLR,\wt{\lambda}^*\R'_\idot,k) \to
\D(O_{\Z/m\Z},i_m^*\wt{\lambda}^*\R'_\idot,k).
\end{align*}
For $\LZ$, these functors were already considered in
Section~\ref{mack.sec}. For either $\LZ$ or $\wLR$, we have an
obvious isomorphism $\xi_m^* \circ i_m^* \cong i_m^* \circ \xi^*$,
where $\xi_m = i_m^*(\xi)$, and the corresponding base change map
$i_m^* \circ \xi_* \to \xi_{m*} \circ i_m^*$ is also an isomorphism,
as in Lemma~\ref{h.bc}. Thus it suffices to prove that for every
$m$, $\xi_m^*$ and $\xi_{m*}$ are mutually inverse equivalences of
categories. This is \cite[Lemma 7.14]{Ka-ma}.
\endproof

\subsection{Reduced $\LR$-graded coalgebras.}

The first corollary of Proposition~\ref{redu.prop} is analogous to
the reduction done in \cite[Subsection 7.5]{Ka-ma}.

\begin{defn}\label{red.coa.defn}
A normalized $\LR$-graded $A_\infty$-coalgebra $\R_\idot$ is called
{\em reduced} if for any map $f$ in $\LR$ of degree $n > 1$, $\R(f)
= 0$ unless $n$ is prime. For a normalized $\LR$-graded
$A_\infty$-coalgebra $R_\idot$, its {\em reduction} $\R^{red}_\idot$
is defined by setting
$$
\R^{red}_\idot(f) =
\begin{cases}
\R_\idot(f) &\quad\text{if the degree of $f$ is $1$ or a prime
  number},\\
0, &\quad\text{otherwise},
\end{cases}
$$
with the $A_\infty$-operations being the same as in $\R_\idot$ when
it makes sense, and $0$ otherwise.
\end{defn}

For any normalized $\LR$-graded $A_\infty$-coalgebra $\R_\idot$, we
obviously have a canonical map $\R^{red}_\idot \to \R_\idot$.

\begin{lemma}\label{redu.equi}
For any normalized $\LR$-graded $A_\infty$-coalgebra $\R_\idot$, the
canonical map $\R^{red}_\idot \to \R_\idot$ induces equivalences
of categories
$$
\D(\wLR,\wt{\lambda}^*\R^{red}_\idot,k) \cong
\D(\wLR,\wt{\lambda}^*\R_\idot,k),\ 
\D(\LZ,\lambda^*\R^{red}_\idot,k) \cong
\D(\LZ,\lambda^*\R_\idot,k).
$$
\end{lemma}

\proof{} By \cite[Lemma 7.15 (ii)]{Ka-ma}, the map $\R^{red}_\idot
\to \R_\idot$ is a strong quasiisomorphism in the sense of
Definition~\ref{str.qua}; the claim then follows from
Proposition~\ref{redu.prop}.
\endproof

We now observe that a reduced $\LR$-graded $A_\infty$-coalgebra is
essentially a linear object: all the potentially non-linear
comultiplication maps are $0$ be definition. To make this precise,
let
$$
\LI_{red} = \coprod_{p \text{ prime}}\Lambda_p \subset \LI,
$$
where $\LI$ is as in \eqref{li.def}.

\begin{defn}
An $A_\infty$-functor from a small category $\C$ to some abelian
category $\Ab$ is {\em normalized} if for any $n$-tuple
$f_1,\dots,f_n$ of composable invertible maps $f_1,\dots,f_n$ in
$\C$, $n \geq 3$, the corresponding $A_\infty$ operation $b_n$ is
equal to $0$.
\end{defn}

\begin{lemma}\label{red.eq}
The category of reduced normalized $\LR$-graded
$A_\infty$-coal\-ge\-b\-ras and $A_\infty$-maps between them is
equivalent to the category of normalized $A_\infty$-functors from
$\LI_{red}^{opp}$ to $\Z\amod$ and $A_\infty$-maps between them.
\end{lemma}

\proof{} Assume given a reduced normalized $\LR$-graded
$A_\infty$-coalgebra $\R_\idot$. For any object $[m] \in \LR$ we
have the slice category $\LR_h/[m]$ of horizontal maps $f:[n] \to
[m]$, $[n] \in \LR$. Since $\R_\idot$ is normalized, for any map
$g:[m] \to [m']$ we can define a functor $\R^g_\idot$ from
$(\LR_h/[m])^{opp}$ to complexes of abelian groups by setting
$$
\R^g_\idot(f) = \R_\idot(g \circ f).
$$
This functor is constant (all transition maps are isomorphisms). Let
$$
P_\idot(g) = \lim_{\gets}\R^g_\idot(f),
$$
where the limit is taken over the category
$(\LR_h/[m])^{opp}$. Moreover, any horizontal map $h:[m] \to [m']$
induces a functor $\LR_h/[m] \to \LR_h/[m']$, $f \mapsto h \circ f$;
for any map $g:[m'] \to [m'']$, restriction with respect to this
functors gives a natural map
$$
h^*:P_\idot(g) \to P_\idot(g \circ h),
$$
and this is associative in the obvious sense.

Recall that we have the equivalence \eqref{li}, and restrict it to
$\LI_{red} \subset \LI$. For every object $a \in \LI_{red}$, let
$$
P_\idot(a) = P_\idot(v(a)).
$$
Then $P_\idot(-)$ has a natural structure of a normalized
$A_\infty$-functor from $\LI_{red}^{opp}$ to complexes of abelian
groups: for every $n$-tuple
$$
\begin{CD}
i_p(a_0) @>{f_1}>> i_p(a_1) @>{f_2}>> \dots @>{f_n}>> i_p(a_n)\\
@V{v_0}VV @V{v_1}VV @. @VV{v_n}V\\
\pi_p(a_0) @>{f'_1}>> \pi_p(a_1) @>{f'_2}>> \dots @>{f'_n}>> \pi_p(a_n)
\end{CD}
$$
of maps of the form \eqref{lzred.maps}, the $A_\infty$ operation
$$
P_\idot(a_n) \to P_\idot(a_0)
$$
is the composition of the map
$$
(f_n \circ \dots \circ f_1)^*:P_\idot(a_n) \to P_\idot(v_n \circ f_n
\circ \dots \circ f_1) = P_\idot(f'_n \circ \dots \circ f'_1 \circ
v_0)
$$
and the map
$$
P_\idot(f'_n \circ \dots \circ f'_1 \circ v_0) \to P_\idot(v_0) =
P_\idot(a_0)
$$
induced by $A_\infty$-operation on $\R_\idot$.

Conversely, assume given a normalized $A_\infty$-functor $P_\idot$
from $\LI_{red}^{opp}$ to complexes of abelian groups. For any map
$f:[m] \to [n]$ of prime degree in $\LR$, let $\C(f)$ be the
category of diagrams
$$
\begin{CD}
[m] @>{h}>> [m'] @>{v}>> [n]
\end{CD}
$$
with vertical $v$, horizontal $h$, and $f = v \circ h$. Then sending
such a diagram to $P_\idot(v)$ gives a functor $P_\idot(f)$ from
$\C(f)$ to complexes of abelian groups, and we set
$$
\R_\idot(f) = \lim_{\gets}P_\idot(f),
$$
where the limit is taken over $\C_f$. We leave it to the reader to
check that the $A_\infty$-functor structure on $P_\idot$ induces a
reduced normalized $A_\infty$-coalgebra structure on $\R_\idot(-)$,
and that both constructions are mutually inverse.
\endproof

Now we can make our final reduction. Say that a normalized
$A_\infty$-functor $M_\idot$ from $\LI_{red}^{opp}$ to $\Z\amod$ is
{\em admissible} if for any object $[a] \in \Lambda_p \subset
\LI_{red}$, we have
$$
M_i([a]) =
\begin{cases}
0, &\quad i < 0,\\
\Z, &\quad i=0,
\end{cases}
$$
and $M_i([a])$ is a free $\Z[\Z/pZ]$-module for $i \geq 1$. We will
say that a reduced normalized $\LR$-graded $A_\infty$-coalgebra
$\R_\idot$ is {\em admissible} if so is the $A_\infty$ functor
$E_\idot$ corresponding to $\R_\idot$ under the equivalence of
Lemma~\ref{red.eq}. Note that the $A_\infty$-coalgebra
$\T_\idot$ is admissible by \cite[Proposition 7.8]{Ka-ma}.

\begin{lemma}\label{adm.equi}
For any two admissible reduced normalized $\LR$-graded
$A_\infty$-coalgebras $\R_\idot$, $\R'_\idot$, and ring $k$, we have
canonical equivalences
$$
\begin{aligned}
\D(\wLR,\wt{\lambda}^*\R'_\idot,k) &\cong
\D(\wLR,\wt{\lambda}^*\R_\idot,k), \\
\D(\LZ,\lambda^*\R'_\idot,k) &\cong
\D(\LZ,\lambda^*\R_\idot,k).
\end{aligned}
$$
\end{lemma}

\proof{} Choose a projective resolution $\wt{P}_\idot$ of the
constant functor $\Z \in \Fun(\LI_{red}^{opp},\Z)$, and let $P_0 =
\Z$, $P_i = \wt{P}_{i-1}$ for $i \geq 1$. Then $P_\idot$ is an
admissible normalized $A_\infty$-functor from $\LI_{red}^{opp}$ to
$\Z\amod$, and $\wt{P}_\idot$ is $h$-projective by
Lemma~\ref{triv.coa}. Therefore for any other admissible normalized
$A_\infty$-functor $P'_\idot$, we have an $A_\infty$-map
$$
\xi:P_\idot \to P'_\idot.
$$
Moreover, for every $[m] \in \Lambda_p \subset \LI_{red}$, the map
$\xi([m])$ is a strong quasiisomorphism with respect to $\Z/p\Z
\subset \Aut([m])$. Let $P'_\idot$ be the $A_\infty$-functor
corresponding to the $A_\infty$-coalgebra $R_\idot$, and let
$\wh{\R}_\idot$ be the $A_\infty$-coalgebra corresponding to the
$A_\infty$-functor $P_\idot$. Then $\xi$ induces a map
$$
\xi:\wh{\R}_\idot \to \R_\idot,
$$
and this map is a strong quasiisomorphism in the sense of
Definition~\ref{str.qua}. Therefore by Proposition~\ref{redu.prop},
the corestriction functors corresponding to $\wt{\lambda}^*\xi$,
resp. $\lambda^*\xi$ are equivalences of categories. Analogously for
$\R'_\idot$.
\endproof

\subsection{Comodules.}\label{como.subs}

Consider now the category $\wLR_h \cong \Lambda \times I$ and the
products $\LI_{red} \times \N^* \subset \LI_{red} \times I$. Define
the functors
$$
i,\pi:\LI_{red} \times I \to \wLR_h
$$ 
by the same formula as in \eqref{i.p.red}. Moreover, let
$\tau:\LI_{red} \times I \to \LI_{red}$ be the tautological
projection.

Any complex $P_\idot$ of functors in $\Fun(\LI_{red}^{opp},\Z)$ is
in particular a normalized $A_\infty$-functor from $\LI_{red}^{opp}$
to $\Z\amod$. Fix such a complex $P_\idot$ so that it is admissible
in the sense of Lemma~\ref{adm.equi}. For any commutative ring $k$,
consider the category of pairs $\langle V_\idot,\phi \rangle$ of
complexes $V_\idot$ in $\Fun(\wLR_h^{opp},k)$ equipped with a map
$$
\phi:\pi^*V_\idot \to i^*V_\idot \otimes \tau^*P_\idot.
$$
Inverting quasiisomorphisms in this category, we obtain a
triangulated category denoted by
$\D(\wLR_h,P_\idot,k)$. Forgetting the map $\phi$ gives a
functor
$$
\wt{h}^*:\D(\wLR_h,P_\idot,k) \to \D(\wLR_h^{opp},k).
$$
Analogously, let $\D(\LZ,P_\idot,k)$ be the category of pairs
$\langle V_\idot,\phi \rangle$ of complexes $V_\idot$ in
$\Fun(\LZ_h^{opp},k)$ equipped with a map
$$
\phi:\pi^*V_\idot \to i^*V_\idot \otimes \tau^*P_\idot,
$$
with inverted quasiisomorphisms, where $i$ and $\pi$ are as in
\eqref{i.p.red}. Restricting from $\wLR_h$ to $\LZ_h$ gives a
forgetful functor
\begin{equation}\label{forg}
\D(\wLR_h,P_\idot,k) \to \D(\LZ_h,P_\idot,k),
\end{equation}
and we have the restriction functor
$$
h^*:\D(\LZ_h,P_\idot,k) \to \D(\LZ_h^{opp},k).
$$
Let $\R_\idot$ be the reduced normalized $\wLR$-graded
$A_\infty$-coalgebra corresponding to $P_\idot$ under the
equivalence of Lemma~\ref{red.eq}, and consider the derived
categories $\D(\LZ,\wt{\lambda}^*\R_\idot,k)$,
$\D(\wLR,\wt{\lambda}^*\R_\idot,k)$.

\begin{lemma}\label{p.eq}
There exists a canonical equivalences of categories
$$
\D(\LZ_h,P_\idot,k) \cong \D(\LZ,\lambda^*\R_\idot,k),\quad
\D(\wLR_h,P_\idot,k) \cong \D(\wLR,\wt{\lambda}^*\R_\idot,k)
$$
commuting with the restriction functors $h^*$, resp. $\wt{h}^*$..
\end{lemma}

\proof{} By Lemma~\ref{triv.coa}, we may modify the definition of
the derived category $\D(\wLR_h,P_\idot,k)$ by replacing complexes
in $\Fun(\wLR_h^{opp},k)$ with $A_\infty$-comodules over the trivial
$\wLR^{opp}_h$-graded $A_\infty$-coalgebra, and analogously for
$\LZ$. Then the resulting category of complexes is tautologically
equivalent the category of $A_\infty$-comodules over
$\wt{\lambda}^*\R_\idot$, resp. $\lambda^*\R_\idot$, and the
equivalence commutes with $\wt{h}^*$, resp. $h^*$. In particular, it
preserves quasiisomorphisms, hence descends to the derived
categories.
\endproof

Let now $\D_c(\wLR_h,P_\idot,k) \subset \D(\wLR_h,P_\idot,k)$ be the
full subcategory spanned by such $M \in \D(\wLR_h,P_\idot,k)$ that
$h^*M \in \D(\wLR_h^{opp},k)$ is locally constant in the sense of
Definition~\ref{norm.defn}, and let $\D_c(\LZ_h,P_\idot,k) \subset
\D(\LZ_h,P_\idot,k)$ be the full subcategory spanned by $M \in
\D(\LZ_h,P_\idot,k)$ with locally constant $h^*M$. Then combining
Lemma~\ref{p.eq} with Lemma~\ref{redu.equi} and
Lemma~\ref{adm.equi}, we obtain the following effective description
of the categories $\DML(k)$, $\DLR(k)$ of $k$-vallued cyclic Mackey
functors and $k$-valued cyclotomic complexes.

\begin{prop}\label{p.all}
For any admissible complex $P_\idot$ of functors from $\LI_{red}$ to
$\Z\amod$, there exist canonical equivalences
$$
\DML(k) \cong \D_c(\LZ_h,P_\idot,k), \qquad
\DLR(k) \cong \D_c(\wLR_h,P_\idot,k)
$$
of triangulated categories.\endproof
\end{prop}

Moreover, restricting from $\wLR$ to $\LZ$ as in \eqref{forg}, we
obtain restriction functors
$$
h^*:\D(\wLR^{opp},k) \to \D(\LZ^{opp},k), \qquad
h^*:\D(\wLR_h,P_\idot,k) \to \D(\LZ^{opp},k).
$$
Let $\D_w(\wLR^{opp}_h,k) \subset \D(\wLR^{opp}_h,k)$,
$\D_w(\wLR_h,P_\idot,k) \subset \D(\wLR_h,P_\idot,k)$ be the full
subcategories spanned by such $M$ that $h^*M \in \D(\LZ^{opp},k)$ is
locally constant.

\begin{prop}\label{reso.bis}
\begin{enumerate}
\item The restriction functor $\wt{h}^*$ has a right-adjoint
$$
\wt{h}_*:\D(\wLR_h^{opp},k) \to \D(\wLR_h,P_\idot,k),
$$
and it sends the subcategory $\D_w(\wLR_h^{opp},k) \subset
\D(\wLR_h^{opp},k)$ into the subcategory $\D_w(\wLR_h,P_\idot,k)
\subset \D(\wLR_h,P_\idot,k)$.
\item Assume given a triangulated subcategory $\D' \subset
\D_w(\wLR_h,P_\idot,k)$ which is closed with respects to arbitrary
products and contains all the objects $h_*M$, $M \in
\D(\LZ^{opp}_h,k)$. Then $\D'=\D_w(\wLR_h,P_\idot,k)$.
\end{enumerate}
\end{prop}

\proof{} On the level of categories of complexes, the forgetful
functor $\wt{h}^*$ has an obvious adjoint given by
$$
\wt{h}_*V_\idot = \pi_*(i^*V_\idot \otimes \tau^*P_\idot),
$$
with the tautological map $\phi$. Since every complex in
$\Fun(\LR_h^{opp},k)$ has an $h$-injective replacement, this
descends to the derived categories. 

\begin{lemma}\label{reso.bis.lemma}
For any $l \geq 0$, let $F^lP_\idot \subset P_\idot$ be the $l$-term
of the stupid filtration on $P_\idot$ (that is, $(F^lP_\idot)_m =
P_m$ for $m \leq l$, and $0$ otherwise). Then for every
$h$-injective complex $V_\idot$ in $\Fun(\wLR_h^{opp},k)$, the
natural map
\begin{equation}\label{lim.l}
\pi_*(i^*V_\idot \otimes \tau^*P_\idot) \to
\lim_{\overset{l}{\to}}R^\hdot\pi_*(i^*V_\idot \otimes
\tau^*F^lP_\idot)
\end{equation}
is a quasiisomorphism.
\end{lemma}

\proof{} We obviously have
$$
i^*V_\idot \otimes \tau^*P_\idot \cong
\lim_{\overset{l}{\to}}i^*V_\idot \otimes \tau^*F^lP_\idot,
$$
so that it suffices to prove that the map
$$
\pi_*(i^*V_\idot \otimes \tau^*F^lP_\idot) \to
R^\hdot\pi_*(i^*V_\idot \otimes \tau^*F^lP_\idot)
$$
is a quasiisomorphism for any $l \geq 0$. By induction, this reduces
to proving that for any integers $m$, $l \geq 0$, $n \geq 1$, we have
$$
R^n\pi_*(i^*V_m \otimes \tau^*P_l) = 0.
$$
This can be checked after evaluating at any object $a \in \wLR_h$;
by base change, we have to show that
$$
H^n(\Z/p\Z,i^*V_m(b) \otimes P_l(\tau(b))) = 0
$$
for every $b \in \Lambda_p \times I \subset \LI_{red} \times I$. But
since $V_\idot$ is $h$-injective, $V_\idot(c)$ is an $h$-injective
complex of $\Z[\Aut(c)]$-modules for any $c \in \wLR_h$, so that
$i^*V_m(b)$ is a free $\Z/p\Z$-module. Hence so is the product
$i^*V_m(b) \otimes \tau^*P_l(b)$.
\endproof

Now, by Proposition~\ref{p.all}, the category $\D(\wLR_h,P_\idot,k)$
does not depend on the choice of an admissible complex
$P_\idot$. Choose $P_\idot$ so that $F^lP_\idot$ is quasiisomorphic
to the shift $\Z[l]$ of the constant functor $\Z \in
\Fun(\LI_{red}^{opp},k)$ for every even integer $l \geq 0$. Let
$\overline{\pi}:\LI_{red} \times \N^* \to \LZ_h$ be the restriction
of the functor $\pi$ to $\LI_{red} \times \N^*$. Then $h^* \circ
R^\hdot \pi_* \cong R^\hdot\overline{\pi}_* \circ h^*$ by base
change, and
$$
R^\hdot\overline{\pi}_*:\D(\LI_{red}^{opp},k) \to \D(\LZ^{opp}_h,k)
$$
obviously sends $\D_c(\LI_{red}^{opp},k)$ into
$\D_c(\LZ^{opp}_h,k)$. Thus restricting to even $l$ in the
right-hand side of \eqref{lim.l}, we see that for every $V_\idot \in
\D_w(\wLR_h^{opp},k)$,
$$
\wt{h}^*\wt{h}_*V_\idot \cong
\lim_{\overset{l}{\to}}R^\hdot\pi_*i^*V_\idot[2l]
$$
indeed lies in $\D_w(\wLR_h^{opp},k)$. This finishes the proof of
\thetag{i}. \thetag{ii} now follows by exactly the same argument as
in Lemma~\ref{reso}.
\endproof

\section{Equivariant homology.}

\subsection{Generalities on equivariant homotopy.}\label{spectra.subs}

To fix notation, we start by recalling some general fact from
equivariant stable homotopy theory. The standard reference here is
\cite{LMS}; we mostly follow the exposition in \cite{HM} which
contains in a concise form everything we will need.

Let $G$ be a compact Lie group. A {\em $G$-CW complex} $X$ is a
pointed CW complex $X$ equipped with a continuous
fixed-point-preserving $G$-action such that for any $g \in G$, the
fixed-points subset $X^g \subset X$ is a subcomplex. Consider the
category of pointed $G$-toplogical spaces and $G$-equivariant maps
between them considered up to a $G$-equivariant homotopy, and let
$G\Top$ be the full subcategory spanned by spaces homotopy
equivalent to $G$-CW complexes. We note that for every closed Lie
subgroup $H \subset G$, sending $X$ to the fixed points subset $X^H
\subset X$ gives a well-defined functor
$$
G\Top \to W_H\Top,
$$
where $W_H = N_H/H$, and $N_H \subset G$ is the normalizer of the
subgroup $H$.

Given a finite-dimensional representation $V$ of the group $G$ over
$\RR$, we will denote by $S^V$ the one-point compactification of
$V$, with infinity being the fixed point. For any $X \in G\Top$, we
will denote $\Sigma^VX = X \wedge S^V$, and we will denote by
$\Omega^VX$ the space of based continuous maps from $S^V$ to
$X$. The functors $\Sigma^V,\Omega^V:G\Top \to G\Top$ are obviously
adjoint.

A {\em $G$-universe} is an $\RR$-vector space $U$ equipped with a
continous linear $G$-action and a $G$-invariant positive-definite
inner product. For any $G$-universe $U$, a {\em $G$-prespectrum $X$
indexed on $U$} is a collection of $G$-CW complexes $X(V)$, one for
each finite-dimensional $G$-invariant subspace $V \subset U$, and
$G$-equivariant continuous maps
\begin{equation}\label{trans.sp}
X(V) \to \Omega^WX(V \oplus W),
\end{equation}
one for every pair of transversal mutually orthogonal
finite-dimensional $G$-invariant subspaces $V,W \subset U$, subject
to an obvious associativity condition. The category of $G$-spectra
indexed on $U$ and homotopy classes of maps between them will be
denoted $G\spp(U)$.

A $G$-prespectrum is a {\em spectrum} if the maps \eqref{trans.sp}
are homeomorphisms. The category of $G$-spectra indexed on $U$ and
homotopy classes of maps between them will be denoted $G\Sp(U)$. We
have the tautological embedding $G\Sp(U) \to G\spp(U)$, and it
admits a left-adjoint {\em spectrification functor} $L$ given by
\begin{equation}\label{spfy}
Lt(V) = \dlim_{\overset{W}{\to}} \Omega^W\Sigma^WX(V \oplus W),
\end{equation}
where the limit is taken over all finite dimensional $G$-invariant
subspaces $W \subset U$ orthogonal to $V$.

For any inclusion $u:U_1 \subset U_2$ of $G$-universes, we have an
obvious restriction functor $\rho^{\#}(u):G\spp(U_2) \to G\spp(U_1)$,
called {\em change of universe}, and this functor has a
left-adjoint $\rho_{\#}(u):G\spp(U_1) \to G\spp(U_2)$ given by
$$
\wt{\rho}_{\#}(u)(t)(V) = \Sigma^{V - (V \cap u(U_1))}t(u^{-1}(V)),
$$
where $V - (V \cap u(U_1)) \subset V$ is the orthogonal to the
intersection $V \cap u(U_1) \subset V$.  The functor $\rho^{\#}(u)$
sends spectra to spectra; the corresponding functor
$\rho^{\#}(u):G\Sp(U_2) \to G\Sp(U_1)$ has a left-adjoint
$\rho_{\#}:G\Sp(U_1) \to G\Sp(U_2)$ given by
$$
\rho_{\#}(u) = L\wt{\rho}_{\#}(u),
$$
where $L$ is the spectrification functor \eqref{spfy}.

In particular, spectra indexed on a trivial universe $U=0$ are just
$G$-spaces, and the restriction $G\Sp \to G\Top$ with respect to the
embedding $0 \hookrightarrow U$ is the forgetful functor sending a
$G$-spectrum to its value at $0$. Its right-adjoint is called the
{\em suspension spectrum functor} and denoted $\Sigma^\infty:G\Top
\to G\Sp(U)$. Explicitly, $\Sigma^\infty X = L
\wt{\Sigma}^\infty X$, where $\wt{\Sigma}^\infty X \in G\spp(U)$ is
given by
$$
\wt{\Sigma}^\infty X(V) = \Sigma^VX.
$$
Of the non-trivial $G$-universes, those of two particular types are
important.
\begin{enumerate}
\item $U=\RR^\infty$ with the trivial $G$-action (where $\infty$ is
  assumed to be countable). Then $G$-spectra indexed on $U$ are
  called {\em naive} equivariant $G$-spectra; we denote the category
  $G\Sp(U)$ by $G\Sp^{naive}$.
\item A $G$-universe is {\em complete} if every finite-dimensional
  representation $V$ of the compact Lie group $G$ appears in $U$ a
  countable number of times. All complete $G$-universes are
  obviously isomorphic. A $G$-spectrum indexed on a complete
  $G$-universe $U$ is known as a {\em genuine} equivariant
  $G$-spectrum; we will denote $G\Sp(U)$ simply by $G\Sp$.
\end{enumerate}
We note that both $G\Sp^{naive}$ and $G\Sp$ are triangulated
categories, with shifts given $X[n](W) = X(W \oplus \RR^n)$. In
addition, $G\Sp$ has an autoequivalence $\Sigma^V:G\Sp \to G\Sp$ for
every finite-dimensional representation $V$, given by $\Sigma^VX(W)
= X(W \oplus V)$ (to make this precise, one has to fix an
isomorphism $U \oplus V \cong U$, or $U \oplus \RR^n \cong U$ in the
naive case, and apply the change of universes functor).

Assume given a closed Lie subgroup $H \subset G$. Then for any
$G$-universe $U$, $U^H$ is a $W_H$-universe, complete if $U$ was
complete. For any $G$-spectrum $X \in G\Sp(U)$, the {\em Lewis-May
fixed points spectrum} $X^H \in W_H\Sp(U^H)$ is given by
$$
X^H(V) = X(V)^H
$$
for any finite-dimensional $V \subset U^H \subset U$. There is a
second fixed points functor $\Phi^H:G\Sp(U) \to H\Sp(U^H)$ called
the {\em geometric fixed points functor}. To define it, one chooses
a finite-dimensional $G$-invariant subspace $W(V) \subset U$ for
every finite-dimensional $W_H$-invariant $V \subset U^H$ such that
$V = W(V)^H$ and
$$
\bigcup_{V \subset U^H}W(V) = U,
$$
and sets
\begin{equation}\label{sp.phi.def}
\phi^Ht(V) = t(W(V))^H
\end{equation}
for any $t \in G\spp(U)$, and
$$
\Phi^HX = L\phi^HX
$$
for any $G$-spectrum $X \in G\Sp(U)$. Here $\phi^H$ depends of the
choice of the subspaces $W(V)$, but the spectrification $\Phi^H$
does not (for a more invariant description of the functor $\Phi^H$
which make this explicit, see \cite[Lemma 1.1]{HM}). For any $X \in
G\Sp$, there is a natural map
$$
X^H \to \Phi^HX,
$$
and this map is functorial in $X$.

For naive $G$-spectra, both fixed points functors obviously
coincide. In the genuine case, let us fix a complete $G$-universe
$U$. Then $U^G \subset U$ is isomorphic to $\RR^\infty$, so that
$G\Sp(U^G)$ is $G\Sp^{naive}$, and the inclusion $u:U^G \to U$
induces a pair of adjoint functors $\rho^{\#}(u):G\Sp \to
G\Sp^{naive}$, $\rho_{\#}(u):G\Sp^{naive} \to G\Sp$. We then have
commutative diagrams
$$
\begin{CD}
G\Top @>{\Sigma^\infty}>> G\Sp^{naive} @>{\rho_{\#}(u)}>> G\Sp\\
@V{(-)^H}VV @V{(-)^H}VV @VV{\Phi^H}V\\
W_H\Top @>{\Sigma^\infty}>> W_H\Sp^{naive} @>{\rho_{\#}(u')}>> W_H\Sp
\end{CD}
$$
and
\begin{equation}\label{psi.unive}
\begin{CD}
G\Sp^{naive} @<{\rho^{\#}(u)}<< G\Sp\\
@V{(-)^H}VV @VV{(-)^H}V\\
W_H\Sp^{naive} @<{\rho^{\#}(u')}<< W_H\Sp,
\end{CD}
\end{equation}
where $u'$ is the embedding $U^G \subset U^H$ (and the
$W_H$-universe $U^H$ is obviously complete).

\subsection{Cyclic sets.}\label{cycl.subs}

From now on, we let $G = S^1 = U(1)$, the unit circle. Then it is
well-known that $G$-spaces are related to cyclic sets. Let us recall
the relation (for details and references, see e.g. \cite{Lo}).

For any object $[n] \in \Lambda$, let $|[n]|$ be its topological
realization: the union of points numbered by vertices $v \in V([n])$
and open intervals $I_e$ numbered by edges $e \in E([n])$, with the
natural topology making $|[n]|$ into a circle. For any function
$a:E([n]) \to \RR$ such that $a > 0$, we can make $|[n]|$ into a
metric space $|[n]|(a)$ by assigning length $a(e)$ to the interval
$I_e$. Let $\Rr([n])^o$ be the space of pairs $\langle a,b \rangle$
of a function $a:E([n]) \to \RR$, $a > 0$, and a metric-preserving
monotonous continuous map $b:|[n]|(a) \to S^1$ to the unit circle
$S^1 \subset \S$. Such a map $b$ exists if and only if $a_1 + \dots
+ a_n = 2\pi$, and the space of all such maps is non-canonically
identified with $S^1$, so that we have a non-canonical homeomorphism
$$
\Rr([n])^o \cong S^1 \times \Tt_{n-1}^o,
$$
where for any $m \geq 0$, $\Tt_m^o \subset \Tt_m$ is the interior of
the standard $m$-simplex $\Tt_m$. Embed $\Rr([n])^o$ into a compact
space $\Rr([n])$ by allowing $a$ to take zero values (and degenerate
metrics on $|[n]|$). We then have
$$
\Rr([n]) \cong S^1 \times \Tt_{n-1}.
$$
This decomposition is not canonical, but the space $\Rr([n])$ itself
is completely canonical, and by construction, it carries a continuos
$G$-action. 

Moreover, for any map $f:[n] \to [m]$ and a pair $\langle a,b
\rangle \in \Rr([m])$, let
$$
f(a)(e) = \sum_{f(e')=e}s(e'), 
$$
so that we have an obvious metric-preserving map $|[n]|(f(a)) \to
|[m]|(a)$, and let $f(b)$ be the composition of this map with
$b$. This makes $\Rr$ into a contravariant functor from $\Lambda$ to
$G\Top$. Turn it into a covariant functor by applying the duality
$\Lambda^{opp} \cong \Lambda$. The result is a functor
\begin{equation}\label{rr.eq}
\Rr:\Lambda \to G\Top.
\end{equation}
By the standard Kan extension procedure, it extends uniquely to a
  colimit-preserving {\em realization functor}
$$
\Real:\Lambda^{opp}\Sets \to G\Top
$$
such that $\Real \circ Y \cong \Rr$, where $Y:\Lambda \to
\Lambda^{opp}\Sets$ is the Yoneda embedding. We also the
right-adjoint functor
$$
S:G\Top \to \Lambda^{opp}\Sets
$$
such that
$$
S(X)([n]) = \Maps_G(\Rr([n]),X),
$$
for any $X \in G\Top$, where $\Maps_G(-,-)$ stands for the set of
$G$-equivariant unbased maps.

It is well known that for any $A \in \Lambda^{opp}\Sets$, the
realization $\Real(A)$ is homeomorphic to the usual geometric
realization of the simplicial set $j^*A \in \Delta^{opp}\Sets$, and
for any $X \in G\Top$, the adjunction map
$$
\Real(S(X)) \to X
$$
is a homotopy equivalence for every $X \in G\Top$ (see \cite{Lo}, or
\cite{Dr} for a modern treatment). In particular, if for any set $A$
we denote by $\Z[A]$ the free abelian group spanned by $A$, and
apply this pointwise to $S(X):\Lambda \to \Sets$, then we have
\begin{equation}\label{equi.eq}
H_\idot(\Lambda^{opp},\Z[S(X)]) \cong H_\idot(X_{hG},\Z),
\end{equation}
where $X_{hG}$ is the homotopy quotient of $X$ by the $G$-action.

However, for our applications the functors $\Z[S(X)] \in
\Fun(\Lambda^{opp},\Z)$ are incovenient because they are not locally
constant in the sense of Definition~\ref{norm.defn}. To remedy this,
we note that the sets $S(X)([n])$ all carry a natural topology.

\begin{defn}\label{cycl.chain}
For any $X \in G\Top$, its {\em chain complex} $C_\idot(X)$ is a
complex in $\Fun(\Lambda^{opp},\Z)$ given by
$$
C_\idot(X)([n]) = C_\idot(S(X)([n]),\Z),
$$
where $C_\idot(-,\Z)$ is the normalized singular chain homology
complex with coefficients in $\Z$.
\end{defn}

Since all the maps $\Rr([n]) \to \Rr([n'])$ are homotopy
equivalences, $C_\idot(X)$, unlike $\Z[S(X)]$, is locally constant
in the sense of Definition~\ref{norm.defn}.  More explicitly,
$C_\idot(X)$ is obtained as follows. We consider the functor
$\wt{\Rr}:\Lambda \times \Delta \to G\Top$ given by
$$
\wt{\Rr}([n] \times [m]) = \Rr([n]) \times \Tt_{m-1},
$$
where $[m] \in \Delta$ is the totally ordered set with $m$
elements. The functor $\wt{\Rr}$ also extends to a pair of adjoint
functors
\begin{gather*}
\wt{\Real}:\Fun(\Lambda^{opp} \times \Delta^{opp},\Sets) \to G\Top,\\
\wt{S}:G\Top \to \Fun(\Lambda^{opp} \times \Delta^{opp},\Sets).
\end{gather*}
to obtain $C_\idot(S(X))$, we take $\Z[\wt{S}(X)] \in
\Fun(\Lambda^{opp} \times \Delta^{opp},\Z)$, and apply the
normalized chain complex construction fiberwise with respect to the
projection $\tau:\Lambda \times \Delta \to \Lambda$.

\medskip

To explain the relation between $\Z[S(X)]$ and $C_\idot(X)$, let
$\iota_m:\Lambda \to \Lambda \times \Delta$ be the embedding given
by $\iota([n]) = [n] \times [m]$, for any $m \geq 1$. Then $\iota_1$
is adjoint to $\tau$, so that for any $A \in \Fun(\Lambda^{opp}
\times \Delta^{opp},\Sets)$, we have a natural adjunction map
\begin{equation}\label{simpl.l}
\tau^*\iota_1^*A \to A.
\end{equation}
Apply this to $\wt{S}(X)$ for some $X \in G\Top$. By definition, we
have $\iota_1^*\wt{S}(X) \cong S(X)$; taking free abelian groups and
apssing to normalized chain complexes, we obtain a natural map
\begin{equation}\label{chain.full}
\Z[S(X)] \to C_\idot(X).
\end{equation}

\begin{lemma}\label{fiberw}
Assme given $A \in \Fun(\Lambda^{opp} \times \Delta^{opp},\Sets)$
such that for any map $f:[m] \to [n]$ in $\Delta$, the corresponding
map
$$
\iota_n^*A \to \iota_m^*A
$$
induces a homotopy equivalence of geometric realizations. Then the
map 
$$
\Real(\iota^*A) \cong \wt{\Real}(\tau^*\iota^*A) \to \wt{\Real}(A)
$$
induced by \eqref{simpl.l} is a homotopy equivalence.
\end{lemma}

\proof{} Let
$$
\Real_{\Lambda}:\Fun(\Lambda^{opp} \times \Delta^{opp},\Sets) \to
\Delta^{opp}\Top
$$
be the fiberwise geometric realization functor with respect to the
projection $\Lambda \times \Delta \to \Delta$. Then the assumption
on $A$ means that \eqref{simpl.l} becomes a homotopy equivalence
already after applying $\Real_\Lambda$.
\endproof

For any $X \in G\Top$, we have $\iota^*_m\wt{S}(X) \cong
S(\Maps(\Delta_{m-1},X))$, so that $\wt{S}(X)$ automatically
satisfies the conditions of Lemma~\ref{fiberw}. Thus in particular,
the natural map $\Real(S(X)) \to \wt{\Real}(\wt{S}(X))$ is a
homotopy equivalence, and the map
\begin{equation}\label{chain.2}
H_\idot(\Lambda,\Z[S(X)]) \to H_\idot(\Lambda,C_\idot(X))
\end{equation}
induced by \eqref{chain.full} is an isomorphism. More generally, if
we denote by $C:\D(\Lambda,\Z) \to \D_c(\Lambda,\Z)$ the
left-adjoint functor to the embedding $\D_c(\Lambda,\Z) \subset
\D(\Lambda,\Z)$, then $C(\Z[S(X)]) \cong C_\idot(X)$.

\medskip

The complexes $C_\idot(X)$ are usually rather big. We will need
smaller models for the standard $G$-orbit spaces $G/C_n$, where $C_n
= \Z/n\Z \subset G$ is the group of $n$-th roots of unity. To obtain
such models, consider the functors $j_n:\Delta \to \Lambda_n^{opp}$
of Subsection~\ref{connes.subs}, with the natural $C_n$-actions on
them. Define functors $J_n:\Lambda^{opp} \times \Delta^{opp} \to
\Sets$ by setting
\begin{equation}\label{J.n}
J_n([m] \times [l]) = \Lambda(i_nj_n([l]),[m])/C_n,
\end{equation}
where $C_n$ acts through its action on $j_n$.

\begin{lemma}\label{J.n.le}
For any $n \geq 1$, we have a natural homotopy equivalence
$$
\wt{\Real}(J_n) \cong G/C_n.
$$
\end{lemma}

\proof{} By Lemma~\ref{fiberw}, it suffices to construct
equivalences
$$
\Real(\iota_l^*J_n) \cong G/C_n
$$
for any $[l] \in \Delta$. By definition, we have
$$
\iota^*_lJ_n \cong Y(j_n([l]))/C_n = Y([nl])/C_n,
$$
where $Y([nl]) \subset \Lambda^{opp}\Sets$ is the Yoneda image of
$[nl] \in \Lambda$. Since $\Real$ is colimit-preserving, we have
homeomorphisms
$$
\Real(Y([nl])/C_n) \cong \Real(Y([nl]))/C_n \cong \Rr([nl])/C_n,
$$
and $\Rr([nl])$ is indeed homotopy-equivalent to $G$.
\endproof

We will now introduce the fixed points subsets into the picture. By
Lemma~\ref{Y.le}, the functor $\Rr$ of \eqref{rr.eq} extends to a
functor
$$
\Rr = \Real \circ Y:\LZ \to G\Top,
$$
so that we obtain an adjoint pair of functors
$$
\Real:\Fun(\LZ^{opp},\Sets) \to G\Top, \quad
S:G\Top \to \Fun(\LZ^{opp},\Sets).
$$
We keep the same notation because the functors are direct extensions
of the functors $\Real$ and $S$ of Subsection~\ref{cycl.subs}.  In
effect, the restriction functor $\Lambda^{opp}\Sets \to
\Fun(\LZ^{opp},\Sets)$ admits a fully faithful left-adjoint
$$
L:\Lambda^{opp}\Sets \to \Fun(\LZ^{opp},\Sets),
$$
and we have $S \cong L \circ S$, $\Real \cong \Real \circ
L$. Explicitly, if we restrict $L(A)$ to $\Lambda \cong \LZ_h^m
\subset \LZ$, we have
\begin{equation}\label{fixed.eq}
L(A)|_{\LZ_h^m} = \pi_{m*}i_m^*A,
\end{equation}
where the adjoint $\pi_{m*}:\Lambda_m^{opp}\Sets \to
\Lambda^{opp}\Sets$ to the pullback functor $\pi_m^*$ is obtained by
taking $\Z/m\Z$-fixed points fiberwise, and the vertical maps in
$\LZ$ act by natural inclusions of fixed points.

As in Subsection~\ref{cycl.subs}, we extend the functors $S$ and
$\Real$ to the product category $\LZ \times \Delta$. For any $X \in
G\Top$, we define its {\em extended chain complex} $\wt{C}_\idot(X)$
of functors in $\Fun(\LZ^{opp},\Z)$ by setting
$$
\wt{C}_\idot(X) = C_\idot(S(X)) = N(\Z[\wt{S}(X)]),
$$
where $N$ is the normalized chain complex functor applied fiberwise
to the projection $\LZ \times \Delta \to \LZ$. Moreover, define the
{\em reduced chain complex} $\overline{C}_\idot(X)$ as
$$
\overline{C}_\idot(X)([n|m]) = \wt{C}_\idot(X)([n|m])/\Z[o],
$$
where $o \in \Maps_{S^1}(\Rr([m|n]),X)$ is the distinguished
point. Then the standard shuffle map induces a natural quasiisomorphism
\begin{equation}\label{shuffle}
\overline{C}_\idot(X) \otimes \overline{C}_\idot(Y) \to
\overline{C}_\idot(X \wedge Y)
\end{equation}
for any $X,Y \in G\Top$, where $X \wedge Y$ is the smash product,
and the tensor product in the right-hand side is the pointwise
tensor product in the category $\Fun(\LZ^{opp},\Z)$.

The reduced chain complex functor $\overline{C}_\idot(-)$ easily
extends to the category $G\Sp^{naive}$ of naive $G$-equivariant
spectra. Namely, for any $X \in G\Sp^{naive}$ and integers $i,j \geq
0$, the transition maps \eqref{trans.sp} induce by adjunction maps
$$
\Sigma^i X(\RR^{\oplus j}) \to X(\RR^{\oplus i+j}),
$$
which give rise to maps
\begin{equation}\label{trans.naive}
\overline{C}_\idot(X(\RR^{\oplus i}))[-i] \cong
\overline{C}_\idot(\Sigma^i X(\RR^{\oplus j})) \to
\overline{C}_\idot(X(\RR^{\oplus i+j})).
\end{equation}

\begin{defn}\label{c.naive.def}
The {\em equivariant homology complex} $C^{naive}_\idot(X)$ of a
naive $G$-spectrum $X \in G\Sp^{naive}$ is given by
$$
C^{naive}_\idot(X) =
  \dlim_{\overset{i}{\to}}\overline{C}_\idot(X(\RR^{\oplus i}))[i] \in
  \D(\LZ^{opp},\Z),
$$
where the limit is taken with respect to the maps
\eqref{trans.naive}.
\end{defn}

\subsection{Cyclic Mackey functors.}\label{chain.subs}

To extend the equivariant homology complex of
Definition~\ref{c.naive.def} to genuine $G$-spectra, we need to pass
to the category $\DML(\Z)$ of cyclic Mackey functors of
Section~\ref{mack.sec}. We will use the model of
Proposition~\ref{p.all}, with an appropriately chosen complex
$P_\idot$.

We begin with the following observation. For any $n \geq 1$, let
$P^n_\idot$ be the complex in $\Fun(\LZ^{opp},\Z)$ obtained as the
cone of the natural augmentation map
$$
N(L(J_n)) \to \Z,
$$
where $\Z \in \Fun(\LZ^{opp},\Z)$ is the constant functor, and
$J_n:\Lambda^{opp} \times \Delta^{opp} \to \Sets$ is as in
\eqref{J.n}. We note that for any $[m|l] \in \LZ$,
$P^n_\idot([m|l])$ is a finite-length complex of finitely generated
free abelian groups. Now, let $\CC(1)$ be the standard complex
representation of the group $G = U(1)$, and let $\CC(n) =
\CC(1)^{\otimes n}$. Treat $\CC(n)$ as a $2$-dimensional real
representation by restriction of scalars, and denote
$$
\Sigma_n = \Sigma^{\CC(n)}:G\Top \to G\Top.
$$

\begin{lemma}\label{susp}
For any $X \in G\Top$, we have a natural functorial quasiisomorphism
$$
P^n_\idot \otimes \overline{C}_\idot(X) \to
\overline{C}_\idot(\Sigma_nX).
$$
\end{lemma}

\proof{} We have $S^{\CC(n)} \cong \Sigma(G/C_n)$, the
non-equivariant suspension of the standard orbit $G/C_n$. The
homotopy equivalence of Lemma~\ref{J.n.le} induces by adjunction a
map
$$
J_n \to \wt{S}(G/C_n),
$$
and this map induces a quasiisomorphism
$$
P^n_\idot \to \overline{C}_\idot(S^{\CC(n)}).
$$
Combining this with \eqref{shuffle}, we get the claim.
\endproof

Now let us fix a complete $G$-universe $U$, by taking
\begin{equation}\label{unive}
U = \bigoplus_{n \geq 0}\CC(n)^\infty,
\end{equation}
where $\infty$ means the sum of a countable number of copies (this
is complete, since $\C(-n) \cong \C(n)$ as real representations).

Say that a map $\nu:\N \to \N$ is admissible if $\nu(i) \neq 0$ for
at most a finite number of $i$. For any admissible sequence $\nu$,
let
$$
V(\nu) = \bigoplus_i \CC(i)^{\oplus \nu(i)} \subset U,
$$
and let
$$
\begin{aligned}
\Sigma_{\nu} &= \Sigma^{V(\nu)} = \Sigma_{i_1})^{\nu(i_1)} \circ
\dots \circ (\Sigma_{i_n})^{\nu(i_n)},\\
P^{\nu}_\idot &= \bigotimes_{1 \leq l \leq
  n}\left(P^{i_l}_\idot\right)^{\otimes\nu(i_l)},
\end{aligned}
$$
where $i_1 \leq \dots \leq i_n$ are all integers such that $\nu(i_l)
\neq 0$. Then Lemma~\ref{susp} immediately gives canonical
quasiisomorphisms
\begin{equation}\label{trans.P}
P^{\nu}_\idot \otimes \overline{C}_\idot(X) \to
C_\idot(\Sigma_{\nu}(X)).
\end{equation}
We will call the subspaces $V(\nu) \subset U$ {\em cellular}; the
collection of all cellular subspaces is obviously cofinal in the
collection of all finite-dimensional $G$-invariant subspaces $V
\subset U$.

\medskip

We are now ready to define our complex $P_\idot$ in
$\Fun(\LI_{red}^{opp},\Z)$. Note that for any $m,l \geq 1$, the
natural $C_m$-action on $i_m^*J_n([l]) = J_n([ml])$ factors through
a free action of the quotient $C_m/(C_m \cap C_n)$. Therefore for
any vertical map $v:a \to b$ in $\LZ$, the corresponding map
$$
v:P^n_\idot(b) \to P^n_\idot(a)
$$
is an isomorphism if $\deg(v)$ divides $n$, and the natural
inclusion $\Z = P^n_0(b) \to P^n_\idot(a)$ otherwise. Moreover, for
any prime $p$ not dividing $n$, the restriction
$$
i_p^*(P^n_\idot)
$$
with respect to the functor $i_p:\Lambda_p \to \Lambda \cong \LZ_h^1
\subset \LZ$ is the constant functor $\Z$ in degree $0$, while
$i_p^*P^n_i([m])$ is a free $\Z/p\Z$-module for any $[m] \in
\Lambda_p$ and any $i \geq 1$. We now define a complex $P_\idot$ of
functors in $\Fun(\LI_{red}^{opp},\Z)$ by setting
$$
P_\idot = \dlim_{\to}i_p^*P^{\nu}_\idot \text{ on } \Lambda_p
\subset \LI_{red},
$$
where the limit with respect to the natural inclusions $\Z \to
P^n_\idot$ is taken over all admissible sequences $\nu$ such that
$\nu(pl) = 0$, $l \in \N$. Then $P_\idot$ is admissible in the sense
of Proposition~\ref{p.all}, so that we have a canonical equivalence
$$
\DML(\Z) \cong \D(\LZ_h,P_\idot,\Z).
$$

\begin{remark}
To help the reader visualize the complex $P_\idot$, we note that on
$\Lambda_p \subset \LI_{red}$, it is essentially given by
$$
i_p^*\lim_{\overset{V}{\to}}\overline{C}_\idot(S^V),
$$
where the limit is taken over all cellular subspaces in $U$
orthogonal to $U^{C_p} \subset U$. The only difference is in that we
use more economical simplicial models for the spheres.
\end{remark}

By definition, $P^{\nu}_\idot$ for any $n$ gives an object in the
category $\D(\LZ_h,P_\idot,\Z)$. The corresponding complex is
$$
h^*P^{\nu}_\idot,
$$
the restriction with respect to the embedding $\LZ_h \to \LZ$, and
the map $\phi:\pi^*h^*P^{\nu}_\idot \to i^*h^*P^{\nu}_\idot$ is
induced by the action of the vertical maps on
$P^{nu}_\idot$. Explicitly, denote by $\nu_p$ the map given by
$$
\nu_p(m) = 
\begin{cases}
\nu(m), &\quad m=pl, l \geq 1,\\
0, &\quad\text{otherwise}.
\end{cases}
$$
Then we have $P^{\nu}_\idot = P^{\nu_p}_\idot \otimes P^{\nu -
  \nu_p}_\idot$, and the $\Lambda_p$-component $\phi_p$ of the map
$\phi$ is the product of an isomorphism 
$$
\pi_p^*h^*P^{\nu_p}_\idot \cong i_p^*h^*P^{\nu_p}_\idot
$$
and the natural inclusion
$$
\Z \cong \pi_p^*h^*P^{\nu-\nu_p}_\idot \to
i_p^*h^*P^{\nu-\nu_p}_\idot.
$$
If we treat the complexes $P_\idot^n$ up to a quasiisomorphism, then
on $\LZ_h^m \subset \LZ_h$, the restriction $h^*P^n_\idot$ is given
by 
$$
\begin{cases}
\Z[2], &\quad n=ml, l \geq 1,\\
\Z, &\quad\text{otherwise}.
\end{cases}
$$
Let us now fix a projective resolution $Q_\idot$ of the constant
functor $\Z \in \Fun(\LZ^{opp}_h,\Z)$, and let us define a complex
$Q^n_\idot$ by setting
$$
Q^n_\idot =
\begin{cases}
\Z[-2], &\quad n=ml, l \geq 1,\\
\Z, &\quad\text{otherwise}
\end{cases}
$$
on $\LZ^m_h \subset \LZ_h$. Then $Q^n_\idot \otimes h^*P^n_\idot$ is
quasiisomorphic to the constant functor $\Z$, and since $Q_\idot$ is
projective, this can be realized by an actual quasiisomorphism
$$
\eps_n:Q_\idot \to Q^n_\idot \otimes h^*P^n_\idot.
$$
Moreover, for any admissible $\nu:\N \to \N$ and any $m \geq 1$, let
$$
d(m,\nu) = \sum_{l \geq 1}\nu(ml),
$$
and let $Q^{\nu}_\idot$ be a complex in $\Fun(\LZ^{opp}_h,\Z)$
given by $Q_\idot[-2d(m,\nu)]$ on $\LZ^m_h \subset \LZ^h$. Then
taking tensor product of the maps $\eps_n$, we obtain a system of
quasiisomorphisms 
\begin{equation}\label{trans.eps}
\eps_{\nu}:\Z \to Q^{\nu}_\idot \otimes h^*P^{\nu}_\idot
\end{equation}
for all admissible $\nu$. For any prime $p$, we have a natural
map
$$
\eps_{\nu-\nu_p}:\pi_p^*Q^{\nu}_\idot \to i_p^*(Q^{\nu}_\idot \otimes
P^{\nu-\nu_p}_\idot),
$$
and composing these maps with the natural embeddings
$i_p^*P^{\nu-\nu_p}_\idot \to P_\idot$, we equip $Q^{\nu}_\idot$
with a natural structure of an object in $\D(\LZ_h,P_\idot,\Z)$.

We are now ready to define our equivariant homology functor. Assume
given a prespectrum $t \in G\spp(U)$. For any two admissible $\nu,
\nu':\N \to \N$, we have a transition map
$$
\Sigma_{\nu'}t(V(\nu)) \to t(V(\nu + \nu'))
$$
adjoint to the map \eqref{trans.sp}. By \eqref{trans.P}, these maps
induce canonical maps
$$
P^{\nu'}_\idot \otimes \overline{C}_\idot(t(V(\nu))) \to
\overline{C}(V(\nu+\nu')).
$$
Let now $\xi^*:\D(\LZ^{opp},\Z) \to \D(\LZ,\lambda^*\T_\idot,\Z)
\cong \DML(\Z)$ be the corestriction functor with respect to the
augmentation map of the $A_\infty$-coalgebra
$\lambda^*\T_\idot$. Then these maps composed with the maps
\eqref{trans.eps} induce transition maps
\begin{equation}\label{trans.mack}
\begin{aligned}
\xi^*\overline{C}_\idot(t(V(\nu)) &\to
Q^{\nu}_\idot \otimes \xi^*P^{\nu}_\idot \otimes
\xi^*\overline{C}_\idot(t(V(\nu))) \to\\
&\to Q^{\nu}_\idot \otimes \xi^*\overline{C}_\idot(t(V(\nu + \nu')))
\end{aligned}
\end{equation}

\begin{defn}
The {\em equivariant chain complex} $C_\idot(X) \in \DML(\Z)$ of a
$G$-prespectrum $t \in G\spp(U)$ is given by
\begin{equation}\label{c.def.eq}
C_\idot(t) = \dlim_{\overset{\nu}{\to}}Q^{\nu}_\idot \otimes
\xi^*\overline{C}_\idot(X(V(\nu))),
\end{equation}
where the limit is taken over all the admissible maps $\nu:\N \to
\N$ with respect to the transition maps \eqref{trans.mack}.
\end{defn}

\begin{lemma}\label{sp.equi}
For any $t \in G\spp(U)$ with spectrification $Lt$, the adjunction
map $t \to Lt$ induces an isomorphism
$$
C_\idot(t) \to C_\idot(Lt).
$$
\end{lemma}

\proof{} Since cellular subspaces are cofinal, in it suffices to
take the limit over cellular subspaces in \eqref{spfy}. Then
substituting \eqref{spfy} into \eqref{c.def.eq}, we see that
$C_\idot(Lt)$ is given by
$$
\dlim_{\overset{\nu \leq \nu'}{\to}}Q^{\nu}_\idot \otimes
\xi^*\overline{C}_\idot(\Omega^{V(\nu'-\nu)}t(V(\nu'))),
$$
where the limit is taken over all pairs $\nu \leq \nu'$ of
admissible maps $\nu,\nu':\N \to \N$. The subset of all pairs with
$\nu = \nu'$ is cofinal in this set, so it suffices to take the
limit over such pairs. This is exactly $C_\idot(t)$.
\endproof

\begin{corr}\label{le.unive}
For any $X \in G\Top$, we have
$$
C_\idot(\Sigma^\infty X) \cong \xi^*\overline{C}(X),
$$
and for any $X \in G\Sp^{naive}$, we have
\begin{equation}\label{c.unive}
C_\idot(\rho_{\#}(u)X) \cong C_\idot(X),
\end{equation}
where $\rho_{\#}(u)$ is the change-of-universe functor associated to
the embedding $\RR^{\oplus\infty} =U^G \subset U$.
\end{corr}

\proof{} By Lemma~\ref{sp.equi}, we can replace $\Sigma^\infty$ with
$\wt{\Sigma}^\infty$ and $\rho_{\#}(u)$ with
$\wt{\rho}_{\#}(u)$. Then for $t = \wt{\Sigma}^\infty X$, all the
transition maps in the filtered limit of \eqref{c.def.eq} are
quasiisomorphisms, and for $t = \wt{\rho}_{\#}(u)X$, the only
possibly non-trivial transition maps are the those of
\eqref{trans.naive}.
\endproof

\subsection{Cyclotomic complexes.}\label{cyclo.subs}

Now for any $m \geq 1$, let $C_m = \Z/m\Z \subset G = U(1)$ be the
group of the $m$-th roots of unity. The $m$-power map gives an
isomorphism $p_m:G/C_m \to G$, and we have $p_m \circ p_n = p_{nm}$,
$m,n \geq 1$. We have obvious canonical $G$-equivariant
isomorphisms
$$
u_m:U^{C_m} \cong U,
$$
where $U$ is the complete $G$-universe \eqref{unive}, and we have
$u_m \circ u_n = u_{nm}$, $n,m \geq 1$.

\begin{defn}\label{ext.fixed}
For any $X \in G\Sp = G\Sp(U)$ and any $m \geq 1$, the {\em extended
geometric fixed points spectrum} $\wPhi^m(X) \in G\Sp$ is given
by
$$
\wPhi^m(X) = \rho_{\#}(u_m)\Phi^{C_m}(X),
$$
and the {\em extended Lewis-May fixed points spectrum} $\wPsi^m(X)
\in G\Sp$ is given by
$$
\wPsi^m(X) = \rho_{\#}(u_m)X^{C_m}.
$$
\end{defn}

For any $m \geq 1$, let $\rho_m:I \to I$ be the multiplication my
$m$, as in \eqref{i.p.red}. It commutes with the $\N^*$-action,
hence induces an endofunctor
$$
\wt{\rho}_m:\LZ \to \LZ
$$
commuting with the projection $\lambda:\LZ \to \LR$. By
\eqref{fixed.eq}, we have
$$
\wt{\rho}_m^*LA \cong LA^{C_m}
$$
for any $A \in \Lambda^{opp}\Sets$. In particular, for every $X \in
G\Top$ we have
\begin{equation}\label{fixed.top}
\wt{\rho}^*_m\overline{C}_\idot(X) \cong
\overline{C}_\idot(X^{C_m}),
\end{equation}
and by \eqref{trans.naive}, these quasiisomorphisms induce
quasiisomorphisms
\begin{equation}\label{fixed.naive}
\wt{\rho}^*_mC_\idot(X) \cong C_\idot(X^{C_m})
\end{equation}
for every naive $G$-spectrum $X \in G\Sp^{naive}$. We want to obtain
a version of this for genuine $G$-spectra. Since $\lambda\circ
\wt{\rho}_m =\lambda$, we tautologically have
$\wt{\rho}^*\lambda^*\T_\idot \cong \lambda^*\T_\idot$, so that we
have a pullback functor
$$
\wt{\rho}_m^*:\DML(Z) \to \DML(Z).
$$

\begin{lemma}\label{phi.chain}
For any $m \geq 1$ and any $X \in G\Sp$, we have a natural
functorial isomorphism
$$
C_\idot(\wPhi^m(X)) \cong \wt{\rho}_m^*C_\idot(X).
$$
\end{lemma}

\proof{} For any map $\nu:\Z \to \Z$, let $r(\nu):\Z \to \Z$ be
given by $r(\nu)(ma + b) = \nu(a)$, $a \in \Z$, $0 \leq b < n$. If
$\nu$ is admissible, then so is $r(\nu)$.

By Lemma~\ref{sp.equi}, we may replace $\wPhi^m$ with the functor
$\phi^m = \wt{\rho}_{\#}(u_m) \circ \phi^{C_m}$. Choose the
subspaces $W(V) \subset U$ so that $W(V(\nu)) = V(r(\nu))$ for any
admissible $\nu:\Z \to \Z$. Then $C_\idot(\phi^mX)$ is given by
\begin{equation}\label{l.1}
\dlim_{\overset{\nu}{\to}}\Q^{\nu}_\idot \otimes
\xi^*\overline{C}_\idot(X(V(r\nu))^{C_m}),
\end{equation}
and since the sequences $r(\nu)$ are cofinal in the set of all
admissible sequences, $\wt{\rho}_m^*C_\idot(X)$ is given by
\begin{equation}\label{l.2}
\dlim_{\overset{\nu}{\to}}\wt{\rho}_m^*Q^{r(\nu)}_\idot
\xi^*\overline{C}_\idot(X(V(r\nu))).
\end{equation}
By definition, the pullback functor commutes with corestriction,
$\wt{\rho}_m^* \circ \xi^* \cong \xi^* \circ \wt{\rho}^*_m$, and we
have the isomorphisms \eqref{fixed.top}. It remains to notice that
by definition, we have
$$
\wt{\rho}_m^*Q^{r(\nu)}_\idot \cong Q^{\nu}_\idot
$$
for any admissible $\nu$.
\endproof

We now consider the right-adjoint $\xi_*:\DML(\Z) \to
\D(\LZ^{opp},\Z)$ to the corestriction $\xi^*$, as in
Subsection~\ref{addi.subs}. The isomorphism \eqref{c.unive} than
induces a base change map
\begin{equation}\label{c.bc}
C_\idot(\rho^{\#}(u)(X)) \to \xi_*C_\idot(X)
\end{equation}
for any $X \in G\Sp$. By virtue of \eqref{fixed.naive} and
\eqref{psi.unive}, we have a natural isomorphism
$$
\wt{\rho}_m^*C_\idot(\rho^{\#}(u)(X)) \cong C_\idot(\wPsi^m(X))
$$
for any genuine $G$-spectrum $X \in G\Sp$ and any $m \geq 1$.

\begin{lemma}\label{psi.bc}
The base change map \eqref{c.bc} is an isomorphism for any $X \in
G\Sp$, so that
$$
C_\idot(\wPsi^m(X)) \cong \wt{\rho}_m^*\xi_*C_\idot(X)
$$
for any $m \geq 1$.
\end{lemma}

\proof{} It suffices to prove that the map \eqref{c.bc} becomes an
isomorphism after evaluating at any object $[n|m] \in \LZ$. By
\eqref{psi.unive} and the definition of the equivariant homology
complex $C_\idot(-)$ of a naive $G$-spectrum, we have a natural
quasiisomorphism
$$
C_\idot(\rho^{\#}(u)(X))([n|m]) \cong C_\idot(X^{C_m}),
$$
where $X^{C_m}$ is treated as a non-equivariant spectrum, and
$C_\idot(-)$ in the right-hand side is obtained by applying the
limit of Definition~\ref{c.naive.def} to the usual non-equivariant
reduced singular chain complex functor $\overline{C}_\idot(-,\Z)$
($n$ is irrelevant since the chain complexes of
Definition~\ref{cycl.chain} are locally constant). Thus we have to
prove that for any $n,m \geq 1$, the natural map
$$
C_\idot(X^{C_m}) \to (\xi_*C_\idot(X))([n|m])
$$
induced by \eqref{c.bc} is an isomorphism.

On the other hand, by the definition of the functor $C_\idot:G\Sp
\to \DML(\Z)$, we have an isomorphism
$$
\Sigma^V \circ C_\idot \cong C_\idot \circ \Sigma^V
$$
for any $V = V(\nu) \subset U$, and the base change map
$$
C_\idot \circ (\Sigma^V)^{-1} \to (\Sigma^V)^{-1} \circ C_\idot
$$
is also an isomorphism. Therefore \eqref{c.bc} is an isomorphism for
some $X \in G\Sp$ if and only if it is an isomorphism for
$\Sigma^VX$. Since every $X \in G\Sp$ is a filtered colimit of
spectra of the form $(\Sigma^V)^{-1}\Sigma^\infty Y$, $Y \in G\Top$,
it suffices to prove that \eqref{c.bc} is an isomorphism for
suspension spectra $\Sigma^\infty Y$. By Corollary~\ref{le.unive},
we have a quasiisomorphism
$$
C_\idot(\Sigma^\infty Y) \cong \xi^*\overline{C}_\idot(Y);
$$
thus what we have to prove is that for any $Y \in G\Top$ and any
$[n|m] \in \LZ$, the natural map
$$
C_\idot((\Sigma^\infty Y)^{C_m}) \to
(\xi_*\xi^*\overline{C}_\idot(Y))([n|m])
$$
induced by \eqref{c.bc} is a quasiisomorphism. The right-hand side
can be computed by \eqref{xi-xi} and Lemma~\ref{i-pi.lemma}; it is
given by
$$
\bigoplus_{l,p \geq 1, lp=m}
C_\idot(\Z/l\Z,\overline{C}_\idot(Y)([nl|p]),
$$
where by definition, we have
$$
\overline{C}_\idot(Y)([nl|p]) \cong \overline{C}_\idot(Y^{C_p},\Z),
$$
the reduced chain complex of the fixed-points set $Y^{C_p}$. The
desired isomorphism then becomes the tom Dieck-Segal splitting
\cite[Section 1]{HM}.
\endproof

We now return to topology, and recall the following fundamental
notion (we again follow \cite{HM}, and refer to that paper for
further information and references).

\begin{defn}
A {\em cyclotomic structure} on a genuine $G$-spectrum $T \in G\Sp$
is given by a collection of homotopy equivalences
$$
r_m:\wPhi^mT \cong T,
$$
one for each integer $m \geq 1$, such that $r_1= \id$ and $r_n \circ
r_m = r_{nm}$ for any two integer $n,m > 1$.
\end{defn}

\begin{exa}\label{loop.exa}
Assume given a pointed CW complex $X$, and let $LX = \Maps(S^1,X)$
be its free loop space. Then for any finite subgroup $C \subset
S^1$, the isomorphism $S^1 \cong S^1/C$ induces a homeomorphism
$$
\Maps(S^1,X)^C = \Maps(S^1/C,X) \cong \Maps(S^1,X),
$$
and these homeomorphism provide a canonical cyclotomic structure on
the suspension spectrum $\Sigma^\infty LX$.
\end{exa}

Then by Lemma~\ref{phi.chain}, a cyclotomic structure $\{r_m\}$ on
$T$ induces a collection of quasiisomorphisms
$$
r_m:\wt{\rho}_m^*C_\idot(T) \to C_\idot(T), \quad m \geq 1
$$
such that $r_m \circ r_n = r_{nm}$. If we treat the equivariant
chain complex $C_\idot(T)$ as an object $\langle h^*C_\idot(T),\phi
\rangle$ in the category $\D(\LZ_h,P_\idot,\Z)$, then the maps $r_m$
extend $h^*C_\idot(T)$ to a complex in the category
$\Fun(\wLR^{opp},\Z)$, and this extension is compatible with the map
$\phi$. Therefore $C_\idot(T)$ canonically defines an object
$\wt{C}_\idot(T)$ in the derived category
$\D(\wLR_h,P_\idot,\Z)$. By Propositon~\ref{p.all}, this turns
$C_\idot(T)$ into a cyclotomic complex. Thus we can finally justify
our terminology by introducting the following definition.

\begin{defn}
The cyclotomic complex $\wt{C}_\idot(T) \in \DLR(\Z)$ is called the
{\em equivariant chain complex} of the cyclotomic spectrum $T$.
\end{defn}

\section{Filtered Dieudonn\'e modules.}

\subsection{Definitions.}\label{fdm.subs}

We now want to compare cyclotomic complexes to a different and much
simpler algebraic notion which appeared earlier in a different
context -- the notion of a filtered Dieudonn\'e module.

\begin{defn}\label{fdm.ori}
Let $k$ be a finite field of characteristic $p$, with it Frobenius
map, and let $W$ be its ring of Witt vectors, with its canonical
lifting of the Frobenius map. A {\em filtered Dieudonn\'e module}
over $W$ is a finitely generated $W$-module $M$ equipped with a
decreasing filtration $F^\hdot M$, indexed by all integers, and a
collection of Frobenius-semilinear maps $\phi_i:F^iM \to M$, one for
each integer $i$, such that
\begin{enumerate}
\item $\phi_i|_{F^{i+1}M} = p \phi^{i+1}$, and
\item the map
$$
\sum \phi_i:\bigoplus_iF^iM \to M
$$
is surjective.
\end{enumerate}
\end{defn}

This definition was introduced by Fontaine and Lafaille \cite{FL} as
a $p$-adic analog of the notion of a Hodge structure. Under certain
assumptions, the de Rham cohomology $H^\hdot_{DR}(X)$ of a smooth
compact algebraic variety $X/W$ has a natural filtered Dieudonn\'e
module structure, with $F^\hdot$ being the Hodge filtration and the
maps $\phi^\hdot$ induced by the Frobenius endomorphism of the
special fiber $X_k = X \otimes_W k$.

The category of filtered Dieudonn\'e modules is obviously additive,
but there is more: just as for mixed Hodge structures, a small
miracle happens, and the category is actually abelian. For this, the
normalization condition \thetag{ii} of Definition~\ref{fdm.ori}
plays the crucial role. If one is prepared to work with non-abelian
additive categories, this condition can be dropped. For the purposes
of present paper, the following notion will be convenient.

\begin{defn}\label{fdm.defn}
A {\em generalized filtered Dieudonn\'e module} (gFDM for short) is
an abelian group $M$ equipped with
\begin{enumerate}
\item a decreasing filtration $F^\hdot M$, indexed by all integers
  and such that $M = \bigcup F^i M$, and
\item for each integer $i$, each positive integer $j \geq 1$, and
  each prime $p$, a map $\phi^p_{i,j}:F^iM \to M/p^jM$,
\end{enumerate}
such that $\phi^p_{i,j+1} = \phi^p_{i,j} \mod p^j$, and
$\phi^p_{i,j} = p\phi^p_{i+1,j}$ on $F^{i+1}M \subset F^iM$.
\end{defn}

This differs from Definition~\ref{fdm.ori} in that we no longer
require the normalization condition \thetag{ii}, we restrict our
attention to prime fields rather than finite fields, and we collect
together the structures for all primes. Note, however, that if $M$
is finitely generated over $\Z_p$, then all other primes act on $M$
by invertible maps, and the extra maps $\phi^l_{\idot,\idot}$ for $l
\neq p$ are all $0$. In effect, for every prime $p$ and every
integer $i$, we can collect all the maps $\phi^p_{i,\idot}$ into a
single map
$$
\phi^p_i:F^iM \to \wh{(M)}_p,
$$
where $\wh{(M)}_p$ means the pro-$p$ completion of the abelian group
$M$. If $M$ is finitely generated over $\Z_p$, we have $\wh{(M)}_p
\cong M$ and $\wh{(M)}_l=0$ for $l \neq p$.

Complexes of gFDMs are defined in the obvious way. A map between
such complexes is a quasiisomorphism if it induces a
quasiisomorphism of the associated graded quotients $\gr^F$.
Inverting such quasiisomorphisms, we obtain a triangulated ``derived
category of gFDMs'' which we will denote by $\FDM$.

For every gFDM $M$ and any integer $i$, we will denote by $M(i)$ the
same $M$ with the filtration $F^\hdot$ twisted by $i$ -- that is, we
set $F^jM(i)=F^{j-i}M$. Under this convention, we introduce the
folowing ``twisted $2$-periodic'' version of the category $\FDM$.

\begin{defn}
The triangulated category $\FDM^{per}$ is obtained by inverting
quasiisomorphisms in the category of complexes of gFDMs $M_\idot$
equipped with an isomorphism $M \cong M[2](1)$.
\end{defn}

We can now formulate the main result of this Section.

\begin{theorem}\label{fdm.cyclo.prop}
There is a natural equivalence
$$
\FDM^{per} \cong \DLR(\Z)
$$
between the twisted $2$-periodic derived category of gFDMs, on one
hand, and the derived category of cyclotomic complexes of abelian
groups in the sense of Definition~\ref{cyclo.defn}, on the other
hand.
\end{theorem}

We will prove this in Subsection~\ref{fdm.pf.subs}, in the more
precise form of Proposition~\ref{exp.equi.fdm}, after finishing the
necessary preliminaries. We start with some generalities on filtered
objects.

\subsection{Filtered objects.}\label{filt.subs}

By a {\em filtered object} in an abelian category $\Ab$ we will
understand an object $E \in \Ab$ equipped with a decreasing
filtration $F^\hdot$ numbered by all integeres. Maps and complexes
of filtered objects are defined in the obvious way.

\begin{defn}
A map $f:E_\idot \to E'_\idot$ between two filetred complexes of
objects in $\Ab$ is a {\em filtered quasiisomorphism} if the
induced map 
$$
f:F^iE_\idot/F^{i+1}E_\idot \to F^iE'_\idot/F^{i+1}E'_\idot
$$
is a quasiisomorphism for every integer $i$.
\end{defn}

\begin{remark}
We do not require that a filtered quuasiisomorphism induces a
quasiisomorphism of the underlying complexes $E_\idot$, $E'_\idot$
of objects in $\Ab$, nor of the induvidual pieces $F^iE_\idot$,
$F^iE'_\idot$ of the filtrations.
\end{remark}

The {\em filtered derived category} $\DF(\Ab)$ is obtained by
inverting filtered quasiisomorphisms in the category of filtered
complexes and filtered maps. The {\em periodic filtered derived
category} $\DF^{per}(\Ab)$ is similarly obtained from the category
of complexes of filtered objects $V_\idot$ in $\Ab$ equipped with an
isomorphism
$$
u:V_\idot \cong V_{\idot-2}(1),
$$
where $V(1)$ means a twist of filtration. Explicitly, such a complex
is given by the two filtered objects $V_0$, $V_1$, and the two
filtered maps
$$
d_1:V_1 \to V_0, \qquad d_0:V_0 \to V_1(1) \cong V_{-1}
$$
such that $d_1 \circ d_0 = 0 = d_0 \circ d_1$; the other terms in
the complex are then given by $V_{2\idot} = V_0$, $V_{2\idot+1} =
V_1$ with the same differentials $d_0$, $d_1$.

If $\Ab$ is the category of abelian groups, we obtain the periodic
filtered derived category $\DF^{per}(\Z)$ of filtered abelian
groups. For any integer $n$, let $\Z(n) \in \DF^{per}$ be the object
$V$ given by $V_0=\Z$, $v_1 = 0$, with the filtration $F^nV_0=V_0$
and $F^{n+1}V_0=0$. Then the objects $\Z(n)$, $n \in \Z$ generate
the category $\DF^{per}(\Z)$ in the following sense.

\begin{lemma}\label{tate.gen}
Any triangulated subcategory $\D' \subset \DF^{per}(\Z)$ closed
under arbitrary sums and products and containing $\Z$ with the
trivial filtration $F^0\Z=\Z$, $F^1\Z=0$ is equal to the whole
$\DF^{per}(\Z)$.
\end{lemma}

\proof{} Since $\D'$ is closed under taking cones, it obviously
contains any filtered complex $\langle V_\idot,F^\hdot \rangle$ with
bounded filtration $F^\hdot$ (that is, $F^iV_\idot=0$,
$F^jV_\idot=V_\idot$ for some integers $i$, $j$). Assume given an
abritrary filtered complex $\langle V_\idot,F^\hdot \rangle$. Then
the natural maps
\begin{equation}\label{d.lim}
\begin{CD}
V_\idot @<<< \dlim_{\overset{i}{\to}}F^{-i}V_\idot @>>>
\dlim_{\overset{i}{\to}}\dlim_{\overset{j}{\gets}}
F^{-i}V_\idot/F^jV_\idot
\end{CD}
\end{equation}
induce isomorphisms on $\gr^F$, thus become isomorphisms in
$\DF^{per}(\Z)$. Thus we may replace $V_\idot$ with the double limit
in the right-hand side of \eqref{d.lim}. The direct limit can be
computed by the telescope construction. Moreover, for any $i$, the
inverse system $F^iV_\idot/F^\hdot V_\idot$ satisfies the
Mittag-Leffler condition, so that the inverse limit can also be
computed by the telescope construction. Since $\D'$ is closed under
products and sums, it is also closed under telescopes, so that it
must contain $V_\idot$.
\endproof

We note that by \eqref{d.lim}, we may represent any object in
$\DF^{per}(\Z)$ by a filtered complex $\langle V_\idot,F^\hdot
\rangle$ which is {\em admissible} in the following sense: both
natural maps
$$
\lim_{\overset{i}{\to}}F^{-i}V_\idot \to V_\idot, \qquad
V_\idot \to \lim_{\overset{i}{\gets}}V_\idot/F^iV_\idot
$$
are isomorphisms of complexes.

To work with filtered abelian groups, it is convenient to use Rees
objects. Consider the algebra $\Z[t]$ of polynomials in one variable
$t$. Say that a module $M$ over $\Z[t]$ is {\em $t$-adically
complete} if the natural map
$$
M \to \lim_{\overset{i}{\gets}} M/t^iM
$$
is an isomorphism (thus in our terminology, ``complete'' includes
``separated''). We turn $\Z[t]$ into a graded ring by assigning
degree $-1$ to the generator $t$.

\begin{lemma}\label{rees}
The filtered derived category $\DF(\Z)$ is equivalent to the full
subcategory in the derived category of the abelian category of
$\Z$-graded $\Z[t]$-modules $M_\idot$ spanned by $t$-adically
complete modules.
\end{lemma}

\proof{} For any filtered abelian group $\langle M,F^\hdot \rangle$,
the corresponding graded $\Z[t]$-module $\wt{M}_\idot$ called the
{\em Rees object} of $M_\idot$ is given by
$$
\wt{M}_\idot = \bigoplus \dlim_{\overset{i}{\to}}F^\hdot
M/F^{\hdot+i}M,
$$
with $t$ induced by the natural embeddings $F^\hdot M \to
F^{\hdot-1}M$. We note that $\wt{M}_\idot$ is automatically
$t$-adically complete, and a filtered quasiisomorphism of complexes
of abelian groups induces a quasiisomorphism of Rees objects.

To get the inverse correspondence, note that every graded
$\Z[t]$-module $M_\idot$ has a finite resolution by modules with no
$t$-torsion, so that it is enough to consider graded modules
$M_\idot$ with injective map $t$. Such a module $M_\idot$ is sent to
$$
\wt{M} = \lim_{\overset{t}{\to}} M_\idot,
$$
with $F^i \wt{M} \subset \wt{M}$ being the image of the natural
embedding $M_i \to \wt{M}$ for any integer $i$.
\endproof

The equivalence of Lemma~\ref{rees} has an obvious periodic version:
the Rees object model of the category $\DF^{per}(\Z)$ is obtained by
inverting quasiisomorphisms in the category of pairs $\langle
\wt{M}_{\idot,\idot},u \rangle$ of a complex $\wt{M}_{\idot,\idot}$
of graded $t$-adically complete $\Z[t]$-modules $M_\idot$ and an
isomorphism
\begin{equation}\label{u.rees}
\wt{u}:\wt{M}_{\idot,\idot} \cong \wt{M}_{\idot+1,\idot-2}.
\end{equation}

\subsection{Cyclic expansion and subdivision.}\label{exp.subs}

Now let $\Ab = \Fun(\Lambda^{opp},\Z)$ be the category of cyclic
abelian groups. Let $I_0 = j^o_!\Z$, $I_1 = j_*\Z$ be as in
\eqref{4.term}, and let $d_1 = B:I_1 \to I_0$, $d_0 = b_0 \circ
b_1:I_0 \to I_1$, where $b_0$, $b_1$ and $B$ are again as in
\eqref{4.term}. Then since \eqref{4.term} is exact, we have $d_1
\circ d_0 = 0 = d_0 \circ d_1$. Moreover, if we define filtrations
$F^\hdot$ on $I_0$ and $I_1$ by setting $F^0I_l = I_l$, $F^1I_l=0$,
$l=0,1$, then both $d_0$ and $d_1$ are filtered maps, so that we
have a periodic filtered complex $I_\idot$ of objects in
$\Fun(\Lambda^{opp},\Z)$.

\begin{defn}
For any object in $\DF^{per}(\Z)$ represented by a periodic complex
$V_\idot$ of admissible filtered abelian groups, its {\em cyclic
expansion} $\Exp(V_\idot)$ is a complex of cyclic abelian groups
given by
$$
\Exp(V_\idot) = F^0(V_\idot \otimes_{\Z[u]} I_\idot),
$$
where $u$ is the periodicity map on $V_\idot$ and $I_\idot$, and
$F^0$ is taken with respect to the product filtration.
\end{defn}

In terms of the corresponding peridoc complex $\wt{V}_{\idot,\idot}$
of Rees objects, cyclic expansion is given by
\begin{equation}\label{exp.expl}
\Exp(V_\idot)_i = (\wt{V}_{0,\idot} \otimes I_1)[1] \oplus
(\wt{V}_{0,\idot} \otimes I_0),
\end{equation}
with the differential $d = d_V \otimes \id + d_I$, where $d_V$ is
the differential on $\wt{V}_{0,\idot}$, and $d_I$ is equal to $\id
\otimes d_1$ on the first summand, and to $t\wt{u} \otimes d_0$ on
the second one, where $\wt{u}:\wt{V}_{\idot,\idot} \cong
\wt{V}_{\idot+1,\idot-2}$ is as in \eqref{u.rees}. In particular,
for every object $[n] \in \Lambda$, $\Exp(V_\idot)([n])$ is the sum
of a finite number of copies of $\wt{V}_{0,\idot}$ and its shift
$\wt{V}_{0,\idot}[1]$; therefore cyclic expansion commutes with
arbitrary sums and arbitrary products.

\begin{lemma}\label{exp.equi}
Cyclic expansion induces an equivalence of categories
$$
\Exp:\DF^{per}(\Z) \cong \D_c(\Lambda^{opp},\Z)
$$
between $\DF^{per}(\Z)$ and the full subcategory
$\D_c(\Lambda^{opp},\Z) \subset \D(\Lambda^{opp},\Z)$ spanned by
objects which are locally constant in the sense of
Definition~\ref{norm.defn}.
\end{lemma}

\proof{} Since \eqref{4.term} is exact, $F^1I_\idot = I_{<0}$ is a
resolution of the constant functor $\Z \in \Fun(\Lambda,\Z)$. By
induction, this immediately implies that for any periodic filtered
complex $V_\idot$, the homology functors $\HH_\idot(\Exp(V_\idot))
\in \Fun(\Lambda^{opp},\Z)$ are given by
$$
\HH_i(\Exp(V_\idot)) \cong (F^0V_i/F^1V_i) \otimes \Z.
$$
In particular, a filtered quasiisomorphism of periodic filtered
complexes induces a quasiisomorphism of their cyclic expansions, so
that $\Exp$ induces a well-defined triangulated functor from the
category $\DF^{per}(\Z)$ to the full subcategory
$\D_c(\Lambda^{opp},\Z) \subset \D(\Lambda^{opp},\Z)$. Moreover, for
any integer $n$, we have
$$
\Exp(\Z(n)) \cong \Z[2n].
$$
Since the fundamental group $\pi_1(|\Lambda|) \cong \pi_1(BU(1))$ is
trivial, every locally constant functor $\Lambda^{opp} \to \Z\amod$
must be constant; therefore $\Z$ generates the triangulated category
$\D_c(\Lambda^{opp},\Z)$ in the same sense as in
Lemma~\ref{tate.gen}. Since $\Exp$ commutes with arbitrary sums and
arbitrary products, it therefore suffices to prove that $\Exp$ is
fully faithful on the objects $\Z(n)$ -- that is, the natural map
$$
\Exp:\RHom^\hdot_{\DF^{per}(\Z)}(\Z(n),\Z(m)) \to
\RHom^\hdot_{\D(\Lambda^{opp},\Z)}(\Z[2n],\Z[2m])
$$
is a quasiisomorphism for any two integers $n$, $m$. This
immediately follows from the isomorphism $H^\hdot(\Lambda^{opp},\Z)
\cong \Z[u]$.
\endproof

Now assume given a positive integer $n \geq 1$, and let
$\delta^n:\Z[t] \to \Z[t]$ be the map which sends $t$ to $nt$. Since
the the ideal $(nt) \subset \Z[t]$ lies inside $(t) \subset \Z[t]$,
the direct image $\delta^n_*M$ of a $t$-adically complete
$\Z[t]$-module $M$ is automatically $t$-adically complete.

\begin{defn}\label{divn.def}
For any positive integer $n \geq 1$ and a filtered abelian group $M$
with the corresponding admissible Rees object $M_\idot$ as in
Lemma~\ref{rees}, the {\em $n$-th subdivision} $\Div_n(M)$ is the
filtered abelian group corresponding to the graded abelian group
$\delta^n_*M_\idot$.
\end{defn}

In other words, $M_\idot$ remains the same as a graded abelian
group, but the map $t$ is replaced by its multiple $nt$. We note
that the underlying filtered abelian group $M$ itself might change
under subdivision: for example, if $M = \Z(0) = \Z$, then
$$
\Div_n(M) \cong \Qq,
$$
with the filtration given by $F^1\Qq=0$, $F^i\Qq = n^i\Z \subset \Qq$
for $i \leq 0$.

For any $n$, $\Div_n$ is obviously an endofunctor of the category of
filtered abelian groups, and it descends to an endofunctor
$$
\Div_n:\DF^{per}(\Z) \to \DF^{per}(\Z)
$$
of the periodic derived category $\DF^{per}(\Z)$.

Fix an integer $n \geq 1$, and recall the two functors
$i_n,\pi_n:\Lambda_n \to \Lambda$ of Subsection~\ref{connes.subs}.

\begin{lemma}\label{exp.div}
For any periodic complex $V_\idot$ of admissible filtered abelian
groups, we have a functorial isomorphism
$$
\pi_{n*}i_n^*\Exp(V_\idot) \cong \Exp(\Div_n(V_\idot)).
$$
\end{lemma}

\proof{} Let $I'_l = \pi_{n*}i_n^*I_l$, $d'_l = \pi_{n*}i_n^*d_l$,
$l=0,1$. Applying the functor $\pi_{n*}i_n^*$ to \eqref{exp.expl},
we see that the complex $\pi_{n*}i_n^*\Exp(V_\idot)$ is given by
$$
\pi_{n*}i_n^*\Exp(V_\idot)_i = (\wt{V}_{0,\idot} \otimes I_1)[1]
\oplus (\wt{V}_{0,\idot} \otimes I_0),
$$
with the differential $d = d_V \otimes \id + d_I$, where $d_V$ is
the differential on $\wt{V}_{0,\idot}$, and $d_I$ is equal to $\id
\otimes d'_1$ on the first summand, and to $t\wt{u} \otimes d'_0$ on
the second one. By \eqref{n.ti}, we have canonical isomorphisms
$I'_0 \cong I_0$, $I'_1 \cong I_1$, and under these isomorphisms, we
have $d'_1 = d_1$ and $d'_0=nd_0$. Thus $t\wt{u} \otimes d_0' =
nt\wt{u} \otimes d_0$, and $\pi_{n*}i_n^*\Exp(V_\idot)$ is exactly
isomorphic to the expression \eqref{exp.expl} for the complex
$\Exp(\Div_n(V_\idot))$.
\endproof

\subsection{Stabilization.}\label{stab.subs}

For any filtered abelian group $M$ with the corresponding admissible
Rees object $M_\idot$, we have a tautological map
$M \to M(1)$; applying the subdivision functor $\Div_n$, we obtain a
natural map
\begin{equation}\label{tau.div}
\Div_n(M) \to \Div_n(M(1)).
\end{equation}

\begin{defn}\label{stab.defn}
In the situation of Definition~\ref{divn.def}, the {\em stabilized
  $n$-th subdivision} $\Stab_n(M)$ is given by
$$
\Stab_n(M) = \lim_{\overset{l}{\to}} \Div_n(M(l)),
$$
where the limit is taken with respect to the tautological maps
\eqref{tau.div}.
\end{defn}

\begin{lemma}\label{stab.comple}
For any filtered abelian group $M$ which is admissible in the sense
of Subsection~\ref{filt.subs} and any prime $p \geq 2$, we have
$$
\Stab_p(M) \cong \wh{M}_p \otimes_{\Z_p} \Q_p,
$$
with the filtration given by
$$
F^l\wh{M}_p = p^l\wh{M}_p,
$$
where $\wh{M}_p = \lim_{\overset{l}{\gets}} M/p^lM$ denotes the
pro-$p$ completion of the group $M$.
\end{lemma}

\proof{} On the level of Rees objects,
$$
M_\idot = \lim_{\to} \Div_p(M)_\idot
$$
is isomorphic to $M$ in every degree $l$, $M_l \cong M$, with $t \in
\Z[t]$ acting by multiplication by $p$. However, this is not
$t$-adically complete. Thus when we apply the equivalence of
Lemma~\ref{rees}, the result is
$$
\lim_{\overline{l}{\gets}} M/t^lM_\idot \cong \wh{M}_p,
$$
as required.
\endproof

As a corollary, we see that the periodic derived category
$\FDM^{per}$ of generalized filtered Dieudonn\'e modules of
Subsection~\ref{fdm.subs} is equivalent to the category of filtered
periodic complexes $V_\idot$ of abelian groups equipped with a map
$$
\phi_p:V_\idot \to \Stab_p(V_\idot)
$$
for any prime $p \geq 2$.

Now let $I_\idot$ be the periodic complex of functors in
$\Fun(\Lambda,\Z)$ of Subsection~\ref{exp.subs}, and let $P_\idot
\subset I_\idot$ be its canonical truncation at $0$. Explicitly, we
have $P_0 = \Z$, the constant functor, and
$$
P_{2l-1} = I^1, \qquad P_{2l} = I^0
$$
for any $l \geq 1$. The complex $P_\idot$ is acyclic. Moreover, the
pullback $i^*P_\idot$ with respect to the functor $i:\LI_{red} \to
\Lambda$ of Subsection~\ref{como.subs} is admissible in the sense of
Lemma~\ref{adm.equi}.

Let $\eta^+:\Z \to F^0 I_\idot$, $\eta^-:\Z \cong P_0 \to P_\idot$
be the natural embeddings, and define a map
$$
\eta:I_\idot \to P_\idot \otimes (F^0I_\idot)
$$
by
\begin{equation}\label{eta.def}
\eta =
\begin{cases}
\eta^- \otimes \id, &\text{ on }I_l, l \leq 0,\\
\id \otimes \eta^+, &\text{ on }I_l, l \geq 0.
\end{cases}
\end{equation}
Let $F^lP_\idot$, $l \geq 0$ be the stupid filtration on the complex
$P_\idot$. Then for any $l \geq 0$, the map $\eta$ induces a map
\begin{equation}\label{eta.qui}
\eta:F^{-l}I_\idot \to F^0I_\idot \otimes F^{2l}P_\idot,
\end{equation}
and this map is a quasiisomorphism (both sides are quasiisomorphic
to $\Z[2l]$).

\begin{lemma}
For any $n \geq 1$, $l \geq 0$, and any complex $V \in
\Fun(\Lambda_n^{opp},Z)$, the natural map
\begin{equation}\label{eta.sta}
\eta:\pi_{n*}(V \otimes i_n^*F^{-l}I_\idot) \to \pi_{n*}(V \otimes
i_n^*F^0I_\idot \otimes i_n^*F^{2l}P_\idot)
\end{equation}
is a quasiisomorphism.
\end{lemma}

\proof{} Both sides are finite-length complexes in
$\Fun(\Lambda_n^{opp},\Z)$, and after evaluating at an object $[m]
\in \Lambda_n$, both sides give complexes of free
$\Z[\Z/n\Z]$-modules. Therefore we may replace $\pi_{n*}$ with its
derived functor $R^\hdot\pi_{n*}$. Then the claim immediately
follows from the fact that \eqref{eta.qui} is a quasiisomorphism.
\endproof

\begin{corr}\label{eta.corr}
For any integer $n \geq 1$ and any twisted periodic complex
$V_\idot$ of filtered abelian groups, the map $\eta$ of
\eqref{eta.def} induces a quasiisomorphism
$$
\begin{CD}
\Exp(\Stab_n(V_\idot)) @>{\sim}>>
\dlim_{\overset{l}{\to}}\pi_{n*}i_n^*(\Exp(V_\idot) \otimes F^lP_\idot).
\end{CD}
$$
\end{corr}

\proof{} By Lemma~\ref{exp.div} and Definition~\ref{stab.defn}, we
have
$$
\Exp(\Stab_n(V_\idot)) \cong
\dlim_{\overset{l}{\to}}\pi_{n*}i_n^*\Exp(V_\idot(l)),
$$
and by definition, we have
$$
\Exp(V_\idot(l)) = F^0(V_\idot \otimes F^{-l}I_\idot)
$$
for any integer $l \geq 0$. 
\endproof

\subsection{Comparison.}\label{fdm.pf.subs}

Consider now the category $\D(\LR_h,i^*P_\idot,k)$ of
Subsection~\ref{como.subs}, with the specific choice of the complex
$P_\idot$ made in Subsection~\ref{stab.subs}. Let $\tau:\LR_h \cong
\Lambda \times I \to \Lambda$ be the tautological projection. Then
for any twisted periodic complex $V_\idot$ of abelian groups, a
collection of maps $\phi_p:V_\idot \to \Stab_p(V_\idot)$ for all
prime $p$ induces a map
$$
\eta \circ \tau^*\Exp(\phi):\tau^*\Exp(V_\idot) \to
\pi_*i^*\tau^*(\Exp(V_\idot) \otimes P_\idot),
$$
which gives by adjunction a map
$$
\wt{\phi}:\pi^*\tau^*\Exp(V_\idot) \to i^*(\tau^*\Exp(V_\idot)
\otimes i^*P_\idot).
$$
Sending $\langle V_\idot,\{\phi_p\}\rangle$ to $\langle
\tau^*\Exp(V_\idot),\wt{\phi} \rangle$ then defines a comparison
functor
\begin{equation}\label{comp}
\Exp:\FDM^{per} \to \D_c(\LR_h,P_\idot,\Z).
\end{equation}
By Proposition~\ref{p.all}, the following result immediately yields
Theorem~\ref{fdm.cyclo.prop}.

\begin{prop}\label{exp.equi.fdm}
The functor $\Exp$ of \eqref{comp} is an equivalence of categories.
\end{prop}

To prove this result, we need to generalize it. Consider the twisted
periodic filtered derived category $\DF^{per}(I^{opp},\Z)$ of
functors from $I^{opp}$ to abelian groups. The cyclic expansion
functor gives a functor
\begin{equation}\label{exp.w}
\Exp:\DF^{per}(I^{opp},\Z) \to \D_w(\LR_h^{opp},\Z),
\end{equation}
where $\D_w(\LR_h^{opp},\Z) \subset \D(\LR_h^{opp},\Z)$ is as in
Proposition~\ref{reso.bis}. Lemma~\ref{exp.equi} immediately shows
that this functor $\Exp$ is an equivalence of categories. Moreover,
for any $n \geq 1$, let $\nu_n:I \to I$ be the functor given by
multiplication by $n$, as in \eqref{i.p.red}, and for any twisted
periodic filtered complex $M \in \DF^{per}(I^{opp},\Z)$, let
$$
\wt{\Stab_n}(M) = \nu_n^*\Stab_n(M).
$$
Let $\wt{\FDM}^{per}$ be the derived category of such complexes $M$
equipped with a collection of maps
$$
\phi_p:M \to \wt{\Stab}_p(M)
$$
for all primes $p$. Denote by $\wt{h}^*:\wt{\FDM}^{per} \to
\DF^{per}(I^{opp},\Z)$ the forgetful functor. It has an obvious
right-adjoint $\wt{h}_*$ given by
\begin{equation}\label{h.stab}
\wt{h}_*(M)(n) = M \oplus \prod_{p | n}\Stab_p(M),
\end{equation}
where the product is taken over all prime divisors of the integer $n
\in I$. Moreover, the construction of the functor $\Exp$ of
\eqref{comp} also gives a functor
\begin{equation}\label{w.comp}
\wt{\Exp}:\wt{\FDM}^{per} \to \D_w(\LR_h,i^*P_\idot,\Z),
\end{equation}
and by construction, we have
\begin{equation}\label{exp.bc}
\wt{h}^*\circ\wt{\Exp} \cong \Exp \circ \wt{h}^*,
\end{equation}
where $\wt{h}^*$ in the left-hand side is the restriction functor of
Proposition~\ref{reso.bis}.

\begin{prop}\label{w.equi}
The functor \eqref{w.comp} is an equivalence of categories.
\end{prop}

\proof{} As in the proof of Proposition~\ref{mack.equi},
Proposition~\ref{reso.bis}~\thetag{ii} implies that the functor
$\wt{\Exp}$ has a right-adjoint functor
$$
K:\D_w(\LR_h,i^*P_\idot,\Z) \to \wt{\FDM}^{per}.
$$
Moreover, by \eqref{h.stab}, Corollary~\ref{eta.corr} and
Lemma~\ref{reso.bis.lemma}, the base change map
$$
\wt{\Exp} \circ \wt{h}_* \to \wt{h}_* \circ \Exp
$$
induced by \eqref{exp.bc} is an isomorphism. Therefore by
adjunction,
$$
\wt{h}^* \circ K \cong \Exp^{-1} \circ \wh{t}^*,
$$
where $\Exp^{-1}$ is the equivalence inverse to =eqref{exp.w}. We
conclude that
$$
\wt{h}^* \circ K \circ \wt{\Exp} \cong \Exp^{-1} \circ \Exp \circ
\wt{h}^*, \qquad \wt{h}^* \circ \wt{\Exp} \circ K \cong \Exp^{-1}
\circ \Exp \circ \wt{h}^*.
$$
As in the proof of Proposition~\ref{mack.equi}, the functor
$\wt{h}^*$ is conservative, and $\Exp$ and $\Exp^{-1}$ are
mutually inverse equivalences of categories. Therefore so are the
functors $\wt{\Exp}$ and $K$.
\endproof

\proof[Proof of Proposition~\ref{exp.equi.fdm}.] The tautological
projection $\tau:I \to \ppt$ induces a functor
$$
\tau^*:\DF^{per}(\Z) \to \DF^{per}(I^{opp},\Z).
$$
This is a full embedding onto the full subcategory spanned by the
locally constant functors; the adjoint functor is given by
$$
\tau_* = R^\hdot\dlim_{\overset{I^{opp}}{\gets}}.
$$
By Lemma~\ref{stab.comple}, stabilized subdivision commutes with
arbitrary products, so that $\tau_* \circ \wt{\Stab}_p \cong \Stab_p
\circ \tau_*$. Therefore $\tau^*$ and $\tau_*$ extend to an adjoint
pair of functors between $\FDM^{per}$ and $\wt{\FDM}^{per}$, so that
$\FDM^{per}$ is identified with the full subcategory in
$\wt{\FDM}^{per}$ spanned by objects $M$ with locally constant
$\wt{h}^*M$. By \eqref{exp.bc}, the equivalence of
Proposition~\ref{w.equi} identifies this subcategory with $\DLR(\Z)
\cong \D_c(\LR_h,i^*P_\idot,\Z) \subset \D_w(\LR_h,i^*P_\idot,\Z)$.
\endproof

\section{Topological cyclic homology.}

We finish the paper with a brief discussion of topological cyclic
homology (we again follow \cite{HM}).

Fix the $G$-universe $U$ as in Subsection~\ref{cyclo.subs}, and
recall that for any $m \geq 1$, we have canonical functors
$\wPsi^m,\wPhi^m$ endofunctors of the category $G\Sp = G\Sp(U)$. We
also have natural maps $\can:\wPsi^m \to \wPhi^m$. For any $T \in
G\Sp$ and a pair of integers $r,s > 1$, one has a natural
non-equivariant map
$$
F_{r,s}:T^{C_{rs}} \to T^{C_r}.
$$
On the other hand, assume that $T$ is equipped with a cyclotomic
structure. Then we have a natural map
$$
\begin{CD}
R_{r,s}:T^{C_{rs}} \cong (\wPsi^s(T))^{C_r} @>{\can}>>
(\wPhi^s(T))^{C_r} @>{r_s}>> T^{C_r},
\end{CD}
$$
where $r_s$ comes from the cyclotomic structure on $T$. Taken
together, the maps $F_{r,s}$ and $R_{r,s}$ define a functor $I(T)$
from the category $\I^{opp}$ of Subsection~\ref{lz.subs} to the
category of non-equivariant spectra: we let $I(T)(n) = T^{C_n}$ for
any $n \in \I$, and we let the morphisms $F_r,R_r:s \to rs$ act by
the maps $F_{r,s}$ resp. $R_{r,s}$.

\begin{defn}\label{Tc.defn}
The {\em topological cyclic homology} $\TC(T)$ of a cyclotomic
spectrum $T$ is a non-equivariant spectrum given by
$$
\TC(T) = \holim_{\I^{opp}} I(T).
$$
\end{defn}

Assume now given a cyclotomic complex $M \in \DLR(\Z)$, and consider
the functor $\wt{\alpha}:\I \to \wLR$ of \eqref{alpha}.

\begin{defn}\label{tc.def}
The {\em topological cyclic homology} $\TC_\idot(M)$ of a cyclotomic
complex $M \in \DLR(\Z)$ is given by
$$
\TC_\idot = H^\hdot(\I^{opp},\wt{\alpha}^*\xi_*M),
$$
where $\xi_*:\DLR(\Z) \cong \D_c(\wLR,\wt{\lambda}^*\T_\idot,\Z) \to
\D_c(\wLR^{opp},\Z)$ is the right-adjoint to the corestriction
functor.
\end{defn}

In particular, assume given a $G$-spectrum $T$, and consider its
equivariant chain complex $C_\idot(T) \in \DML(\Z)$ of
Subsection~\ref{chain.subs}. 

\begin{prop}
For any cyclotomic spectrum $T$, we have a natural isomorphism
$$
\TC_\idot(C_\idot(T)) \cong H_\idot(\TC(T),\Z),
$$
where $H_\idot(-,\Z)$ denotes the homology of a non-equivariant
spectrum with coefficients in $\Z$.
\end{prop}

\proof{} By Lemma~\ref{psi.bc}, we have
$$
\wt{\alpha}^*\xi_*C_\idot(T)(m) \cong C_\idot(T^{C_m})
$$
for any $m \in \I$, so that $\wt{\alpha}^*\xi_*C_\idot(T) \in
\D(\I^{opp},\Z)$ is isomorphic to $C_\idot(I(T),\Z)$, the
non-equivariant chain homology complex of the system of spectra
$I(T)$ used in Definition~\ref{Tc.defn}.
\endproof

In view of the equivalence of Theorem~\ref{fdm.cyclo.prop}, it would
be desirable to express the topological cyclic homology functor
$\TC_\idot$ of Definition~\ref{tc.def} in terms of filtered
Dieudonn\'e modules. A natural notion of homology for filtered
Dieudonn\'e modules is the following.

\begin{defn}
{\em Syntomic cohomology} of a generalized filtered Dieudonn\'e
module $M \in \FDM^{per}(\Z)$ is given by
$$
\RHom^\hdot(\Z,M),
$$
where $\Z \in \FDM^{per}$ is the trivial filtered Dieudonn\'e
module.
\end{defn}

In general, syntomic cohomology and topological cyclic homology are
different. However, to have the following
comparison result.

\begin{theorem}
Assume that $M \in \DLR(\Z) \cong \FDM^{per}(\Z)$ is profinitely
complete. Then we have a natural isomorphism
$$
\TC_\idot(M) \cong \RHom^\hdot(\Z,M).
$$
\end{theorem}

\proof{} In term of cyclotomic complexes, the trivial Dieudonn\'e
module $\Z$ corresponds to the corestriction $\wt{\xi}^*\Z$ of the
constant functor $\Z \in \Fun(\wLR^{opp},\Z)$. Then by adjunction,
we have
$$
\RHom^\hdot(\Z,M) \cong H^\hdot(\wLR^{opp},\xi_*M).
$$
Thus it suffices to prove that for any profinitely complete $M \in
\D(\wLR^{opp},\Z)$, the natural map
$$
H^\hdot(\I^{opp},\wt{\alpha}^*M) \to H^\hdot(\wLR^{opp},M)
$$
is an isomorphism. Equivalently, we have to prove that the map
$$
H^\hdot([1/\N^*],\wt{\lambda}_*\wt{\alpha}^*M) \to
H^\hdot(\LR^{opp},\wt{\lambda}_*M)
$$
is an isomorphism. By base change, $\wt{\lambda}_*\wt{\alpha}^*
\cong \alpha^*\wt{\lambda}_*$, and since right-derived Kan
extensions commute with profinite completions, we apply
Proposition~\ref{profini} and replace $\LR$ with $\DeR$. But we know
that $\wt{j}^*\wt{\lambda}_*M$ becomes locally constant after
restricting to $\Delta \subset \DeR$, and we claim that for any $M' \in
\D(\DeR,\Z)$ with locally constant $h^*M' \in \D(\Delta,k)$, the map
$$
H^\hdot([1/\N^*],\alpha^*M') \to H^\hdot(\DeR,M')
$$
is an isomorphism. Indeed, by Lemma~\ref{der.le}, we can compute the
direct image $\delta_!$ fiberwise, and since $h^*M'$ is locally
constant, the adjunction map
$$
\delta^*\delta_!M' \to M'
$$
is an isomorphism. Then
$$
H^\hdot(\DeR,M') \cong H^\hdot(\DeR,\delta^*\delta_!M') \cong
H^\hdot([1/\N^*],\delta_!M'),
$$
and since $\delta \circ \alpha = \id$, the right-hand side is
exactly $H^\hdot([1/\N^*],\alpha^*M')$.
\endproof

In conclusion, let me say that to obtain an analogous comparison
isomorphism for an arbitrary cyclotomic spectrum $T$, one has to
modify the definition of topological cyclic homology $\TC(T)$ by
replacing the fixed-points functors $T^{C_m}$ in the system $I(T)$
with their homotopy fixed points
$$
\left(T^{C_m}\right)^{hG}
$$
with respect to the residual action of the group $G = U(1)$ (this
should be related to the usual $\TC$ by a cofiber sequence analogous
to \eqref{connes.lr}). I do not know whether it makes sense to do it
from the topological point of view. I am very grateful to
L. Hesselholt for explaining to me why it does not matter in the
profinitely complete case.

\section{Appendix.}
\def\thesection{A}

In this Appendix, we collect some facts used throughout the paper;
most of the facts are well-known, but our terminology may be
non-standard.

\medskip

\noindent
First of all, a piece of general nonsense. Assume given a square
$$
\begin{CD}
A @>{b}>> B\\
@V{a}VV @VV{c}V\\
C @>{d}>> D,
\end{CD}
$$
and assume that $a$ and $c$ admit left-adjoint functors $a_!$,
$c_!$. Then an isomorphism $d \circ a \cong c \circ b$ induces by
adjunction a map
$$
c_! \circ b \to a_1 \circ b.
$$
We call this map and its various adjoints the {\em base change map}
induced by the isomorphism $d \circ a \cong c \circ b$.

For any category $\C$ with objects $c,c' \in \C$, we denote by
$\C(c,c')$ the set of maps from $c$ to $c'$, and we denote by
$\C^{opp}$ the opposite category, $\C^{opp}(c,c') = \C(c',c)$. For a
small category $\C$ and a category $A$, we denote by $\Fun(\C,A)$
the category of cuntors from $\C$ to $A$. For a functor $f:\C \to
\C'$, $f^*:\Fun(\C',A) \to \Fun(\C,A)$ is the pullback, and $f_!$,
$f_*$ are the left and the right Kan extensions (when they
exists). If $A$ is abelian, we denote by $\D(\C,A)$ the derived
category of the category $\Fun(\C,A)$, and by abuse of notation, we
use $f_!$ and $f_*$ to denote the derived functors of the Kan
extensions. If $A = k\amod$, the category of modules over a ring
$k$, we shorten $\Fun(\C,k\amod)$, $\D(\C,k\amod)$ to $\Fun(\C,k)$,
$\D(\C,k)$. The homology $H_\idot(\C,k)$ resp. cohomology
$H^\hdot(\C,k)$ of a small category $\C$ is defined by taking
dervied Kan extensions with respect to the projection $\C \to
\ppt$. As usual, $H^\hdot(\C,k)$ is an algebra, and $H_\idot(\C,k)$
is a module over $H^\hdot(\C,k)$. The homology can be computed by an
explicit {\em bar complex} $C_\idot(\C,k)$.

We also introduce the following slightly non-standard definition.

\begin{defn}\label{norm.defn}
Assume given a small category $\C$.  An object $E_\idot \in
\D(\C,k)$ is {\em locally constant} if for any map $f:a \to b$ in
$\C$, the map $E_\idot(f):E_\idot(a) \to E_\idot(a)$ is a
quasiisomorphism.
\end{defn}

We freely use the notions from \cite{SGA} (Cartesian map, fibration,
cofibration, bifibration); for a brief overview with exactly the
same notation as in this paper, see \cite[Section 1]{K}. We also
freely use the base change isomorphism and projection formula of
\cite[Lemma 1.7]{K}. We also use the following notion from
\cite{bou}.

\begin{defn}\label{facto.def}
A {\em factorization system} on a category $\C$ is given by two
subcategories $\C_v,\C_h \subset \C$ such that all isomorphisms in
$\C$ lie both in $\C_v$ and in $\C_h$, and any morphism $f$ in $\C$
decomposes as $f = v \circ h$, $v \in \C_v$, $h \in \C_h$, and such
a decomposition is unique up to a unique isomorphism.
\end{defn}

\begin{exa}
If $\gamma:\C \to \C'$ is a fibration, then fiberwise and Cartesian
maps form a factorization system on $\C$.
\end{exa}

Factorization systems are actually very useful gadgets, although
they are traditionally relegated to appendices and introductions
(and we follow the tradition). Definition~\ref{facto.def} has
several corollaries and/or equivalent reformulations. For example,
one can show that $\C_v \cap \C_h$ exactly consists of all the
isomorphisms in $\C$; moreover, maps in $\C_v$ have a unique lifting
property with respect to maps in $\C_h$, and vice versa. We refer
the reader to \cite{bou} for discussion. We will need one result
which is not in \cite{bou}. Assume that the category $\C$ is small,
and let $\overline{\C}$ be the category of all objects in $\C$ and
all isomorphisms between them, so that we have a Cartesian square
$$
\begin{CD}
\overline{\C} @>{\overline{v}}>> \C_h\\
@V{\overline{h}}VV @VV{h}V\\
\C_v @>{v}>> \C,
\end{CD}
$$
where $v$, $h$, $\overline{v}$, $\overline{h}$ are the embedding
functors. Then the isomorphism $\overline{v}^* \circ h^* \cong
\overline{h}^* \cong v^*$ induces a base change map
\begin{equation}\label{fact.bc.eq}
\overline{h}_! \circ \overline{v}^* \to v^* \circ h_!
\end{equation}
of functors from $\D(\C_h,k)$ to $\D(\C_v,k)$.

\begin{lemma}\label{fact.bc}
The base change map \eqref{fact.bc.eq} is an isomorphism.
\end{lemma}

\proof{} Since the category $\D(\C_h,k)$ is generated by
representable functors of the form
$$
k_c(c') = k[\C_h(c,c')],
$$
$c,c' \in \C_h$, it suffices to prove that \eqref{fact.bc.eq}
becomes an isomorphism after applying to some such functor $k_c \in
\Fun(\C_h,k)$, $c \in \C_h$. Indeed, $h_!$ sends representable
functors into representable ones, so that we have
\begin{equation}\label{v-h}
v^*h_!k_c(c') = h_!k_c(c') = k[\C(c,c')]
\end{equation}
for any $c' \in \C$. But by the definition of a factorization
system, we have a natural isomorphism
\begin{equation}\label{fac.eq}
\C(c,c') = \coprod_{c'' \in \C} (\C_h(c,c'') \times
\C_v(c'',c'))/\Aut(c''),
\end{equation}
and the actions of the groups $\Aut(c'')$ are free. Therefore
\eqref{v-h} is isomorphic to
$$
\overline{h}_! \bigoplus_{c'' \in \overline{\C}}k[\C_h(c,c'')],
$$
and the right-hand side is exactly $\overline{v}^*k_c$.
\endproof

We also use the $A_\infty$-technology; a brief introduction to the
relevant part of it is contained in \cite[Section 1.5]{Ka-ma}. Here
we just recall that an $A_\infty$-algebra is an algebra over a
certain asymmetric operad $\Ass_\infty$ of complexes of abelian
groups, a cofibrant resolution of the associative asymmetric operad
$\Ass$. An $A_\infty$-coalgebra is an $A_\infty$-algebra in the
opposite category. Since the operads are asymmetric, one can define
algebras in an arbitrary tensor category. In particular, for any set
$S$, one can consider the category of abelian groups graded by $S
\times S$, with tensor product given by
$$
(V \otimes V')_{s,s'} = \bigoplus_{s'' \in S} V_{s,s''} \otimes
V'_{s'',s}.
$$
An $A_\infty$-algebra in this category is a small $A_\infty$-category
with the set of objects $S$. More generally, given a small category
$\C$ with the of objects $\C_0$ and the set of morphisms $\C_1$, one
can consider the category of $\C_1$-graded vector spaces, with the
tensor product given by
$$
(V \otimes V')_f = \bigoplus_{f = f' \circ f''}V_{f'} \otimes
V'_{f''}.
$$
A {\em $\C$-graded $A_\infty$-coalgebra} is an $A_\infty$-coalgebra
in this category.

The derived category $\D(\B_\idot,\Ab)$ of $A_\infty$-functors from
a small $A_\infty$-ca\-te\-go\-ry $\B_\idot$ to the category of
complexes of objects in an abelian category $\Ab$ is defined by
considering the DG category of all such functors and $A_\infty$-maps
between them, and inverting quasiisomorphisms. The invertion
procedure presents no problems, since the DG category is
well-behaved (at least if $\Ab$ is large enough, for example
$\Ab=k\amod$ for some ring $k$). Every object has an $h$-projective
and an $h$-injective replacement. For every $A_\infty$-morphism
$f:\B_\idot \to \B'_\idot$, we have the pullback functor
$f^*:\D(\B'_\idot,\Ab) \to \D(\B_\idot,\Ab)$, and it has the right
and left-adjoint functors $f_*,f_!:\D(\B_\idot,\Ab) \to
\D(\B'_\idot,\Ab)$.

Given a $\C$-graded $A_\infty$-coalgebra $\R_\idot$, one considers
the DG category of $\C$-graded $k$-valued $A_\infty$-comodules over
$\R_\idot$. One can still invert quasiisomorphisms to obtain the
triangulated derived category $\D(\C,\R_\idot,k)$. For any
$A_\infty$-map $f:\R_\idot \to \R'_\idot$, we have the corestriction
functor $f^*:\D(\C,\R_\idot,k) \to \D(\C,\R'_\idot,k)$. For any
functor $f:\C' \to \C$, we have the pullback $\C'$-graded
$A_\infty$-coalgebra $f^*\R_\idot$ and a pullback functor
$f^*:\D(\C,\R_\idot,k) \to \D(\C',f^*\R_\idot,k)$. However,
$h$-projective replacements usually do not exist at all, and there
is no general procedure for constructing $h$-injective
ones. Therefore the existence of adjoints is non-trivial One case
where an adjoint does exist is the embedding functor $i:\ppt \to \C$
of the point category onto an object $c \in \C$. In this case, an
adjoint to the pullback functor $i^*:\D(\C,\T_\idot,k) \to \D(k)$ is
given by the cofree $A_\infty$-module functor, sending $M \in \D(k)$
to an $A_\infty$-comodule $M_c$ with
$$
M_c(c') \cong M \otimes \bigoplus_{f:c' \to c}\T_\idot(f).
$$
Another situation where things are easy is the {\em trivial
  $\C$-graded $A_\infty$-coalgebra} $\R_\idot$ given by $\R_0(f)
= \Z$ for any $f \in \C_1$, and $R_i = 0$ for $i \neq 0$.

\begin{lemma}\label{triv.coa}
Assume given a small category $\C$, and let $\R$ be the trivial
$\C$-graded $\A_\infty$-coalgebra. Then every $h$-projective complex
of functors from $\C^{opp}$ to $k\amod$ is $h$-projective as an
$A_\infty$-comodule over $\R$, and the natural embedding
$$
\D(\C^{opp},k) \to \D(\C,\R,k)
$$
is an equivalence of categories.
\end{lemma}

\proof{} Let $\B_\idot(c,c') = \Z[\C(c,c')]$ be the free additivie
category generated by $\C$, and treat it as an $A_\infty$-category
in the obvious way. Then by definition, an $A_\infty$-comodule over
the trivial $A_\infty$-coalgebra $\R$ is the same thing as an
$A_\infty$-functor from $\B_\idot^{opp}$ to complexes of
$k$-modules. Thus it suffices to prove the claim for
$A_\infty$-categories, wehre it is well-know.
\endproof

\begin{lemma}\label{cofree}
Assume given a small category $\C$, a $\C$-graded
$A_\infty$-coalgebra $\T_\idot$, and a functor $\rho:\C' \to \C$
from a small category $\C'$ such that $\rho^*\T_\idot$ is isomorphic
to the trivial $\C'$-graded $A_\infty$-coalgebra. Then the pullback
functor $\rho^*:\D(\C,\T_\idot,k) \to \D(\C',k)$ admits a
right-adjoint functor
$$
\rho_*:\D(\C',k) \to \D(\C,\T_\idot,k).
$$
\end{lemma}

\proof{} By definition, we have to prove that for any object $M \in
\D(\C',k)$, the functor
\begin{equation}\label{cor.adj}
N \mapsto \Hom(\rho^*N,M)
\end{equation}
from $\D(\C,\T_\idot,k)$ to the category of $k$-modules is
representable. By the cobar construction, every object $M \in
\D(\C',k)$ is the cone of an endomorphism of an object $M' \in
\D(\C',k)$ of the form
$$
M' = \prod M_i,
$$
where each $M_i = M'_{a}$ is the corepresentable functor from $\C'$
to $k\amod$ corresponding to an object $a \in \C'$ and some $M' \in
k\amod$. Thus it suffices to prove that the functor \eqref{cor.adj}
is representable for $M = M'_{a}$. The representing object is
given by the cofree $A_\infty$-comodule $M'_{\rho(a)} \in
\D(\C,\T_\idot,k)$.
\endproof

\begin{lemma}\label{cofree.bc}
Assume given a small category $\C$ equipped with a factorization
system, as in Lemma~\ref{fact.bc}. Moreover, assume given a
$\C$-graded $A_\infty$-coalgebra $\T_\idot$ such that $h^*\T_\idot$
is a trivial $\C_h$-graded $A_\infty$-coalgebra, as in
Lemma~\ref{cofree}. Then the base change map
$$
v^* \circ h_* \to \overline{h}_* \circ \overline{v}^*
$$
of functors from $\D(\C_h^{opp},k)$ to $\D(\C_v,v^*\T_\idot,k)$
induced by the obvious isomorphism $\overline{v}^* \circ h^* \cong
\overline{h}^* \circ v^*$ is itself an isomorphism.
\end{lemma}

\proof{} The same proof as for Lemma~\ref{fact.bc} works, with the
adjoints constructed by Lemma~\ref{cofree}.
\endproof

All the categories of $A_\infty$-comodules that we consider in this
paper ought to be symmetric tensor categories. However, to construct
the tensor product, one would need to equip the
$A_\infty$-coalgebras with some sort of Hopf algebra structure, and
this is too heavy technically. Therefore in general, we avoid tensor
products. We do need them in one easy case. Say that a $\C$-graded
$\A_\infty$-coalgebra $\R_\idot$ is {\em augmented} if $\R_i = 0$
for $i < 0$, and $\R_0(f) = \Z$ for any morphism $f$ in $\C$. We
them have an obvious augmentation map $\xi:\Z_{\C} \to \R_\idot$, and
the corresponding corestriction functor
$$
\xi^*:\D(\C^{opp},k) \to \D(\C,\R_\idot,k).
$$
However, we have more: for every complex $M_\idot$ in
$\Fun(\C^{opp},k)$ and every $A_\infty$-comodule $E_\idot$ over
$\R_\idot$, we can define the tensor product $E_\idot \otimes
\xi^*M_\idot$ by setting
$$
(E_\idot \otimes \xi^*M_\idot)(c) = E_\idot(c) \otimes M_\idot(c),
\qquad c \in \C,
$$
with the $A_\infty$-operations given by the products of the
$A_\infty$-operations in $E_\idot$ and the structure maps of the
functor $M_\idot$. This construction is obivously associative in
$M_\idot$.

\medskip

Finally, in order to construct $A_\infty$-categories and
$A_\infty$-coalgebras, we use various categorical constructions
associative ``up to an isomorphism''. Here is the prototype example
(for more details, see \cite[Subsection 1.6]{Ka-ma}). Assume given a
small monoidal category $\C$, with an associativity isomorphism
satisfying the usual pentagon equation. Then in effect, $\C$ is an
algebra over the following asymmetric operad.

\begin{defn}\label{I.n.ope}
The {\em monodial category operad} $I_n$ is an operad of groupoids
defined as follows:
\begin{enumerate}
\item on objects, $I_n$ is the free operad generated by a single
  binary operation,
\item on morphisms, there exists exactly one morphisms between any
  two objects of the groupoid $I_n$.
\end{enumerate}
\end{defn}

The bar complex $C_\idot(\C,k)$ is automatically an algebra over the
operad $C_\idot(I_\idot,\Z)$. But the operad $C_\idot(I_\idot,\Z)$
is a resolution of the associative operad $\Ass$, and the
$A_\infty$-operad $\Ass_\infty$ is another such resolution, and a
cofibrant one. Therefore the augmentation map $\Ass_\infty \to \Ass$
factors through a map $\Ass \to C_\idot(I_\idot,\Z)$. Fixing such a
factorization once and for all, we turn the bar complex
$C_\idot(\C,k)$ into an $A_\infty$-algebra over $k$.

\bigskip

\noindent
{\sc
Steklov Math Institute\\
Moscow, USSR}

\bigskip

\noindent
{\em E-mail address\/}: {\tt kaledin@mccme.ru}

\end{document}